\documentclass[12pt]{article}
\usepackage{amsmath}
\usepackage{mathrsfs} 
\usepackage{fixmath}
\usepackage{amsthm}
\usepackage{amsfonts}
\usepackage{slashed}
\usepackage{amssymb}
\usepackage{mathtools}
\usepackage{commath}
\usepackage{latexsym}
\usepackage{tikz} 
\usepackage{esint} 
\usepackage{bm}
\usepackage{cite}
\usepackage{tikz}
\usepackage{tikz-cd}

\usepackage[matrix,tips,graph,curve]{xy}



\makeatletter 
\@addtoreset{equation}{section}
\makeatother

\theoremstyle{plain}
\newtheorem{theorem}[equation]{Theorem}
\newtheorem{corollary}[equation]{Corollary}
\newtheorem{lemma}[equation]{Lemma}
\newtheorem{proposition}[equation]{Proposition}
\newtheorem{conjecture}[equation]{Conjecture}

\theoremstyle{definition}
\newtheorem{definition}[equation]{Definition}

\newtheorem{remark}[equation]{Remark}

\newcommand{\IC}{\mathbb{C}}

\newcommand{\IH}{\mathbb{H}}

\newcommand{\IL}{\mathbb{L}}

\newcommand{\IR}{\mathbb{R}}

\newcommand{\IZ}{\mathbb{Z}}

\newcommand{\coker}{\mathrm{coker}}

\newcommand\union{\bigcup} 

\def\d/{/\mspace{-6.0mu}/}

\newcommand{\e}{\mathbf{e}}

\let\oldphi\phi \let\phi\varphi \let\varphi\oldphi

\begin{document}
\title{Solutions of the Bogomolny Equation on $\IR^3$ with Certain Type of Knot Singularity I}
\author{Weifeng Sun}
\date{}

\maketitle

\tableofcontents

\section*{Abstract}

Moduli space of the Bogomolny equation on $\IR^3$ with certain boundary condition at infinity has been well studied for a long time. Only several examples are listed here for interested readers: \cite{Hitchin1983OnMonopoles} ,  \cite{Donaldson1984NahmsMonopoles} and \cite{Jaffe1980VorticesTheories}.  This paper studies the moduli space of solutions to the Bogomolny equation on $\IR^3$ with a knot singularity. The author hopes such kind of moduli spaces have potential applications in low-dimensional topology and knot theory in the future.\\

Fix a constant $\gamma \in (0, \dfrac{1}{2})$ and a smooth knot $K$  in $\IR^3$. There is a model solution to the Bogomolny equations with a knot singularity of monodromy $\gamma$ near the knot. Approaching infinity, there is also a (in some sence, standard) model solution with ``mass" $M > 0$ and ``charge" $k$. Any solution that is close enough to the model solution both near the knot and approaching infinity is what we study in this paper. The definition of ``close enough" is, the difference lies in a Hilbert space $\IH$ which will be defined. In fact, weaker assumptions imply that the difference is in $\IH$. And for a solution of the Bogomolny equations, we have better descriptions on this difference both near the knot and approaching infinity than merely lying in $\IH$.\\

Section 1-3 set up the Fredholm theory for the linearization of the Bogomolny equations with a knot singularity. We define Hilbert spaces $\IH$ and $\IL$. Then when $\gamma \neq \dfrac{1}{4}$ the linearization is a Fredholm operator from $\IH$ to $\IL$. Since the cokernel may not be $0$, this is not enough to give the moduli space a smooth structure, but only an analytical structure. The author suspects that the cokernel is actually $0$ but cannot prove it. On the other hand, the case that $\gamma = \dfrac{1}{4}$ is particular interesting and out of the scope of this paper. The author wishes to study it in the future.\\

Section 4 proves the existence of a solution with a knot singularity of monodromy $\gamma$ for each $\gamma \neq \dfrac{1}{4}$. The method is to glue several standard charge 1 monopoles (a solution with $k = 1$) onto the model solution with a knot singularity.  More precisely, assuming the ``mass" $M$ (or more explicitly, the limit of the norm of the Higgs field $|\Phi|$ as $|x| \rightarrow +\infty$) is large, one can carefully make the linearization of a glued ``approximate solution" to have a $0$ cokernel. And the an actual solution exists nearby by implicit function theorem. The idea of gluing monopoles originates in Taubes' thesis \cite{Jaffe1980VorticesTheories}, where he glued several charge 1 monpoles onto the trivial solution without a knot singularity.\\

Section \ref{Section:Behaviors of a solution when approaching infinity} and section \ref{Section: A Ulenbeck type of regularity theory near the knot}   give a regularity theorem of a solution to the Bogomolny equations. The main result is, if we only assume a smooth solution  to the Bogomolny equations on $\IR^3 \backslash K$ has a curvature with finite $\IL$ norm, then one can define a mass ``$M \geq 0$" and an integer charge ``k" for the solution. If $M > 0$, then after a gauge equivalent, the solution differs from the model solution by an element in $\IH$. This argument is divided into two aspects: Near the knot and approaching infinity. Section  \ref{Section:Behaviors of a solution when approaching infinity}  proves the aspect when approaching infinity, which is a modification of Taubes' arguments in \cite{Jaffe1980VorticesTheories}.  Section \ref{Section: A Ulenbeck type of regularity theory near the knot} proves the aspect near the knot, which is essentially a Ulenbeck type of arugment inspired by Sibner and Sibner's paper \cite{Sibner1992ClassificationHolonomy}.\\

Section $\ref{Section: Limit behaviors near the knot}$ studies the behavior of our solutions near the knot. It is conjectured that any solution has a poly-homogeneous expansion near the knot. This is not unexpected in the viewpoint of micro-local analysis. In fact, Mazzeo and Witten \cite{Mazzeo2013TheCondition}\cite{Mazzeo2017TheKnots} have proved the poly-homogeneous expansion of a solution to certain type of the extended Bogomolny equations in a different setting.  It is likely that their argument also works in our situation. But the author chooses to use a straightforward and different approach and get a weaker statement: only the leading term in the poly-homogeneous expansion.\\

\textbf{Remark on the second edition}~\\

More that two years after the author put the first edition of this paper on arxiv, the author decided to re-write it completely. The newer version differs from the older version in the following aspects:

\begin{itemize}
    \item Many notations are modified and many arguments are simplified. The author hopes the new edition is much easier to read.
    \item To establish a Fredholm theory, in the old version of this paper, we have assumed that the monodromy $\gamma \in (0, \dfrac{1}{8}) \cup (\dfrac{3}{8}, \dfrac{1}{2})$ due to technical difficulties. In the new version, the technical difficulty is resolved and the new assumption is $\gamma \neq \dfrac{1}{4}$. The case $\gamma = \dfrac{1}{4}$ is different in an essential way. The author hopes to address the case when $\gamma = \dfrac{1}{4}$ in a future paper.
    \item More details on the regularity theory are added, both near the knot and approaching infinity. 
    \item More studies on the asymptotic behaviors of a solution near the knot and approaching infinity are added.
    \end{itemize}

\textbf{Acknowledgements}  I would like to express my deep gratitude to professor C.Taubes for his patient guidance, enthusiastic encouragement and useful critiques of this research work. I also want to thank Boyu Zhang, Siqi He and Donghao Wang for many helpful discussions. For the second edition, I also owe professor R. Mazzeo a lot for introducing his theory about elliptic edge operators to me and kindly encouraging and supporting me all the time.\\

\section{Introduction and preliminary set-ups}\label{Section: Introduction and preliminary set-ups}
In this paper,  we use the standard convention to let $C$ represent a large constant which may change from line by line. 

\subsection{The Bogomolny equations and its linearization}\label{Subsection: The Bogomolny equations and its linearization}

Suppose the ambient space is a $3$-dimensional manifold. Suppose $A$ is an $su(2)$-valued $1$-form (representing an $SU(2)$ connection in a local coordinate). Suppose $\Psi$ is an $su(2)$-valued function (representing a local section of the adjoint $su(2)$ bundle). Suppose the curvature is $F_A = dA + \dfrac{1}{2}[A \wedge A]$. Then the Bogomolny equation is written as
$$V(A, \Phi) = *F_A - d_A \Phi = 0,$$
where $d_A = d + [A \wedge \cdot], ~~ *$ is the Hodge star opeartor.\\

Fix a fiducial configuration $\Psi = A + \Phi$ (which can be viewed as an $su(2)$-valued ``0 + 1" form). Then a variation $\psi = a + \phi$ can be written as

\begin{equation*}
  \begin{array}{ll}
V(\Psi + \psi) & = V(\Psi) +  *d_A a - d_A \phi - a \wedge \Phi  + \dfrac{1}{2}*[a \wedge a]  -[a \wedge \phi] \\
 & = V(\Psi) + L_{\Psi}(\psi)  + Q(\psi, \psi),
\end{array}
\end{equation*}
where $L_{\Psi}$ and $Q$ are the linear part and the quadratic part respectively. Note that $Q$ doesn't depend on the fiducial configuration $\Psi$.\\

The Bogomolny equations are $SU(2)$ gauge invariant. Typically, if two configurations are gauge equivalent, we regard them as the same. To make  $L_{\Psi}$ a candidate for a Fredholm operator, we need to add a gauge fixing condition
$$* d_A * a + [\Phi, \phi] =0.$$

\textbf{Remark:} This gauge fixing condition kills all the ``infinitesimal gauge transformations" but $1$-direction. The infinitesimal  guage transformations that cannot be detected by the gauge fixing condition is given by the Lie bracket with $\Phi$. We will discuss this extra 1-dimensional gauge transformation with more details in subsection \ref{Subsection: The extra 1 dimensional gauge transformation}.\\

The combination of $L_{\Psi}$ and the gauge fixing operator is written as $\tilde{L}_{\Psi}$. Then
$$\tilde{L}_{\Psi}(\psi) =  *d_A a - d_A \phi - [a \wedge \Phi]  + * d_A * a + [\Phi, \phi] = D_A + [\Phi, \cdot],$$ where $D_A \psi = D_A (a + \phi)= *d_Aa - d_A \phi + *d_A*a$ can be viewed as a twisted Dirac operator. The formal ad-joint of $ \tilde{L}_{\Psi}$ is
$$\tilde{L}_{\Psi}^{\dagger} = D_A - [\Phi, \cdot].$$

\textbf{Remark:} In this paper, the ``formal-adjoint" of an operator $A$ is the operator $A^{\dagger}$ which satisfies, for any two smooth and compactly supported functions/sections on $\IR^3 \backslash K$, namely $\psi$ and $\psi'$,
$$\int_{\IR^3 \backslash K} A(\psi)\psi' d^3 x = \int_{\IR^3 \backslash K}\psi A^{\dagger}(\psi') d^3x.$$

One goal of this paper is to find two suitable Hilbert spaces $\IH$ and $\IL$ such that:

\begin{itemize}
    \item $\tilde{L}$ is a Fredholm operator from $\IH$ to $\IL$.
    \item $Q$ is a bounded map from $\IH \times \IH$ to $\IL$.
    \item $\IH$ has some other good properties which will be specified later.
\end{itemize}

Before constructing the Hilbert spaces, we need to set up the ambient space and the model configuration first.

\subsection{The ambient space}\label{Subsection: The ambient space}

In this paper, the ambient space is $X = \IR^3 \backslash K$, where $K$ is a smooth knot embedded in $\IR^3$.  Any $SU(2)$-bundle over $X$ is (topologically) the product bundle. So in fact, we can use a ``0 + 1" form $A + \Phi$ to represent a configuration globally.\\

We need two sets of local parametrizations of $X$: one near the knot and one far away from the knot (``near infinity").

\paragraph{Parametrization near the knot}~\\

Throughout this paper, we assume $N_{\epsilon}$ is a small enough $\epsilon$ tubular neighbourhood of $K$. Let $N$ be the normal bundle of $K$. Choose an orthonormal frame of N, namely $\{\e_1(s), \e_2(s)\}$. Then $(s, z_1, z_2)$ gives a parametrization of $N_{\epsilon}$ by
$$(s, z_1, z_2) \rightarrow K(s) + z_1 \e_1(s) + z_2 \e_2(s) \in \IR^3.$$ 

We would also introduce the notations $z = z_1 + iz_2$ and $z = \rho e^{i\theta}$. We assume $\rho$ is the distance to the knot $K$. And $\rho$ is smooth with no critical points in $N_{\epsilon} \backslash K$. \\

We assume $|K'(s)| = 1$ and $s \in [0, l]$, where $l$ is the length of the knot. So $s(0) = s(l)$. Moreover, we may assume there are two smooth functions $k(s)$ and $\tau(s)$ such that

$$K''(s) = k(s) \textbf{e}_1(s), ~~ \textbf{e}_1'(s) = - k(s) K'(s) + \tau(s)\text{e}_2(s), ~~ \text{e}_2'(s) = -\tau(s) \text{e}_1(s). $$

Note that the Euclidean measure on $N_{\epsilon}$ is not the same as $\rho d\rho d\theta ds$. But they are equivalent. We will frequently use $\rho d\rho d\theta ds$ for convenience.

\paragraph{Parametrization ``near infinity"}~\\

We fix a large enough $R > 0$. Consider the ball $B_R$ of radius $R$ centered at the origin. We assume $K$ is in the interior of $B_R$. Then we use the coordinates $(r, S)$ for $\IR^3 \backslash B_R$ (the standard spherical coordinates), where $r$ is the distance to the origin, $S$ is the spherical coordinates for $\partial B_r$.\\

We use $dS$ to denote the measure on $\partial B_r$. Note that the Euclidean measure on $\IR^3 \backslash B_R$ is $r^2 dr dS$.

\subsection{The model configuration}\label{Subsection: The model configuration}

In this subsection, we set up a model configuration $\Psi_{\gamma, M, k}$, where all the three subscripts will be specified later. Roughly speaking, the model solution has a knot singularity of monodromy $\gamma$, and has ``mass" $M$ and charge $k$ near infinity. An interesting fiducial configuration (or a solution to the Bogomolny equations) should look like (will be specified later) the model configuration both near the knot and approaching infinity after a gauge transformation.

\paragraph{The model solution near the knot}~~\\

Suppose $\gamma \in (0, \dfrac{1}{2})$ is a fixed constant represent. Then here is a model solution to the Bogomolny equation with a knot singularity of monodromy $\gamma$:

\begin{definition}\label{definition of the model solution near the knot}
Fix a constant element $\sigma \in su(2)$ with norm $1$. Consider the differential $1$-form $d\theta$ in $N_{\epsilon}$. This $1$-form can be extended to the enire $\IR^3 \backslash K$ as a closed $1-$form with bounded support in $B_R$, still denoted as $d\theta$. Let $M$ be a fixed positive constant. Then the model configuration is $\Psi_{\gamma, M, 0} = A_{\gamma, M, 0} + \Phi_{\gamma, M, 0}$:
$$A_{\gamma, M, 0} = \gamma \sigma d\theta, ~~~ \Phi_{\gamma, M, 0} = M \sigma.$$

\end{definition}

It can be easily verified that the model solution satisfies the Bogomolny equations. The subscripts ``$M$" and "$0$" indicate that this model configuration has ``mass" $M$ and ``charge" $0$, which we will defined later.\\

Note that near the knot $|d\theta| = \dfrac{1}{\rho}$ is unbounded, illustrating that $A_{\gamma, M, 0}$ has a knot singularity.

\begin{remark}
A natural question is, how to define in general a configuration with the same type of singularity as the model one. The precise definition will be deferred later in subsection \ref{Subsection: Admissible configurations} after we set up the Hilbert spaces. But here are some intuitions: \\

Naively thinking, a configuration that looks like 
$$\Psi = \Psi_{\gamma, M, 0} + O(\rho^{-1 + \delta}) $$
should have the same type of singularity, where $O(\rho^{-1 + \delta})$ is a lower order term. But then both the curvature $F_A$ and $D_A \Phi$ would be $O(\rho^{-2 + 2\delta})$, that is, unbounded. However, if we use the following definition:

$$ \Psi = \Psi_{\gamma, M, 0} + \text{smaller terms} ~~~ \text{and} ~~~~ |F_A| ~\text{is bounded}, $$
then it would be too rigid to expect any solution except the model one. To compromise, we require the following integral
$$\int_{N_{\epsilon}} \rho^2 |F_A|^2 d\rho d\theta ds $$
to be finite. This would imply (under some other assumptions) that
$$\Psi = \Psi_{\gamma, M, 0} + O(\rho^{-\frac{1}{2} + \delta})$$
for some small $\delta > 0$.
Luckily, this is enough to do Fredholm theory for $\tilde{L}_{\Psi}$ (that is to say, roughly, the moduli space has finite dimension and is not an isolated point in general). Sadly, we have to accept the fact that $|F_A|$ may not be bounded near the knot.

\end{remark}

\paragraph{The model solution near infinity}~\\

It is well-known that for each $M > 0$ and non-negative integer $k$, there are smooth solutions to the Bogomolny equations on $\IR^3$ with ``mass" $M$ and ``charge" $k$. Suppose $\Psi$ is such a solution (may not be unique). Then it has the following properties:

\begin{itemize}
    \item  Fix a large enough $R > 0$. $$ \int_{\IR^3 \backslash B_R} ( |M - |\Phi||^6 + |\nabla_A \Phi|  + |F_A|^2)r^2 dr dS$$
is finite, where $B_R$ is the ball of radius $R$ centered at the origin. (This characterizes the mass ``M").
\item  Suppose $E$ is the product $su(2)$ bundle over $B_R$. Suppose $E = l \oplus l^{\perp}$, where $l$ is the $1$-dimensional sub line bundle generated by a unit vector field $e_0$, and  $l^{\perp}$ can be naturally viewed as a complex line bundle whose Chern class is $k$. Then on $\IR^3 \backslash B_R$,  under a suitable gauge, $\Psi$ can be written as $M e_0 + \psi$, where $\psi$ is a configuration such that the following integral is finite:
$$\int_{\IR^3 \backslash B_R} (|\nabla \psi|^2 + |[e_0, \psi]|^2) r^2 dr dS. $$
(This characterizes the charge ``k".)
\end{itemize}

It is convenient to specify a ``model solution" $\Psi_{M, k}$ on $\IR^3\backslash B_R$ as follows:
    $$A = [de_0, e_0], ~~ \Phi = (M - \dfrac{k}{r})e_0 ~~~ \text{when} ~~~k \geq 1, $$
where $$e_0 = (\cos(k\theta)\sigma_1 + \sin(k\theta) \sigma_2)\cos\phi + \sigma_3\sin\phi,$$ and $\sigma_1, \sigma_2, \sigma_3$ are a standard basis of $su(2)$, $(r, \theta, \phi)$ is the standard spherical coordinates of $\IR^3$ restricted on $\IR^3 \backslash B_{R}$. We may also assume $e_0 = \sigma$ when $k = 0$, where $\sigma$ is the constant element in $su(2)$ used in definition \ref{definition of the model solution near the knot}. One can check directly that the pair $(A, \Phi)$ satisfies the Bogomolny equations and all the properties that listed above.

\paragraph{The global model configuration}~\\

Throughout this paper, we fix a smooth cut-off function $\chi$ which is $1$ on $N_{\epsilon}$ and equals $0$ on $N_{2\epsilon}$. We also fix another smooth cut-off function $\chi'$ which is $1$ on $\IR^3 \backslash B_{2R}$ and $0$ on $B_R$. (Occasionally we may also write $\chi_R = 1 - \chi'$). Then
\begin{definition}
Here is the definition of the global model configuration $\Psi_{\gamma, M, k}$:
$$\Psi_{\gamma, M, k} = (1 - \chi') \Psi_{\gamma, M, 0} + \chi' \Psi_{M, k}.$$
\end{definition}

Note that when $k = 0$, it agrees with $\Psi_{\gamma, M, 0}$ defined in definition \ref{definition of the model solution near the knot}, which is a solution to the Bogomolny equations on $\IR^3 \backslash K$. However when $k \geq 1$, $\Psi_{\gamma, M, k}$ may not be a global solution. Neverthelss, 
$$\Psi_{\gamma, M, k} = \Psi_{\gamma, M, 0} ~ \text{on} ~ N_{\epsilon} \backslash K, ~~~ \Psi_{\gamma, M, k} = \Psi_{M, k} ~ \text{on} ~ \IR^3 \backslash B_{2R}.$$

So $\Psi_{\gamma, M, k}$ satisfies the Bogomolny equations on $( \IR^3 \backslash B_{2R}) \cup N_{\epsilon}$.

\subsection{The Banach spaces $\IH$ and $\IL$} \label{Subsection: The Banach spaces H and L}

Now it is time to define the Banach spaces $\IH$ and $\IL$ that we use in this paper. Here are some candidates:

\paragraph{Spaces of Sobolev types near the knot}~\\

The operator $\tilde{L}_{\Psi}$ is an elliptic operator of ``edge type" near the knot, which has been well-studied in microlocal analysis. One possible reference is Mazzeo's work \cite{Mazzeo1991EllipticI}. To make $\tilde{L}_{\Psi}$ a Fredholm operator, it is standard to add ``weights" to the ordinary Sobolev spaces. To be more precise, for $\IL$, we use the following definition near the knot

$$\int_{N_{\epsilon}\backslash K} \rho^{2\delta} |\psi|^2 \rho d\rho d\theta ds$$

and for $\IH$, we use the following definition
$$\int_{N_{\epsilon}\backslash K} \rho^{2\delta} (|\partial_{\rho} \psi|^2 + |\partial_s \psi|^2 + \dfrac{1}{\rho^2}(|\partial_{\theta}\psi|^2 + |\psi|^2)) \rho d\rho d\theta ds,$$
where $\delta$ is some constant.\\

Then it terms out that, $\tilde{L}$ is Fredholm (near the knot) when $$\min\{2\gamma, 1 - 2\gamma\} < \delta < \max\{2\gamma, 1 - 2\gamma\}.$$ In this paper, we simply choose $\delta = \dfrac{1}{2}$ for $\gamma \neq \dfrac{1}{4}$. The $\gamma = \dfrac{1}{4}$ case is especially interesting and tricky, which the author hopes to address in a future paper. We have two remarks:

\begin{enumerate}
    \item This definition is equivalent to Mazzeo's weighted Sobolev spaces for elliptic edge operators, see \cite{Mazzeo1991EllipticI}. In fact, $\delta$ actually lies in the region corresponding to ``Fredholm weights" for $\tilde{L}_{\Psi}$. So this Fredholm theory near the knot is not unexpected. But instead of simply quoting Mazzeo's theory, we give a direct proof in this paper when $\delta = \dfrac{1}{2}$ to make it self-contained.
    \item If $\delta = \dfrac{1}{2}$, then locally $Q$ is a bounded map from $\IH \times \IH$ to $\IL$. This is also not unexpected: It is the ``weighted spaces" version of a standard Sobolev multiplication theorem.  The precise statement and the  proof is in subsection \ref{Subsection: The quadratic operator Q}.
    \item In fact, $|\partial_{\rho}\psi|^2 + |\partial_s\psi|^2 + \dfrac{1}{\rho^2}|\partial_{\theta}\psi|^2$ is equivalent with $|\nabla \psi|^2$ near the knot. So finally our $\IH$ norm near the knot can also be defined as
    $$\int_{N_{\epsilon} \backslash K}(\rho^2 |\nabla \psi|^2 + |\psi|^2) d\rho d\theta ds.$$
\end{enumerate}

\paragraph{Spaces of Sobolev types near infinity}~\\

There are two methods to define $\IH$ and $\IL$ near infinity. The first method, which is more standard, is to add a weight near infinity. But in this paper, we use another method which is essentially due to Taubes in his thesis \cite{Jaffe1980VorticesTheories} (but with some modifications). And here it is:\\

Near infinity, we can use the following definition for $\IH$:

$$\int_{\IR^3 \backslash B_R} (|\nabla_A \psi|^2 + |[e_0, \psi]|^2) r^2 dr dS.  $$

\textbf{Remark:}~
\begin{enumerate}
    \item The above definition of $\IH$ doesn't bound $L^2$ near infinity. (In fact, if one decomposes $\psi$ into two components by parallel or orthogonal to $\Phi$ near infinity, then $\IH$ only bounds the $L^2$ norm of the orghogonal component of $\psi$.)
    \item Strictly speaking, the definition of $\IH$ depends on $e_0$, which depends on $k$. So we may also write it as $\IH_k$ to emphasize this dependence occasionally.
\end{enumerate}

And we use the ordinary $L^2$ norm for $\IL$ near infinity.

\paragraph{Spaces of Holder type}~\\

One may alternatively uses spaces of Holder type (with suitable weights near the knot and near infinity), following the spirit of \cite{Mazzeo1991EllipticI}. The potential advantage of these spaces is:

\begin{itemize}
    \item The quadratic map $Q$ sends $\IH \times \IH$ to $\IL$ for granted. (There is no need to carefully establish a weighted version of ``Sobolev multiplication inequality".)
\end{itemize}

But we have the following disadvantages:

\begin{itemize}
    \item It is hard (at least for the author) to establish a Ulenbeck type of argument like section \ref{Section: A Ulenbeck type of regularity theory near the knot} in this paper.
    \item In order to prove the existence of an actual solution near a glued approximate solution, either the cokernel of $\tilde{L}_{\Psi}$ is $0$, or more generally, certain obstruction vanishes. But this is a hard task (at least for the author) in the realm of Holder spaces.
\end{itemize}

In this paper, we choose not to use the spaces of the Holder types.

\paragraph{Conclusion}

After all, here are the definitions that we choose for $\IH$ and $\IL$ in this paper:\\

The $\IH$ norm of a configuration $\psi$ is the square root of 

$$\int_{N_{\epsilon} \backslash K} (\rho^2 |\nabla \psi| + |\psi|^2 )d\rho ds d\theta$$ $$ + \int_{B_R \backslash N_{\epsilon}} (|\nabla \psi|^2 + |\psi|^2) d^3x + \int_{\IR^3 \backslash B_R} (|\nabla \psi|^2 + |[e_0, \psi]|^2) r^2 dr dS. $$

And $\IL$ norm of a configuration $\psi$ is the square root of
$$\int_{N_{\epsilon} \backslash K} \rho^2 |\psi|^2 d\rho ds d\theta + \int_{\IR^3 \backslash N_{\epsilon}} |\psi|^2 d^3 x. $$

Strictly speaking, $\IH$ depends on $e_0$, which further depends on the charge ``$k$". But typically we don't emphasize on this dependence.\\

The Hilbert spaces $\IH$ and $\IL$ are defined to be the completions of smooth configurations with bounded support in $\IR^3 \backslash K$ under the corresponding norms.\\

Occasionally, we will use an equivalent way to represent $|\nabla \psi|^2$ as $$|\partial_{\rho}\psi|^2 + |\partial_s\psi|^2 + \dfrac{1}{\rho^2}|\partial_{\theta}\psi|^2$$ near the knot. (They are not equal in general. But they bound each other and define equivalent norms.) We will sometimes freely make this substitution implicitly in our arguments if there is no ambiguity.

\subsection{Admissible configurations}\label{Subsection: Admissible configurations}

Fix $\gamma \in (0, \dfrac{1}{2})$, $M > 0$ and a non-negative integer $k$.

\begin{definition}
If $\Psi$ is a configuration on $\IR^3 \backslash K$ such that 
$$\Psi - \Psi_{\gamma, M, k} \in \IH,$$
then we say $\Psi$ is an admissible configuration with a knot singularity of monodromy $\gamma$, mass $M$ and charge $k$, or an admissible configuration for simplicity.\\

If $\Psi$ is gauge equivalent to an admissible configuarion, then we also say $\Psi$ has a knot singularity of monodromy $\gamma$, mass $M$ and charge $k$. 
\end{definition}

Recall that $\chi$ is the cut-off function that is $1$ on $N_{\epsilon}$ and $0$ on $\IR^3 \backslash N_{2\epsilon}$ defined in subsection \ref{Subsection: The model configuration}. 

\begin{definition}
To say that a configuration $\Psi$ has a knot singularity with monodromy $\gamma$  means after a gauge transformation, $\Psi$ satisfies:
$$\chi (\Psi - \Psi_{\gamma, M, 0}) \in \IH.$$
\end{definition}

Despite the fact that $\Psi_{\gamma, M, 0}$ depends on $M$ and $\IH$ depends on $k$, the above definition doesn't depend on $M$ or $k$ because $\chi M \sigma$ itself is in $\IH$ and because this definition only depends on $\Psi$ on $N_{2\epsilon} \backslash K$.\\

In this paper, we are mainly interested in admissible configurations. \\

In fact, $\Psi = A + \Phi$ is a smooth solution to the Bogomolny equations on $\IR^3 \backslash K$ such that  $F_A \in \IL$, then one of the three statements is true:

\begin{itemize}
    \item $\lim\limits_{|x| \rightarrow +\infty} |\Phi| = 0$. This is the situation that is not studied in this paper. 
    \item $\Psi$ can be extended as a global smooth solution to the Bogomonly equatuions on $\IR^3$.
    \item $\Psi$ is gauge equivalent to an admissible configuration for some $\gamma \in (0, \dfrac{1}{2})$, $M > 0$ and non-negative integer $k$.
\end{itemize}

So solutions are essentially admissible under very weak assumptions.

\subsection{Some analytical lemmas}\label{Subsection: Some analytical lemmas}

This subsection is a digression from our main topic. We establish some analytical lemmas that will be used in this paper.

\begin{lemma}\label{boundary convergent lemma}
Suppose $\psi$ is smooth on $N_{\epsilon}\backslash K$ with finite $\IH$ norm. Then
$$ \lim\limits_{\rho \rightarrow 0} \int_0^l\int_0^{2\pi} \rho | \psi|^2 d\theta ds = 0.$$
\end{lemma}

\begin{proof}
This lemma is actually a special case of proposition \ref{proposition of limit behavior of H and L}. But here is a simple proof:
Fix $\theta, s$, then
\begin{equation*}
    \begin{array}{ll}
         |\psi(2\rho) -  \psi(\rho)| & \leq \displaystyle \int_{\rho}^{2\rho} |\partial_{t} \psi(t)| dt  \\
         & \leq (\rho)^{\frac{1}{2}} (\displaystyle \int_{\rho}^{2\rho} |\partial_t \psi(t)|^2)^{\frac{1}{2}} \\
         & \leq \rho^{-\frac{1}{2}} \displaystyle (\int_{\rho}^{2\rho} t^2 |\partial_t \psi(t)|^2)^{\frac{1}{2}}.
    \end{array}
\end{equation*}
So $$\rho^{\frac{1}{2}}|\psi(\rho)| \leq \dfrac{(2\rho)^{\frac{1}{2}}|\psi(2\rho)|}{\sqrt{2}} +  \displaystyle (\int_{\rho}^{2\rho} t^2|\partial_t \psi(t)|^2 dt)^{\frac{1}{2}}. $$

Using the following inequaily: $(a+b)^2 \leq \dfrac{3}{2}a^2 + 9b^2,$
$$\rho|\psi(\rho)|^2 \leq \dfrac{3}{2}(\dfrac{(2\rho)^{\frac{1}{2}}|\psi(2\rho)|}{\sqrt{2}})^2 +  9(\displaystyle (\int_{\rho}^{2\rho} t^2|\partial_t \psi(t)|^2 dt)^{\frac{1}{2}})^2 $$
$$=\dfrac{3}{4}(2\rho|\psi(2\rho)|^2) +  9\displaystyle (\int_{\rho}^{2\rho} t^2|\partial_t \psi(t)|^2 dt). $$
Hence
$$\int_0^l \int_0^{2\pi}  \rho|\psi(\rho)|^2 d\theta ds \leq \dfrac{3}{4}\int_0^l \int_0^{2\pi}  (2\rho)| \psi(2\rho)|^2d\theta ds + 9 \int_0^l \int_0^{2\pi}  \int_{\rho}^{2\rho} t^2|\partial_t \psi(t)|^2 dt d\theta ds. $$
$$\int_{\partial N_{\rho}} |\psi|^2 \rho d\theta ds \leq \dfrac{3}{4}\int_{\partial N_{2\rho}} |\psi|^2 \rho d\theta ds + 9 \int_0^l \int_0^{2\pi}  \int_{\rho}^{2\rho} t^2|\partial_t \psi(t)|^2 dt d\theta ds. $$
Since $\psi \in \IH$,  $$ \int_0^l \int_0^{2\pi} \int_0^{\epsilon} t^2 |\partial_t \psi(t)|^2 dt d\theta ds < \infty,$$
then $$\lim\limits_{\rho \rightarrow 0} 9 \int_0^{2\pi} \int_0^{\epsilon} \int_{\rho}^{2\rho} t^2|\partial_t \psi(t)|^2 dt d\theta ds = 0, $$ which implies
$$\lim\limits_{\rho \rightarrow 0} \int_0^l \int_0^{2\pi}  \rho|\psi(\rho)|^2 d\theta ds = 0. $$
\end{proof}

\begin{lemma}\label{lemma of a preliminary inequality on N_epsilon}
Suppose $\psi \in \IH$ is a smooth configuration. For any $\eta \in (0,1)$,
\begin{equation*}
         \displaystyle   \eta (1-\eta) \int_0^l \int_0^{2\pi}\int_0^{\epsilon} |\psi|^2 d\rho d\theta ds  \leq \int_0^l \int_0^{2\pi}\int_0^{\epsilon} \rho^2 |\partial_{\rho} \psi|^2 d\rho d\theta ds + \eta \int_0^l \int_0^{2\pi} \epsilon |\psi(\epsilon)|^2 d\theta ds,
\end{equation*}
\begin{equation*}
         \displaystyle \eta \int_0^l \int_0^{2\pi} \epsilon |\psi(\epsilon)|^2 d\theta ds
          \leq \displaystyle  \int_0^l \int_0^{2\pi} \int_{0}^{\epsilon} \rho^2|\partial_{\rho} \psi|^2 d\rho d\theta ds + \eta (1 + \eta) \int_0^l \int_0^{2\pi} \int_{0}^{\epsilon} |\psi|^2 d\rho d\theta ds.
\end{equation*}
\end{lemma}

\begin{proof}

Using the previous lemma
\begin{equation*}
    \begin{array}{ll}
         \displaystyle \int_0^l \int_0^{2\pi} \epsilon |\psi(\epsilon)|^2 d\theta ds & = \displaystyle \int_0^l \int_0^{2\pi} \int_{0}^{\epsilon} \partial_{\rho}(\rho  |\psi(\rho)|^2) d\rho d\theta ds \\
         & = \displaystyle \int_0^l \int_0^{2\pi} \int_{0}^{\epsilon} \rho \partial_{\rho}(  |\psi(\rho)|^2) + |\psi(\rho)|^2 d\rho d\theta ds \\
         & \geq \displaystyle  \int_0^l \int_0^{2\pi} \int_{0}^{\epsilon} ((1 - \eta)|\psi(\rho)|^2 - \dfrac{1}{\eta} \rho^2|\partial_{\rho} \psi|^2   )d\rho d\theta ds, \\
    \end{array}
\end{equation*}
and
\begin{equation*}
    \begin{array}{ll}
         \displaystyle \int_0^l \int_0^{2\pi} \epsilon |\psi(\epsilon)|^2 d\theta ds
         & = \displaystyle \int_0^l \int_0^{2\pi} \int_{0}^{\epsilon} \rho \partial_{\rho}(  |\psi(\rho)|^2) + |\psi(\rho)|^2 d\rho d\theta ds \\
         & \leq \displaystyle \int_0^l \int_0^{2\pi} \int_{0}^{\epsilon} (\dfrac{1}{\eta}\rho^2|\partial_{\rho} \psi|^2 + (1 + \eta)|\psi(\rho)|^2) d\rho d\theta ds.
    \end{array}
\end{equation*}
\end{proof}

\begin{lemma}\label{Appendix1}
Let $\Psi = A + \Phi$ be a fixed fiducial configuration. Suppose $\psi$ is a smooth configuration on $\IR^3 \backslash K$ with bounded support. Suppose $B_R$ is the ball of radius $R$ centered at origin. Suppose $R$ is large enough. Recall that $d^3x = r^2 drdS$ on $\IR^3 \backslash B_R$. Then
$$2\int_{\IR^3 \backslash B_R} |\nabla_A \psi|^2 d^3x\geq \displaystyle  \dfrac{1}{R} \int_{\partial B_R} |\psi(R)|^2 R^2dS + \dfrac{1}{2}\int_{\IR^3 \backslash B_R} \dfrac{1}{r^2}|\psi|^2 d^3x. $$
\end{lemma}
\begin{proof} Suppose $\nabla_{A_r}$ is the covariant derivative $\nabla_A$ in $\partial_r$ direction. Then
\begin{equation*}
    \begin{array}{ll}
         ~ & ~~~~\displaystyle 2\int_{\IR^3 \backslash B_R} |\nabla_A \psi|^2 r^2 dr dS + \dfrac{1}{2}\int_{\IR^3 \backslash B_R} |\psi|^2 dr dS  \geq \displaystyle - \int_R^{+\infty} \int_{\partial B_r} 2 r |\nabla_{A_r} \psi||\psi| dS dr \\ & \geq \displaystyle - \int_R^{+\infty} \int_{\partial B_r} r \partial_r (|\psi|^2) dS dr = \displaystyle  \dfrac{1}{R} \int_{\partial B_R} |\psi(R)|^2 R^2dS + \int_{\IR^3 \backslash B_R} |\psi|^2 drdS.
    \end{array}
\end{equation*}
\end{proof}

\begin{theorem}\label{Appendix2}
Same condition as the lemma \ref{Appendix1}. Suppose $\epsilon \leq \delta$ are two small enough constants. There exists a constant $C$ which depends on $\epsilon$ and the knot $K$, but not on $A$, such that $$\int_{\partial N_{\epsilon}} |\psi|^2 \rho d\theta ds + \int_{B_R\backslash N_{\epsilon}} |\psi|^2 d^3x \leq C \int_{\IR^3 \backslash N_{\epsilon}} |\nabla_A \psi|^2 d^3x. $$
\end{theorem}

\begin{proof}
From lemma \ref{Appendix1}, $$\displaystyle \int_{\partial B_R} |\psi|^2 R^2dS$$  can be bounded  (up to a constant) by $$\displaystyle \int_{\IR^3 \backslash N_{\epsilon}} |\nabla_A \psi|^2 d^3x.$$ By a standard trace embedding theorem, $$\displaystyle \int_{\partial N_{\epsilon}} |\psi|^2 \rho d\theta ds + \int_{\partial B_R} |\psi|^2 R^2dS$$ can be bounded (up to a constant) by $$\displaystyle \int_{B_R\backslash N_{\epsilon}} (|\nabla (|\psi|^2)| + |\psi|^2)d^3x \leq \displaystyle \int_{\IR^3\backslash N_{\epsilon}} |\nabla_A \psi|^2 d^3x + 5 \displaystyle \int_{B_R\backslash N_{\epsilon}}  |\psi|^2 d^3x.$$ So it remains to show $ \displaystyle \int_{B_R\backslash N_{\epsilon}}  |\psi|^2 d^3x$  can be bounded (up to a constant) by  $ \displaystyle \int_{\IR^3\backslash N_{\epsilon}} |\nabla_A \psi|^2 d^3x + \displaystyle \int_{\partial B_R} |\psi|^2 R^2dS$.\\

Suppose on the contrary, there exists a sequence $\{ \psi_n\}$ such that
$$\lim\limits_{n \rightarrow +\infty} \displaystyle  (\int_{\IR^3\backslash N_{\epsilon}} |\nabla_A \psi_n|^2 d^3x + \displaystyle \int_{\partial B_R} |\psi_n|^2 R^2dS) = 0, ~~ \text{but} ~ \displaystyle \int_{B_R\backslash N_{\epsilon}}  |\psi_n|^2 d^3x = \text{Vol} (B_R \backslash N_{\epsilon}) > 0.$$
Then 
$$ \limsup\limits_{n \rightarrow +\infty} \displaystyle  \int_{B_R \backslash N_{\epsilon}} |\nabla (|\psi_n|^2)| d^3x  \leq 2 \limsup\limits_{n \rightarrow +\infty} ~ ((\int_{B_R \backslash N_{\epsilon}} |\nabla_A \psi_n|^2 d^3x)^{\frac{1}{2}} (\int_{B_R \backslash N_{\epsilon}} |\psi_n|^2| d^3x)^{\frac{1}{2}}) = 0. $$

So $|\psi_n|^2$ have bounded $L^1_1(B_R \backslash N_{\epsilon})$ norms. By a Rellich lemma, there exists a subsequence of $\{|\psi_n|^2\}$, namely $\{|\psi_{n_k}|^2\}$, that converges to some $f \in L^1_1(B_R \backslash N_{\epsilon})$. \\

Hence $\displaystyle  \int_{B_R \backslash N_{\epsilon}} |\nabla f| d^3x = 0,$ and $\displaystyle \int_{B_R\backslash N_{\epsilon}}  f d^3x = \text{Vol} (B_R \backslash N_{\epsilon})$, so $f $ is the constant function $1$.\\

By a trace embedding thereom,
$$ \lim\limits_{k \rightarrow + \infty} \int_{\partial B_R} ||\psi_{n_k}|^2 - 1| R^2dS = 0, $$
which contradicts the assumption
$$ \lim\limits_{n \rightarrow + \infty} \int_{\partial B_R} |\psi_n|^2 R^2dS = 0. $$
\end{proof}

\begin{corollary}\label{Appendix3}
Same condition as lemma \ref{Appendix1}. There exists a constant $C$ that doesn't depend on $\epsilon$, such that, supposing $\epsilon \leq \delta$, and supposing $\eta > 0$, then
$$ \displaystyle \epsilon^{\eta-1}\int_{\partial N_{\epsilon}} |\psi|^2 \rho d\theta ds  \leq (C\delta^{\eta-1} + \dfrac{\delta^{\eta }}{\eta})\displaystyle  \int_{\IR^3 \backslash N_{\epsilon}}|\nabla_A \psi|^2 d^3x.$$
\end{corollary}

\begin{proof}
Suppoce $C$ is the constant from theorem \ref{Appendix2} such that
$$\int_{\partial N_{\delta}}|\psi|^2 \rho d\theta ds \leq C \int_{\IR^3 \backslash N_{\delta}}|\nabla_A \psi|^2 d^3x. $$
Then for any $\epsilon \leq \delta$,
\begin{equation*}
    \begin{array}{ll}
         ~ & ~~~~ \displaystyle \epsilon^{\eta-1}\int_{\partial N_{\epsilon}} |\psi|^2 \rho d\theta ds =   \delta^{\eta-1} \int_{\partial N_{\delta}}|\psi|^2 \rho d\theta ds -   \int_{\epsilon}^{\delta} \int_{\partial N_{\rho}} \partial_{\rho}(\rho^{\eta}|\psi|^2) d\rho d\theta ds\\
         & \leq \displaystyle C \delta^{\eta-1} \int_{\IR^3 \backslash N_{\delta}}|\nabla_A \psi|^2 d^3x - (\eta + 1)\int_{\epsilon}^{\delta} \int_{\partial N_{\rho}}\rho^{\eta-1}|\psi|^2 d\rho d\theta ds +  2\int_{\epsilon}^{\delta} \int_{\partial N_{\rho}}\rho^{\eta}|\nabla_A \psi||\psi|  d\rho d\theta ds \\
         & \leq \displaystyle C \delta^{\eta-1} \int_{\IR^3 \backslash N_{\delta}}|\nabla_A \psi|^2 d^3x +  \dfrac{1}{(\eta + 1)}\int_{\epsilon}^{\delta} \int_{\partial N_{\rho}}\rho^{\eta }|\nabla_A \psi|^2 \rho d\rho d\theta ds  \\
         & \leq (C\delta^{\eta-1} + \dfrac{\delta^{\eta }}{\eta})\displaystyle  \int_{\IR^3 \backslash N_{\delta}}|\nabla_A \psi|^2 d^3x.
    \end{array}
\end{equation*}
\end{proof}

\begin{theorem}\label{A Rellich theorem}(A Rellich theorem)
 Suppose $S$ is a compact subset of $\IR^3 \backslash K$, and suppose $\psi \in \IH$ for some fixed smooth configuration $\Psi$.\\

(1)The following integral is compact relative to $\norm{\psi}^2_{\IH}$: ~~$\displaystyle \int_{S} |\psi|^2 d^3x,$\\
which means, any sequence in $\IH$ with bounded $\IH$ norms has a sub-sequence that the above integral converges.\\

(2)For any $\epsilon > 0$, the following integral is also compact relative to $\norm{\psi}_{\IH}$: ~~ $\displaystyle\int_{\partial N_{\epsilon}} |\psi|^2 \rho d\theta ds.$\\
\end{theorem}

\begin{proof}
(1) Using theorem \ref{Appendix2} by assuming $S \subset B_R\backslash N_{\epsilon}$ for some $R$ and $\epsilon$, it is clear that the  norm $L^{1,2}(S)$ of $|\psi|$ on $S$ is bounded above by $\norm{\psi}_{\IH}$. Hence (1) follows by a Rellich theorem on $S$.\\

(2) Letting $S$ be the completion of $N_{\epsilon} \backslash N_{\frac{\epsilon}{2}}$. Then the norm $L^{1,2}(S)$ of $\psi$ on S is bounded above by $\norm{\psi}_{\IH}$. Since $\partial N_{\epsilon}$ is one component of $\partial S$, (2) follows from a compact trace embedding theorem. 
\end{proof}

\begin{lemma}\label{weight version of Sobolev embedding lemma}
There exists a constant $c_1 > 0$ which doesn't depend on $\epsilon$ such that, supposing $\psi$ is a smooth configuration on $\IR^3  \backslash K$ with finite $\IH$ norm, then  \begin{equation*}
    \begin{array}{ll}
         \displaystyle 
         (\int_{N_{\epsilon} \backslash K} \rho^4 |\psi|^6  d\rho d\theta ds)^{\frac{1}{3}} \leq c_1 \int_{N_{\epsilon}\backslash K} (\rho^2 |\nabla \psi|^2  +  |\psi|^2) d\rho d\theta ds.
    \end{array}
\end{equation*}
An immediate corollary is, for any $1 \leq p \leq 3$ and a possibly slightly larger $c_1$
\begin{equation*}
    \begin{array}{ll}
         \displaystyle 
         (\int_{N_{\epsilon} \backslash K} \rho^{-2} (\rho^2|\psi|^2)^p  d\rho d\theta ds)^{\frac{1}{p}} \leq c_1 \int_{N_{\epsilon}\backslash K} (\rho^2 |\nabla \psi|^2  +  |\psi|^2) d\rho d\theta ds.
    \end{array}
\end{equation*}
\end{lemma}

\begin{proof}
In fact, this lemma is just a special case of the first bullet of proposition proposition \ref{proposition of generilized weighted Sobolev embedding} with $\delta = -\dfrac{1}{2}$. But we give a straightforward and simple proof here.\\

Consider the function $\rho \psi$. Clearly the $L^{1, 2}$ norm of $\rho^{\frac{1}{2}} \psi$
$$(\int_{N_{\epsilon} \backslash K} \rho (|\nabla(\rho^{\frac{1}{2}}\psi)|^2 + |\rho^{\frac{1}{2}}\psi|^2 ) d\rho d\theta ds)^{\frac{1}{2}}$$
is bounded above by a constant times
$$(\int_{N_{\epsilon} \backslash K} (\rho^2 |\nabla \psi|^2 + |\psi|^2) d\rho d\theta ds)^{\frac{1}{2}}. $$
So by the standard Sobolev embedding theorem, the $L^6$ norm of $\rho^{\frac{1}{2}}\psi$ is bounded above by the same amount up to a multiplication of a constant. But hte $L^6$ norm of $\rho^{\frac{1}{2}}\psi$ is just
$$(\int_{N_{\epsilon} \backslash K} \rho|\rho^{\frac{1}{2}}\psi|^6 d\rho d\theta ds)^{\frac{1}{6}} = (\int_{N_{\epsilon} \backslash K} \rho^4|\psi|^6 d\rho d\theta ds)^{\frac{1}{6}}. $$
 \end{proof}

\begin{theorem}\label{Appendix7}
(Sobolev embedding on $\IR^3 \backslash N_{\epsilon}$) Suppose $\psi$ is a smooth configuration on $\IR^3 \backslash N_{\epsilon}$ with bounded support. Let $A$ be a fiducial connection. Then there exists a constant $C > 0$ which doesn't depend on $A$ and $\psi$ such that 
$$\int_{\IR^3 \backslash N_{\epsilon}} |\psi|^6 d^3 x \leq C \int_{\IR^3 \backslash N_{\epsilon}} |\nabla_A \psi|^2 d^3x. $$
\end{theorem}

\begin{proof}
Note that $\nabla(|\psi|^2) \leq |\nabla_A \psi|^2$, then the theorem is just a standard 3-dimensional Sobolev inequality.
\end{proof}

\begin{theorem}\label{Appendix8}(A Hardy's inequality) Suppose $\eta \in (0,1)$. Then there exists a constant $C_{\eta}$. For any ball $B_R \subset \IR^3$ of radius $R$ for some $R > 0$, let $B_{\eta R}$ be the ball with the same center as $B_R$. Suppose $f$ is a smooth function on $B_{R}$. Let $r$ be the distance to the center of the ball. Then
$$\int_{B_{\eta R}} \dfrac{1}{r}|f|^2 d^3x \leq C_{\eta}(\dfrac{1}{R}\int_{B_{R}} |f|^2 d^3x + R \int_{B_R} |\nabla f|^2 d^3x). $$
\textbf{Remark:} A direct corollary of the above theorem is $\displaystyle \int_{B_{R}} \dfrac{1}{r}|f|^2 d^3x$ can be bounded above (up to a constant) by  $\displaystyle C_{\eta}(\dfrac{1}{R}\int_{B_{R}} |f|^2 d^3x + R \int_{B_R} |\nabla f|^2 d^3x).$\\

Note that since the Green's function to the Laplacian operator in 3-dimensional space has a point singularity of type $\sim \dfrac{1}{r}$, the above inequality can be used to give a point-wise estimate of a function on a bounded open set in $\IR^3$ based on its Laplacian and certain boundary conditions.
\end{theorem}

\begin{proof}
Let $d^3x = r^2 dS$, where $dS$ be the spherical volume form. Suppose $\chi(r)$ is a cut-off function which equals $1$ when $0 \leq r \leq \eta R$ and $0$ when $r = R$. It can be assumed that $|\chi'(r)| \sim \dfrac{1}{R}$ up to a constant which depends on $\eta$ but independent with $R$. Then
$$ \displaystyle \int_{B_{\eta R}} \dfrac{1}{r}|f|^2 d^3 x =  \displaystyle \int_{\partial B_{r}}\int_{0}^{\eta R}  |f|^2 rdrdS \leq \displaystyle \int_{\partial B_r} \int_0^{R}  \chi(r)^2 |f(r)|^2 rdr dS, $$
and 
\begin{equation*}
    \begin{array}{ll}
         ~  & ~~~~ \displaystyle  \int_{\partial B_r} \int_0^{R}  \chi(r)^2 |f(r)|^2 rdr dS\\
         & = \displaystyle - \int_{\partial B_r} \int_0^R  \dfrac{d}{dr}(\chi(r) |f(r)|)\chi(r) |f(r)| r^2 dr dS \\
         & \leq \displaystyle \dfrac{1}{2}\int_{\partial B_r} \int_0^{R} \chi(r)^2 |f(r)|^2 r dr dS + \dfrac{1}{2}\int_{\partial B_r} \int_0^R  |\dfrac{d}{dr}(\chi(r) |f(r)|)|^2 r^3 dr dS\\
         & = \displaystyle \dfrac{1}{2}\int_{\partial B_r} \int_0^{R} \chi(r)^2 |f(r)|^2 r dr dS + \int_{\partial B_r} \int_0^R  (|\chi'(r) f(r)|^2 + |\chi(r) |f'(r)|^2) r^3 dr dS\\
         & \leq \displaystyle \dfrac{1}{2}\int_{\partial B_r} \int_0^{R} \chi(r)^2 |f(r)|^2 r dr dS + R\int_{\partial B_r} \int_0^R  (|\chi'(r) f(r)|^2 + |\chi(r) |f'(r)|^2) r^2 dr dS\\
         & \leq \displaystyle \dfrac{1}{2}\int_{\partial B_r} \int_0^{R} \chi(r)^2 |f(r)|^2 r dr dS + \dfrac{C_{\eta}}{2} (\dfrac{1}{R}\int_{B_{R}} |f|^2 d^3x + R \int_{B_R} |\nabla f|^2 d^3x).
    \end{array}
\end{equation*}
\end{proof}

\subsection{The quadratic operator $Q$}\label{Subsection: The quadratic operator Q}

Recall (in subsection \ref{Subsection: The Bogomolny equations and its linearization}) that the quadratic operator $Q$ comes from the quadratic part of the variation of the Bogomolny equations. And here is a precise definition:

\begin{definition}
Suppose $\psi_1 = a_1 + \phi_1$ and $\psi_2 = a_2 + \psi_2$ are two configurations, then 
$$Q(\psi_1, \psi_2) = \dfrac{1}{4}*([a_1 \wedge a_2] + [a_2 \wedge a_1]) - \dfrac{1}{2}([a_1, \phi_2] + [a_2, \phi_2]). $$
\end{definition}

It is convenient to give the ``0 + 1" forms a Clifford algebra structure.\\

Suppose $\tau_1 = dx_1, \tau_2 = dx_2, \tau_3 = dx_3$. Then $\{1, \tau_1, \tau_2, \tau_3\}$ form an orthonormal basis of ``0 + 1" differential forms. We define the Clifford multiplication as the follows:
$$\tau_1^2 = \tau_2^2 = \tau_3^3 = -1, ~~ \tau_1\tau_2 = - \tau_2 \tau_1 = \tau_3, ~~ \tau_2\tau_3 = - \tau_3 \tau_2 = \tau_1, ~~ \tau_3\tau_1 = - \tau_1 \tau_3 = \tau_2. $$

Then we have 

$$Q(\psi_1, \psi_2) = \dfrac{1}{2} [\psi_1, \psi_2], $$
where the Lie bracket is taken with respect to the Clifford multiplication.\\

We prove the following statement in this subsection:

\begin{proposition}\label{proposition of quadratic map is bounded}
$Q$ is a bounded map from $\IH \times \IH$ to $\IL$.
\end{proposition}

\begin{proof}
The $\IL$ norm of $Q(\psi_1, \psi_2)$ is give by the square root of
$$ \dfrac{1}{4}\int_{N_{\epsilon} \backslash K} \rho^2 |[\psi_1, \psi_2]|^2 d\rho d\theta ds + \dfrac{1}{4} \int_{\IR^3 \backslash N_{\epsilon}} |[\psi_1, \psi_2]|^2 d^3x.$$

First of all,
\begin{equation*}
    \begin{array}{ll}
         & ~~~ \displaystyle \int_{N_{\epsilon} \backslash K} \rho^2|[\psi_1, \psi_2]|^2 d\rho d\theta ds \leq \int_{N_{\epsilon} \backslash K} \rho^2 |\psi_1|^2 |\psi_2|^2 d\rho d\theta ds   \\
         & \leq \displaystyle (\int_{N_{\epsilon} \backslash K} |\psi_1|^2 d\rho d\theta ds)^{\frac{1}{4}}(\int_{N_{\epsilon} \backslash K} |\psi_2|^2 d\rho d\theta ds)^{\frac{1}{4}}(\int_{N_{\epsilon} \backslash K} \rho^4 |\psi_1|^6)^{\frac{1}{4}}(\int_{N_{\epsilon} \backslash K} \rho^4|\psi_2|^6)^{\frac{1}{4}}\\
         & \leq C ||\psi_1||_{\IH}^2 ||\psi_2||_{\IH}^2.
    \end{array}
\end{equation*}
The last step is by lemma \ref{weight version of Sobolev embedding lemma}.\\

Then by lemma \ref{weight version of Sobolev embedding lemma}

$$\int_{B_R \backslash N_{\epsilon}} |[\psi_1, \psi_2]|^2 d^3x \leq \int_{B_R \backslash N_{\epsilon}} |\psi_1|^2 |\psi_2|^2 d^3x $$
$$\leq (\int_{B_R \backslash N_{\epsilon}} |\psi_1|^4 d^3x)^{\frac{1}{2}}(\int_{B_R \backslash N_{\epsilon}} |\psi_2|^4 d^3x)^{\frac{1}{2}} \leq C ||\psi_1||^2_{\IH}|| \psi_2||^2_{\IH}.$$

Finally on $\IR^3 \backslash B_R$, any configuration can be decomposed into the part that is parallel to $e_0$ and the part that is perpendicular to $e_0$. So be more precise, $\psi = \psi^{\perp} + \psi^{\parallel}$, where $\psi^{\perp} =- [[\psi, e_0], e_0]$, $\psi^{\parallel} = <\psi, e_0> e_0 = \psi - \psi^{\perp}$.\\

Then one observes that 

$$[\psi_1, \psi_2] = [\psi_1^{\perp}, \psi_2^{\perp}] + [\psi_1^{\parallel}, \psi_2^{\perp}] + [\psi_1^{\perp}, \psi_2^{\parallel}].$$

So

$$\int_{\IR^3 \backslash B_R} |[\psi_1, \psi_2]|^2 d^3x \leq \int_{\IR^3 \backslash B_R} (|\psi_1|^2|\psi_2^{\perp}|^2 + |\psi_2|^2|\psi_1^{\perp}|^2) d^3x $$
$$\leq (\int_{\IR^3 \backslash B_R} |\psi_1|^6 d^3x)^{\frac{1}{3}}(\int_{\IR^3 \backslash B_R} |\psi_2|^6 d^3x)^{\frac{1}{6}} (\int_{\IR^3 \backslash B_R} |\psi_2^{\perp}|^2 d^3x)^{\frac{1}{2}} $$
$$ + (\int_{\IR^3 \backslash B_R} |\psi_2|^6 d^3x)^{\frac{1}{3}}(\int_{\IR^3 \backslash B_R} |\psi_1|^6 d^3x)^{\frac{1}{6}} (\int_{\IR^3 \backslash B_R} |\psi_1^{\perp}|^2 d^3x)^{\frac{1}{2}}.$$
Here by lemma \ref{weight version of Sobolev embedding lemma}
$$\int_{\IR^3 \backslash B_R} |\psi_1|^6 d^3x \leq C ||\psi_1||^2_{\IH}, ~~~ \int_{\IR^3 \backslash B_R} |\psi_2|^6 d^3x \leq C ||\psi_2||^2_{\IH}.$$
And from the definition of $\IH$,

$$\int_{\IR^3 \backslash B_R} |\psi_1^{\perp}|^2 d^3 x \leq ||\psi_1||_{\IH}^2, ~~~  \int_{\IR^3 \backslash B_R} |\psi_2^{\perp}|^2 d^3 x \leq ||\psi_2||_{\IH}^2.$$

So we have

$$\int_{\IR^3 \backslash B_R} |[\psi_1, \psi_2]|^2 d^3x \leq C ||\psi_1||_{\IH}^2||\psi_2||_{\IH}^2. $$
\end{proof}

The following observation is useful in later arguments:\\

Suppose $\Psi$ and $\Psi_0$ are two fiducial configurations. Suppose $\Psi = \Psi_0 + \psi_0$. Then
$$V(\Psi_0 + \psi) = V(\Psi_0) +  \tilde{L}_{\Psi_0}(\psi) + Q(\psi, \psi) = V(\Psi) + \tilde{L}_{\Psi}(\psi_0 + \psi) + Q(\psi_0 + \psi, \psi_0 + \psi).$$
So
$$\tilde{L}_{\Psi_0}(\psi) = \tilde{L}_{\Psi}(\psi) +  2Q(\psi_0, \psi).$$

So the difference between $\tilde{L}_{\Psi}$ and $\tilde{L}_{\Psi_0}$ is given by $2Q(\psi_0, \psi)$.

\section{The Fredholm theory near the knot}\label{Section: The Fredholm theory near the knot}

In this section, we establish a local Fredholm theory for $\tilde{L}_{\Psi}$ near the knot.

\subsection{The statement of the local Fredholm theory}\label{Subsection: The statement of the local Fredholm theory}

Recall that we have fixed a smooth cut-off function $\chi$ which equals $1$ in $N_{\epsilon}$ and equals $0$ in $\IR^3 \backslash N_{2 \epsilon}$. In addition, in this section, for any $\delta < \epsilon$, suppose $\chi_{\delta}$ is a smooth cut-off function that is $1$ on $N_{\delta}$ and $0$ on $\IR^3 \backslash N_{2\delta}$ with $|\nabla \chi_{\delta}| \leq 
\dfrac{2}{\delta}$.\\

Throughout this section, we assume $\Psi$ is a configuration such that 

$$\chi(\Psi - \Psi_{\gamma, M, 0}) \in \IH.$$

The goal of this section is to prove the following statement:\\

Roughly speaking, supposing $\gamma \neq \dfrac{1}{4}$, then when $\delta$ is small enough $\tilde{L}_{\Psi}$ is locally an ``isomorphism" from $\IH$ to $\IL$ in $N_{\delta}$.  Precisely, we have

\begin{proposition}\label{proposition of local Fredholmness}
Suppose $\gamma \neq \dfrac{1}{4}$. And suppose $\Psi$ has a knot singularity with monodromy $\gamma$ in $K$. Assume $\delta$ is small enough, which may depend on $\Psi$. Then there is a small $c > 0$ such that
$$||\tilde{L}_{\Psi}(\chi_{\delta} \psi)||_{\IL} \geq c ||\chi_{\delta} \psi||_{\IH}, ~~~ ||\tilde{L}_{\Psi}^{\dagger} ( \chi_{\delta} \psi)||_{\IL} \geq c ||\chi_{\delta} \psi||_{\IH},$$
where $\tilde{L}_{\Psi}^{\dagger}$ is the formal ad-joint operator of $\tilde{L}_{\Psi}$.
\end{proposition}

\begin{proposition}\label{proposition of formal adjoint and cokernel}
Suppose $\psi$ is an element in $\IL$ that satisfies
$$\tilde{L}^{\dagger}_{\Psi}(\rho \psi) = 0 ~~~ \text{on} ~~ N_{2\epsilon}\backslash K. $$
Then $\rho \chi \psi \in \IH$.  
\end{proposition}

The proof of the above two propositions occupies the remaining of this section.

\subsection{Compare the general case with the model case}\label{Subsection: Compare the general case with the model case}

This subsection compares $\tilde{L}_{\Psi}$ and $\tilde{L}_{\Psi_{\gamma, M, 0}}$ near the knot, where $\Psi$ is the model solution. Recall that $\chi(\Psi - \Psi_{\gamma, M, 0}) \in \IH$.

\begin{lemma}\label{lemma that R delta is a small operator}
 Let $R_{\delta}$ be the following operator from $\IH$ to $\IL$:

$$R_{\delta}(\psi) = \tilde{L}_{\Psi}(\chi_{\delta} \psi) - \tilde{L}_{\Psi_{\gamma, M, 0}}(\chi_{\delta} \psi).$$

Then for any $c > 0$, there is a $\delta_0 > 0$ such that when $\delta < \delta_0$, for any $\psi \in \IH$
$$||R_{\delta}(\psi)||_{\IL} \leq c||\psi||_{\IH} $$
The same statement is also true if we replace $R_{\delta}$ with $R^{\dagger}_{\delta}$.
\end{lemma}

\begin{proof}
We have $\tilde{L}_{\Psi} (\psi)= D_A \psi + [\Phi, \psi] ~~~ \text{and} ~~\tilde{L}_{\Phi_{\gamma, M, 0}} (\psi) = D_{A_{\gamma, M, 0}} \psi + [\Psi_{\gamma, M, 0}, \psi]$.\\

So their difference is

$$[(A - A_{\gamma, M, 0}) \wedge \psi] + [\Phi - \Phi_{\gamma, M, 0}, \psi].$$

Thus $$R_{\delta}(\psi) = [\chi_{\delta}(A - A_{\gamma, M, 0}) \wedge (\psi)] + [\chi_{\delta}(\Phi - \Phi_{\gamma, M, 0}), \psi].$$

Similarly, we have
$$R^{\dagger}_{\delta}(\psi) = [\chi_{\delta}(A - A_{\gamma, M, 0}) \wedge (\psi)] - [\chi_{\delta}(\Phi - \Phi_{\gamma, M, 0}), \psi].$$

Both are bilinear expressions in terms of $\chi_{\delta}(\Psi - \Psi_{\gamma, M, 0})$ and $\psi$ and are supported in $N_{\epsilon}$. By lemma \ref{weight version of Sobolev embedding lemma}, there is a constant $c_1$ which doesn't depend on $\delta$ such that, 
$$(\int_{N_{2\delta}}\rho^3|\chi_{\delta}(\Psi - \Psi_{\gamma, M, 0})|^4 d\rho d\theta ds)^{\frac{1}{4}} \leq c_1 ||\chi_{\delta}(\Psi - \Psi_{\gamma, M, 0})||_{\IH} $$
and
$$(\int_{N_{2\delta}}\rho^3|\psi|^4 d\rho d\theta ds)^{\frac{1}{4}} \leq c_1 ||\psi||_{\IH}.$$

So by Cauchy-Schwartz inequality,
$$||R_{\delta}(\psi)||_{\IL} \leq  2(\int_{N_{2\delta}}\chi^2_{\delta}|\Psi - \Psi_{\gamma, M, 0}|^2|\psi|^2 \rho^3d\rho d\theta ds)^{\frac{1}{2}} \leq 2 c_1^2 ||\chi_{\delta}(\Psi - \Psi_{\gamma, M, 0})||_{\IH} ||\psi||_{\IH}.$$

Note that 
$$||\chi_{\delta}(\Psi - \Psi_{\gamma, M, 0})||_{\IH}^2 \leq \int_{N_{2\delta}} (\rho^2|(\nabla (\chi_{\delta}(\Psi - \Psi_{\gamma, M, 0}))|^2 + \chi_{\delta}^2|\Psi - \Psi_{\gamma, M, 0}|^2) d\rho d\theta ds$$ $$\leq 16\int_{N_{2\delta}} (\rho^2|\nabla(\Psi - \Psi_{\gamma, M, 0})|^2 + |\Psi - \Psi_{\gamma, M, 0}|^2) d\rho d\theta ds $$
which goes to $0$ as $\delta \rightarrow 0$, since $\chi(\Psi - \Psi_{\gamma, M, 0})$ is in $\IH$. Thus when $\delta$ is small enough, 

$$||R_{\delta}(\psi)||_{\IL} \leq c||\psi||_{\IH}~~ \text{and} ~~ ||R^{\dagger}_{\delta}(\psi)||_{\IL} \leq c||\psi||_{\IH}.$$ 
\end{proof}

\subsection{Compare $\tilde{L}_{\Psi_{\gamma, M, 0}}$ with the operator $N$}\label{Subsection: Compare L with the operator N}

In this subsection, we assume the fiducial configuration $\Psi = \Psi_{\gamma, M, 0}$. We further simplifie $\tilde{L}_{\Psi}$ to an operator $N$. Before we define the operator $N$, we need to introduce some notations first.\\ 

Suppose a configuration is $\psi = \phi - a_1 dz_1 - a_2 dz_2 - a_3 ds$. Then near the knot, it is convenient to use a pair of complex-valued functions $(\alpha, \beta)$ to represent $\psi$, where $\alpha = a_1 - ia_2, ~~ \beta = a_3 + i\phi$. Both $\alpha, \beta$ are $sl(2, \IC) = su(2) \otimes \IC$-valued functions. \\

 Consider the following three operators: $\partial_{\theta}, ~ \partial_s, ~[i \sigma, \cdot]$ acting on $sl(2, \IC)$-valued functions. They are pairwise commutative and they have eigenvalues $(im, ik, 2j)$ on
the common ``eigenvector" $e^{-im\theta - iks}h_{j}$, where $m, \dfrac{kl}{2\pi} \in \IZ$ and $j \in \{\pm 1, 0\}$. Here $l$ is the length of the knot.

\begin{definition}

Any smooth $sl(2, \IC)$-valued funtion $\alpha$ on $N_{\epsilon} \backslash K$ can be decomposed into infinite sums of those common ``eigenvectors" with coefficients that depend only on $\rho$. Each such common eigenvector with its coefficients is called a Fourier component of $\alpha$.
\end{definition}

A typically Fourier component of $\begin{pmatrix}
 \alpha \\ \beta
\end{pmatrix}$ looks like $$\begin{pmatrix}
 \alpha \\ \beta
\end{pmatrix} = \begin{pmatrix}
 u(\rho) e^{-i(m + 1)\theta - iks} \\  v(\rho) e^{-im \theta - iks}
\end{pmatrix} .$$

Note that we have paired the $e^{- (m+1)\theta}$ component for $\alpha$ with the $e^{-im\theta}$ component for $\beta$. Because we want to define the operator $N$ as follows: Let $N_j$ be the operator $N$ acting on a Fourier component corresponding to $j$. (Recall that the Fourier component is a common eigenvector of $\partial_{\theta}, \partial_s, [i\sigma, \cdot]$  with eigenvalues $im, ik, 2j$ respectively). Then
$$N_j\begin{pmatrix}
\alpha \\ \beta
\end{pmatrix}= \begin{pmatrix}
 - (\partial_s + 2Mj) \alpha  + (\partial_{\rho} - \dfrac{i\partial_{\theta} }{\rho} - \dfrac{2\gamma j}{\rho}) \beta  \\
      (\partial_{\rho} + \dfrac{i\partial_{\theta}}{\rho} + \dfrac{ 2\gamma j}{\rho})\alpha + (\partial_s -2Mj)\beta 
\end{pmatrix}.$$
Or more precisely,
$$N_j\begin{pmatrix}
ue^{-i(m+1)\theta - iks} \\ v e^{-im\theta - iks}
\end{pmatrix}= \begin{pmatrix}
 (( \partial_{\rho} - \dfrac{(m  + 2\gamma j)}{\rho}) v +(ik - 2Mj)u) e^{-i(m+1)\theta - iks}   \\
      ((\partial_{\rho} + \dfrac{(m+1 + 2\gamma j)}{\rho})u + (- ik -2Mj) v)e^{-im\theta - iks} 
\end{pmatrix},$$
where the subscript $j$ indicates that we are working on the Fourier component such that $[i\sigma, \cdot ] = 2j$ for $j \in \{\pm 1, 0\}$.\\

\begin{remark}
Here is intuition behind the operator $N$: If we replace the Euclidean metric with the following cylindrically flat metric on $N_{\epsilon}$: 
$$g' = dz_1^2 + dz_2^2 + ds^2 = d\rho^2 + \rho^2 d\theta^2 + ds^2, $$
then one can easily check that $N$ is the  linearization $\tilde{L}_{\Psi_{\gamma, M, 0}}$ using the cylindrically flat metric. 
\end{remark}

Come back to the Euclidean metric 
$$g = dx_1^2 + dx_2^2 + dx_3^2.$$ The following proposition describes the difference between $\tilde{L}_{\Psi}$ and $N$. 

\begin{proposition}\label{proposition that N and L differ by lower order terms}
$$N (\psi) - \tilde{L}_{\Psi}(\psi)  = \rho A (\nabla \psi) + B(\psi), $$
where both $A$ and $B$ are linear operators with bounded smooth coefficients. A similar proposition holds for $N^{\dagger} - \tilde{L}_{\Psi}^{\dagger}$ as well.
\end{proposition}

\begin{proof}

Under the metric $g'$, consider an orthonormal basis $\{ \partial_s,  \partial_{z_1}, \partial_{z_2}\}$. Recall that in subsection \ref{Subsection: The ambient space}, $(s, z_1, z_2)$ represents the point $K(s) + z_1\e_1(s) + z_2\e_2(s)$. So 
$$\partial_s = K'(s) + z_1\e_1'(s) + z_2\e_2'(s) = (1 - z_1k)K' + z_1 \tau \e_2 - z_2\tau \e_1, ~~ \partial_{z_1} = \e_1(s), ~~\partial_{z_2} = \e_2(s). $$

Under the Euclidean metric $g$,
$$|\partial_s|^2 = (1 - z_1 k)^2 + \rho^2 \tau^2 , ~~ |\partial_{z_1}|^2 = |\partial_{z_2}|^2 = 1, ~~ <\partial_s, \partial_{z_1}> = - z_2 \tau , ~~ <\partial_s, \partial_{z_2}> = z_1 \tau. $$
Thus 
$$|g - g'| = O(\rho). $$
Moreover
$$\nabla (|\partial_s|^2) = O(1), ~ \nabla |\partial_{z_1}|^2 = \nabla |\partial_{z_2}|^2 = 0, ~~
 \nabla <\partial_s, \partial_{z_1}> = O(1),~~ \nabla <\partial_s, \partial_{z_2}> = O(1). $$
 So $$|\nabla g - \nabla g'| = O(1).$$
 
 Recall that the definition of $$\tilde{L}_{\Psi} = D_A + [\Phi, \cdot], ~~~ \tilde{L}^{\dagger}_{\Psi} = D_A - [\Phi, \cdot], $$ where $D_A$ is the Dirac operator twisted by the connection $A$, which is the only part that depends on the metric. Thus the difference between $N$ and $\tilde{L}_{\Psi}$ (or $N^{\dagger}$ and $\tilde{L}^{\dagger}_{\Psi}$) comes from the difference between the Dirac operator under two metrics $g$ and $g'$. Thus the proposition follows directly.
\end{proof}

Finally we have

\begin{lemma}\label{lemma that R' is a small operator}
Let
$$R_{\delta}' (\psi) = N(\chi_{\delta} \psi) - \tilde{L}_{\Psi}(\chi_{\delta} \psi), ~~~ R^{\dagger}{}'_{\delta}(\psi) = N^{\dagger}(\chi_{\delta} \psi) - \tilde{L}^{\dagger}_{\Psi}(\chi_{\delta} \psi).$$
Note that $R^{\dagger}{}'_{\delta}$ is not the formal ad-joint of $R'_{\delta}$.\\

Then for any $c > 0$, there is a $\delta_0 > 0$ such that when $\delta < \delta_0$, for any $\psi \in \IH$
$$||R'_{\delta}(\psi)||_{\IL} \leq c||\psi||_{\IH}, ~~ ||R^{\dagger}{}'_{\delta}(\psi)||_{\IL} \leq c||\psi||_{\IH}. $$
\end{lemma}

\begin{proof}
From proposition \ref{proposition that N and L differ by lower order terms}
$$|R'_{\delta}(\psi)| = |\rho A(\nabla(\chi_{\delta}\psi)) + B(\chi_{\delta}\psi)| \leq C (\rho |\nabla \psi| + |\psi|).$$
So
$$||R'_{\delta}(\psi)||_{\IL}^2 \leq C \int_{N_{2\delta} \backslash K} (\rho^4 |\nabla \psi|^2 + \rho^2|\psi|^2) d\rho d\theta ds \leq 2C \delta \int_{N_{2\delta} \backslash K} (\rho^2 |\nabla \psi|^2 + |\psi|^2) d\rho d\theta ds \leq 2C \delta ||\psi||^2_{\IH}.$$

Since $C$ doesn't depend on $\delta$, the proposition follows directly. And the same statement for $R^{\dagger}{}'_{\delta}$ follows for the same reason.

\end{proof}

\subsection{The operator $N$ as a local isomorphism}\label{Subsection: The operator N as a local isomorphism}

This subsection shows that $N$ is a ``local isomorphism" from $\IH$ to $\IL$ near the knot. Here is the precise statement.

\begin{proposition}
Suppose $\psi$ is a smooth configuration on $\IR^3 \backslash K$ whose support is in $N_{\epsilon}$. And suppose $\psi$ has a finite $\IH$ norm. Then there exists a small constant $c > 0$ which doesn't depend on $\psi$ or $\epsilon$ (assuming that $\epsilon$ is small enough), such that

$$||N (\psi)||_{\IL} \geq 3c ||\psi ||_{\IH}, ~~~ \text{and} ~~ ||N^{\dagger} \psi ||_{\IL} \geq 3c ||\psi||_{\IH},$$
where $N^{\dagger}$ is the formal ad-joint of $N$.

\end{proposition}

\begin{proof}
We only prove that the $\IL$ norm of $N(\psi)$ bounds the $\IH$ norm of $\psi$. And the same argument also works for $N^{\dagger}$ which is omitted.\\

It suffices to assume that $\psi = \begin{pmatrix}
 \alpha \\ \beta
\end{pmatrix}$ has only one Fourier component, that is to say $$\begin{pmatrix}
 \alpha \\ \beta
\end{pmatrix} = \begin{pmatrix}
 u(\rho) e^{-i(m + 1)\theta - iks} \\  v(\rho) e^{-im \theta - iks}
\end{pmatrix} .$$.\\

The $\IH$ norm of $\psi$ is equivalent to (the square root of)

$$\int_0^{\epsilon} \rho^2 (|\partial_{\rho} u|^2 + |\partial_{\rho} v|^2) + (k^2\rho^2 + m^2 + 1)( |u|^2 + |v|^2) d\rho. $$

In addition, since the above integral is finite (because we assume $\psi \in \IH$ a prior), by lemma \ref{boundary convergent lemma}, $$\lim\limits_{\rho \rightarrow 0} \rho (|u|^2 + |v|^2)  = 0.$$

The $\IL$ norm of $N(\psi)$ is equivalent to (the squre root of)

$$\int_0^{\epsilon} \rho^2 (|\partial_{\rho}u + \dfrac{\lambda + 1}{\rho}u + (-ik + 2Mj)v|^2 + |\partial_\rho v - \dfrac{\lambda}{\rho} v + (ik + 2Mj)u |^2) d\rho,$$
where $\lambda = m + 2\gamma j$.\\

The $\IL$ norm of $N(\psi)$ can be re-written as
$$\int_0^{\epsilon} (\rho^2(|\partial_{\rho} u|^2 + |\partial_{\rho} v|^2) + (k^2 + 4M^2j^2)\rho^2(|u|^2 + |v|^2) + ((\lambda + 1)^2|u|^2 + \lambda^2 |v|^2) + 2\rho <u, (-ik + 2Mj)v)>$$
$$+ \rho \partial_{\rho}((\lambda + 1)|u|^2 - \lambda |v|^2) + 2 \rho^2\partial_{\rho}(< u, (-ik + 2Mj)v>)) d\rho$$
$$= \int_0^{\epsilon} (\rho^2(|\partial_{\rho} u|^2 + |\partial_{\rho} v|^2) + (k^2 + 4M^2j^2)\rho^2 (|u|^2 + |v|^2) + \lambda (\lambda + 1)(|u|^2 +  |v|^2) $$
$$ - 2 \rho (< u, (-ik + 2Mj)v>)) d\rho + (\text{boundary terms}).$$

The boundary terms vanish because $\lim\limits_{\rho \rightarrow 0} \rho (|u|^2 + |v|^2)  = 0$ and because $\psi$ is supported in $N_{\epsilon}$. The boundary terms in the remaining derivations will also vanish for the same reason. And we'll not address them again.\\

The $\IL$ norm of $N(\psi)$ is greater or equal to

$$\int_0^{\epsilon} (\rho^2(|\partial_{\rho} u|^2 + |\partial_{\rho} v|^2) + ((k^2 + 4M^2j^2)\rho^2 + \lambda (\lambda + 1) - \rho \sqrt{k^2 + 4M^2})(|u|^2 +  |v|^2) d\rho.$$

Suppose $0 < \mu < \dfrac{1}{2}$ and suppose $t$ is a positive number. The key observation is the following inequality:

$$\int_0^{\epsilon} \rho^2 |\partial_{\rho} u|^2 + (\mu^2 - \mu + t^2 \rho^2 - (2 - 2\mu) t\rho)|u|^2 d\rho = \int_0^{\epsilon}|\rho \partial_{\rho}u + \mu u + t \rho u|^2 d\rho \geq 0.$$

Similar inequality also holds for $v$. It can be directly verified that, one can choose $\mu$ that is slightly less than $\dfrac{1}{2}$ (but greater than $\min\{\gamma, \dfrac{1}{2} - \gamma\}$), $t$ slightly less that $\sqrt{k^2 + 4M^2j^2} > 0$ or $t = 0$ if $k = Mj = 0$ and conclude that, there is a small constant $c > 0$ which doesn't depend on $\lambda$ or $k$, such that

$$\int_0^{\epsilon} (\rho^2(|\partial_{\rho} u|^2 + |\partial_{\rho} v|^2) + ((k^2 + 4M^2j^2)\rho^2 + \lambda (\lambda + 1) - \rho \sqrt{k^2 + 4M^2})(|u|^2 +  |v|^2) d\rho$$
$$\geq c \int_0^{\epsilon} \rho^2 (|\partial_{\rho} u|^2 + |\partial_{\rho} v|^2) + (k^2\rho^2 + m^2 + 1)( |u|^2 + |v|^2) d\rho.  $$
\end{proof}

\subsection{Proof of the main statements (proposition \ref{proposition of local Fredholmness} and proposition \ref{proposition of formal adjoint and cokernel})}\label{Proof of two propositions}

Now the proofs are straightforward:

\paragraph{proof of proposition \ref{proposition of local Fredholmness}}

We have shown that 

$$||N(\chi_{\delta} \psi)||_{\IL}  \geq 3c ||\chi_{\delta} \psi||_{\IH}.$$
And we have $$\tilde{L}_{\Psi}(\chi_{\delta} \psi) = N(\chi_{\delta} \psi) - R'_{\delta}(\psi) - R_{\delta}(\psi).$$

By lemma \ref{lemma that R delta is a small operator} and lemma \ref{lemma that R' is a small operator}, when $\delta$ is small enough, 
$$||R_{\delta}(\psi)||_{\IL} \leq c||\psi||_{\IH}, ~~~  ||R'_{\delta}(\psi)||_{\IL} \leq c||\psi||_{\IH}.$$

So

$$||\tilde{L}_{\Psi}\psi||_{\IL} \geq c  ||\psi||_{\IH}.$$

Thus the first part of proposition \ref{proposition of local Fredholmness} is true. The second part (the same statement for $\tilde{L}_{\Psi}^{\dagger}$) is also true for exactly the same reason.

\paragraph{proof of proposition \ref{proposition of formal adjoint and cokernel}}
Suppose $\psi \in \IL$. And suppose $\tilde{L}^{\dagger}_{\Psi}(\rho \psi) = 0$ on $N_{2\epsilon \backslash K}$. 
By a standard elliptic regularity argument, $\rho \psi$ is smooth on any compact subset of $N_{2\epsilon }\backslash K$. Suppose $\delta \ll \epsilon$. Let $\chi_{\delta}$ be a smooth cut-off function that is $1$ on $N_{\delta}$ and equals $0$ on $N_{2\delta}$. Moreover, we assume $|\nabla \chi_{\delta}| \leq \dfrac{2}{\delta}$. Then 

$$||\tilde{L}_{\Psi} (\chi (1 - \chi_{\delta}) \rho \psi)||_{\IL}^2 \leq C \int_{N_{2\epsilon} \backslash K}\rho^4(|\nabla \chi|^2 + |\nabla \chi_{\delta}|^2)|\psi|^2 d\rho d\theta ds \leq C ||\psi||_{\IL}^2, $$
where $C$ doesn't depend on $\delta$, but changes in different expressions. Form proposition \ref{proposition of local Fredholmness},
$$||\tilde{L}_{\Psi} (\chi (1 - \chi_{\delta}) \rho \psi)||_{\IL} \geq c ||\chi (1 - \chi_{\delta}) \rho \psi||_{\IH}.$$
So $||\chi (1 - \chi_{\delta}) \rho \psi||_{\IH}$ is bounded above uniformly in $\delta$.
Finally, 
$$||\chi(1 - \chi_{\delta}) \rho \psi||_{\IH}^2 \geq \int_{N_{2\epsilon} \backslash N_{2\delta}} \rho^2(|\partial_{\rho}(\chi \rho \psi)|^2 + |\partial_s(\chi \rho \psi)|^2 ) + |\partial_{\theta}(\chi \rho \psi)|^2 + |\chi \rho \psi|^2 d\rho d\theta ds.$$
So $$\liminf\limits_{\delta \rightarrow 0}||\chi(1 - \chi_{\delta}) \rho \psi||_{\IH}^2 \geq ||\chi \rho\psi||_{\IH}^2.$$
So $\chi \rho \psi \in \IH$.

\section{The global Fredholm theory and the moduli space}\label{Section: The golbal Fredholm theory and the moduli space}

\subsection{The Fredholm theory near infinity}\label{Subsection: The Fredholm theory near infinity}
This subsection is standard based on Taubes' thesis \cite{Jaffe1980VorticesTheories} with mild modifications. But we state it here to make this paper self-contained. In fact, the ideas here are also similar with the arguments in Section \ref{Section: The Fredholm theory near the knot} but are much simpler.\\

Suppose $\chi_r$ is a fixed cut-off function which is $1$ on $B_r$ and $0$ on $\IR^3 \backslash B_{2r}$. Suppose $|\nabla \chi_r| \leq \dfrac{2}{r}$. Recall that we have fixed large enough $R$. The goal of this subsection is to prove the following:

\begin{proposition}\label{proposition of Fredholmness at infinity}
There is a constant $c > 0$ such that and a large enough $r_0 > 0$ (which may depend on $\Psi$) such that, if $r > r_0$, then for any $\psi \in \IH$,
$$||\tilde{L}_{\Psi}((1 - \chi_r)\psi)||_{\IL} \geq c ||(1 - \chi_r)\psi||_{\IH}, ~~ ||\tilde{L}^{\dagger}_{\Psi}((1 - \chi_r)\psi)||_{\IL} \geq c ||(1 - \chi_r)\psi||_{\IH},  $$
where $\tilde{L}_{\Psi}^{\dagger}$ is the formal ad-joint operator of $\tilde{L}_{\Psi}$.
\end{proposition}

First of all, we need a Weitzenbock type formula for $\tilde{L}_{\Psi}$.\\

We may write ``$1, \tau_1 = dx_1, \tau_2 = dx_2, \tau_3= dx_3$" as a basis of a Clifford algebra and we use the standard Clifford multiplication. Then the Dirac operator can be written as
$$D_A \psi = (\sum\limits_{i=1}^3 \tau_i (\nabla_{A_i}\psi)), $$
where $\nabla_{A_i}$ is the covariant derivative $\nabla_A$ in the $x_i$ direction.

\begin{proposition}\label{proposition: Weitzenbock identity}
We have the following identities
$$\tilde{L}^{\dagger}_{\Psi}\tilde{L}_{\Psi}(\psi) = \nabla_A^{\dagger}\nabla_A \psi + [*F_A + \nabla_{A}\Phi, \psi] - [\Phi, [\Phi, \psi]],$$ $$\tilde{L}_{\Psi} \tilde{L}^{\dagger}_{\Psi} (\psi) =  \nabla_A^{\dagger}\nabla_A \psi + [*F_A - \nabla_{A}\Phi, \psi] - [\Phi, [\Phi, \psi]],$$
where $F_A$ is the curvature of $A$ and the Lie bracket is defined using the Clifford multiplication.
\end{proposition}

\begin{proof}
Recall that $$\tilde{L}_{\Psi} \psi = D_A \psi + [\Phi, \psi], ~~~ \tilde{L}^{\dagger}_{\Psi}\psi = D_A \psi - [\Phi, \psi].$$
Then
$$\tilde{L}^{\dagger}_{\Psi}\tilde{L}_{\Psi}(\psi) = (\sum\limits_{i = 1}^3  \sum\limits_{j = 1}^3 \tau_i \tau_j \nabla_{A_i} \nabla_{A_j} \psi) + \sum\limits_{i = 1}^3 (\tau_i \nabla_{A_i} ([\Phi, \psi]) - [\Phi, \nabla_{A_i}\psi])) - [\Phi, [\Phi, \psi]]$$
$$= - (\sum\limits_{i = 1}^3 \nabla_{A_i}^2 \psi) + \sum\limits_{i = 1}^3(\tau_i[F_i + \nabla_i \Phi, \psi])- [\Phi, [\Phi, \psi]]$$ $$= \nabla_A^{\dagger}\nabla_A \psi + [*F_A + \nabla_{A}\Phi, \psi] - [\Phi, [\Phi, \psi]], $$
where $F_A = *(F_1dx_1 + F_2 dx_2 + F_3 dx_3)$ is the curvature of $A$.\\

The other identity for $\tilde{L}_{\Psi} \tilde{L}^{\dagger}_{\Psi}$ follows for the same reason.
\end{proof}

\begin{lemma}\label{lemma that P is small}
Let $$P_r(\psi) =  \tilde{L}_{\Psi}((1 - \chi_r)\psi) - \tilde{L}_{\Psi_{\gamma, M, k}}((1 - \chi_r)\psi). $$
For any $c > 0$, there exists a large enough $r_0 > 0$ such that, if $r > r_0$ and $\psi \in \IH$ then
$$||P_r(\psi)||_{\IL} \leq c || \psi||_{\IH}, ~~~ ||P^{\dagger}_r(\psi)||_{\IL} \leq c || \psi||_{\IH}.$$
\end{lemma}

\begin{proof}
Suppose $\Psi = A + \Phi$, $\Psi_{\gamma, M, k} = A_{\gamma, M, k} + \Phi_{\gamma, M, k}$. Then
$$P_r(\psi) = 2Q((1 - \chi_r) (\Psi - \Psi_{\gamma, M, k}), \psi) + P_0, $$
where $Q$ is given by subsection \ref{Subsection: The Bogomolny equations and its linearization} and has been studied in subsection \ref{Subsection: The quadratic operator Q}.
We have assumed that $\Psi - \Psi_{\gamma, M, k} \in \IH$. So
$$||(1 - \chi_r)(\Psi - \Psi_{\gamma, M, k})||_{\IH}^2 = \int_{\IR^3 \backslash B_r} |\nabla ((1 - \chi_r)(\Psi - \Psi_{\gamma, M, k}))|^2 + (1 - \chi_r)^2|[e_0, \Psi - \Psi_{\gamma, M, k}]|^2 d^3x $$
$$\leq 2 \int_{\IR^3 \backslash B_r} (|\nabla (\Psi - \Psi_{\gamma, M, k})|^2 + |[e_0, \Psi - \Psi_{\gamma, M, k}]|^2) d^3 x + 8 \int_{B_{2r} \backslash B_r} \dfrac{1}{r^2}|\Psi - \Psi_{\gamma, M, k}|^2 d^3x.$$
Since $\Psi - \Psi_{\gamma, M, k} \in \IH$, the above term goes to $0$ when $r \rightarrow +\infty$. Thus when $r$ is large enough, by proposition \ref{proposition of quadratic map is bounded},

$$||Q((1 - \chi_r) (\Psi - \Psi_{\gamma, M, k}), \psi)||_{\IL} \leq c ||\psi||_{\IH}.$$

\end{proof}

\begin{lemma}\label{Lemma of local isomorphism at infinity}
There is a small enough $c > 0$, there exists an $r_0 > R$ such that when $r \geq r_0$, for any $\psi \in \IH$, we have
$$|| \tilde{L}_{\Psi_{\gamma, M, k}} ((1 - \chi_r)\psi)||_{\IL} \geq 2c ||(1 - \chi_r) \psi||_{\IH}. $$
Same inequality also holds for $\tilde{L}^{\dagger}_{\gamma, M, k}$.
\end{lemma}

\begin{proof}
We may assume $\psi$ has bounded support. Note that the support of $(1 - \chi_r)\psi$ is also away from the knot. So there are no boundary terms when doing integrals on $(1 - \chi_r)\psi$ over the entire $\IR^3$. In particular, by Weitzenbock identities in proposition \ref{proposition: Weitzenbock identity},
$$\int_{\IR^3} |\tilde{L}_{\Psi_{\gamma, M, k}} (1 - \chi_r) \psi|^2 d^3x = \int_{\IR^3} (1 - \chi_r)<\tilde{L}^{\dagger}_{\Psi_{\gamma, M, k}}\tilde{L}_{\Psi_{\gamma, M, k}} ((1 - \chi_r) \psi), \psi> d^3x $$
$$= \int_{\IR^3}(|\nabla_A ((1 - \chi_r)\psi)|^2 + |[\Phi_{\gamma, M, k}, (1- \chi_r)\psi]|^2 + (1 - \chi_r)^2<[*F_A - d_A \Phi, \psi], \psi>)d^3x. $$
Note that one can check directly that both $|F_A|$ and $|d_A \Phi|$ are in $L^2(\IR^3 \backslash B_r)$. And moreover, $$|<[*F_A - d_A \Phi, \psi], \psi>| \leq (|F_A| + |d_A\Phi|)[\psi \wedge \psi] \leq 2(|F_A| + |d_A \Phi|) Q(\psi, \psi), $$
where $Q$ is defined in subsection \ref{Subsection: The Bogomolny equations and its linearization} and studied in subsection \ref{Subsection: The quadratic operator Q}. So

$$\int_{\IR^3} (1 - \chi_r)^2|*F_A - d_A \Phi| |Q(\psi, \psi)| d^3 x \leq C(\int_{\IR^3 \backslash B_r} (|F_A| + |d_A\Phi|)^2 d^3x )^{\frac{1}{2}} ||\psi||_{\IH}^2,$$
where $C$ doesn't depend on $r$ and 
$$ \lim\limits_{r \rightarrow +\infty} \int_{\IR^3 \backslash B_r} (|F_A| + |d_A\Phi|)^2 d^3x  = 0.$$

Since $r$ is large, we may assume
$$C\int_{\IR^3} (1 - \chi_r)^2|*F_A - d_A \Phi| |Q(\psi, \psi)| d^3 x \leq \dfrac{1}{4} \min\{1, |M|^2\}.$$

Moreover, we may assume $\Phi_{\gamma, M, k} = M e_0$ on $\IR^3 \backslash B_r$. Thus
$$\int_{\IR^3} |\tilde{L}_{\Psi_{\gamma, M, k}} (1 - \chi_r) \psi|^2 d^3x \geq \dfrac{1}{2} \int_{\IR^3}(|\nabla_A ((1 - \chi_r)\psi)|^2 + M^2|[e_0, (1 - \chi_r)\psi]|^2)d^3x. $$
Clearly letting $c = \dfrac{1}{4} \min\{1, |M|\}$ will do the job. The same inequality holds for $\tilde{L}_{\Psi_{\gamma, M, k}}$ for the same reason.
\end{proof}

\paragraph{proof of proposition \ref{proposition of Fredholmness at infinity}}

From lemma \ref{Lemma of local isomorphism at infinity},

$$|| \tilde{L}_{\Psi_{\gamma, M, k}} ((1 - \chi_r)\psi)||_{\IL} \geq 2c ||(1 - \chi_r) \psi||_{\IH}. $$

From lemma \ref{lemma that P is small}

$$||P_r(\psi)||_{\IL} \leq c || \psi||_{\IH}.$$

So 
$$|| \tilde{L}_{\Psi} ((1 - \chi_r)\psi)||_{\IL} \geq c ||(1 - \chi_r) \psi||_{\IH}. $$
Same inequality also holds for $\tilde{L}^{\dagger}_{\Psi}$. Thus proposition \ref{proposition of Fredholmness at infinity} is true.\\

\begin{corollary}\label{corollary of cokernel in H near infinity}
If $\psi \in \IL$ is a smooth configuration such that $\tilde{L}^{\dagger}_{\Psi} \psi = 0$ on $\IR^3 \backslash B_R$. Then for large enough $r$, $(1 - \chi_r) \psi \in \IH$.
\end{corollary}

\begin{proof}
We cannot use proposition \ref{proposition of Fredholmness at infinity} directly since we don't know a priori that $(1 - \chi_R)\psi$ is in $\IH$. But assuming $r_1 \gg r \gg R$, $\chi_{r_1}(1 - \chi_{r})\psi$ is smooth and has bounded support. Thus $\chi_{r_1}(1 - \chi_{r})\psi$  is in $\IH$. And by proposition \ref{proposition of Fredholmness at infinity},
$$|| \tilde{L}^{\dagger}_{\Psi} (\chi_{r_1}(1 - \chi_r)\psi)||_{\IL} \geq c ||\chi_{r_1}(1 - \chi_r) \psi||_{\IH}. $$

On the one hand,
$$|| \tilde{L}^{\dagger}_{\Psi} (\chi_{r_1}(1 - \chi_r)\psi)||_{\IL}^2 \leq 2||\tilde{L}^{\dagger}_{\Psi} ((1 - \chi_r)\psi) ||_{\IL}^2 + 2\int_{B_{2r_1} \backslash B_{r_1}} |\nabla \chi_{r_1}|^2 |\psi|^2 d^3x $$
$$\leq 2||\tilde{L}^{\dagger}_{\Psi} ((1 - \chi_r)\psi) ||_{\IL}^2 + \dfrac{C}{r_1^2}\int_{B_{2r_1} \backslash B_{r_1}} |\psi|^2 d^3x. $$
Note that since $\psi \in \IL$, we have 
$$\lim\limits_{r_1 \rightarrow +\infty} \dfrac{C}{r_1^2}\int_{B_{2r_1} \backslash B_{r_1}} |\psi|^2 d^3x = 0.$$

And since $\tilde{L}_{\Psi}^{\dagger} \psi = 0$, $\tilde{L}^{\dagger}_{\Psi} ((1 - \chi_r)\psi)$ has bounded support and hence
$$||\tilde{L}^{\dagger}_{\Psi} ((1 - \chi_r)\psi)||_{\IL} < +\infty $$
So $$\limsup\limits_{r_1 \rightarrow +\infty} || \tilde{L}^{\dagger}_{\Psi} (\chi_{r_1}(1 - \chi_r)\psi)||_{\IL}^2 < +\infty. $$
On the other hand,
$$||\chi_{r_1}(1 - \chi_r) \psi||_{\IH}^2 \geq \int_{B_{r_1}} |\nabla ((1 - \chi_r)\psi)|^2  + |[e_0, (1 - \chi_r)\psi]|^2 d^3 x. $$
So let $r_1 \rightarrow +\infty$, $B_{r_1} \rightarrow \IR^3$, and
$$\liminf\limits_{r_1 \rightarrow +\infty} ||\chi_{r_1}(1 - \chi_r) \psi||_{\IH}^2 \geq ||(1 - \chi_r) \psi||_{\IH}^2.  $$

Thus 
$$ ||(1 - \chi_r) \psi||_{\IH} < +\infty.$$

\end{proof}

\subsection{The global Fredholm theory}\label{Subsection: The global Fredholm theory}

In this subsection, we assume $\gamma \neq \dfrac{1}{4}$ and $M > 0$. The goal of this subsection is to prove that $\tilde{L}_{\Psi}$ is a Fredholm operator from $\IH$ to $\IL$.\\

To start, we modify $\IL$ a little bit as the follows: Let $w = \chi \rho + (1 - \chi)$. Then $w = \rho$ on $N_{\epsilon}$ and $w = 1$ on $\IR^3 \backslash N_{2\epsilon}$. Define the $\IL_{w}$ norm of a configuration $\psi$ to be the square root of 
$$\int_{\IR^3 \backslash K} w |\psi|^2 d^3x.$$

It is obvious that $\IL_{w}$ and $\IL$ are equivalent norms.\\

First of all, we show that
\begin{lemma}
Suppose $\delta$ is small enough and $r$ is large enough. Then there exists a small constant $c$ and a large constant $C$ such that for any $\psi \in \IH$
$$||\tilde{L}_{\Psi}(\psi)||^2_{\IL} \geq c||\psi||^2_{\IH} - C\int_{B_{2r} \backslash N_{\delta}} |\psi|^2 d^3x. $$
The same statement is also true if we replace $\tilde{L}_{\Psi}$ by $\tilde{L}^{\dagger}_{\Psi}$.
\end{lemma}
\begin{proof}
Recall that $\chi_{\delta}$ is the cut-off function that is $1$ on $N_{\delta}$ and $0$ on $\IR^3 \backslash N_{2 \epsilon}$. And $\chi_r$ is the cut-off function that is $1$ on $B_r$ and $0$ on $\IR^3 \backslash B_{2r}$. Thus we can divide $\psi$ into three parts:

$$\psi = \chi_{\delta} \psi + (1 - \chi_{\delta})\chi_r \psi + (1 - \chi_r)\psi.$$

Assume $c$ is a small enough constant. And assume $\delta$ is small enough and $r$ is large enough. Then from proposition \ref{proposition of local Fredholmness},

$$||\tilde{L}_{\Psi}(\chi_{\delta}\psi) ||_{\IL} \geq c ||\chi_{\delta} \psi||_{\IH}.$$

From lemma \ref{Lemma of local isomorphism at infinity},

$$||\tilde{L}_{\Psi}((1 - \chi_r)\psi)||_{\IL} \geq c ||(1 - \chi_r)\psi||_{\IH}.$$

From Weitzenbock formula, the square of the $L^2$ norm of $\tilde{L}_{\Psi}((1 - \chi_{\delta}) \chi_r \psi)$ is
$$\int_{\IR^3} |\tilde{L}_{\Psi}((1 - \chi_{\delta}) \chi_r \psi)|^2 d^3x = \int_{\IR^3} (1 - \chi_{\delta}) \chi_r <\tilde{L}_{\Psi}^{\dagger} \tilde{L}_{\Psi} ( (1 - \chi_{\delta}) \chi_r \psi), \psi> d^3 x $$
$$\geq \int_{\IR^3} |\nabla_A ((1 - \chi_{\delta}) \chi_r \psi)|^2 + |[\Phi, (1 - \chi_{\delta}) \chi_r \psi]|^2  - (|F_A| + |d_A \Phi|) |\psi|^2d^3x$$
$$ \geq c ||(1 - \chi_{\delta}) \chi_r \psi ||^2_{\IH} - C \int_{B_{2r} \backslash N_{\delta}} |\psi|^2 d^3x. $$

Finally, it can be checked directly that
$$\dfrac{c}{2} ||\psi||^2_{\IH} \leq c (||\chi_{\delta}\psi ||^2_{\IH} + ||(1 - \chi_r) \psi||^2_{\IH} + ||(1 - \chi_{\delta}) \chi_r \psi ||^2_{\IH}) - C \int_{B_{2r}\backslash N_{\delta}} |\psi|^2 d^3x.$$
Thus the statement is true. The same statement for $\tilde{L}_{\Psi}^{\dagger}$ follows by the same reason.
\end{proof}

\begin{corollary}
Both $\tilde{L}_{\Psi}$ and $\tilde{L}_{\Psi}^{\dagger}$ are closed operators from $\IH$ to $\IL$ with finite dimensional kernels.
\end{corollary}

\begin{proof}
This is because the term
$$(\int_{B_{2r} \backslash N_{\delta}} |\psi|^2 d^3x)^{\frac{1}{2}} $$
is relatively compact with respect to $||\psi ||_{\IH}$ by theorem \ref{A Rellich theorem}.
\end{proof}

To prove that $\tilde{L}_{\Psi}$ has finite dimensional cokernel, we need a lemma first

\begin{lemma}\label{lemma that equivalent cokernel in H and L}
View $\tilde{L}_{\Psi}$ as a map from $\IH$ to $\IL_{w}$. Suppose $\psi$ is an element of the cokernel of $\tilde{L}_{\Psi}$ in $\IL_w$, then $w \psi$ is an element of the kernel of $\tilde{L}_{\Psi}^{\dagger}$ in $\IH$.
\end{lemma}

\begin{proof}
 For any $\psi' \in \IL_w$, we have
$$\int_{\IR^3 \backslash K} w\tilde{L}_{\Psi}(\psi) \psi' d^3 x = 0.$$ Thus $w\psi'$ is a weak solution of $\tilde{L}^{\dagger}_{\Psi} (w\psi') = 0$. By a standard local elliptic regularity argument, since $\tilde{L}^{\dagger}$ is an elliptic operator, we know that $w \psi$ is smooth on any compact subset of $\IR^3 \backslash K$ and is in the kernel of $\tilde{L}^{\dagger}$. The task is to prove that it is in $\IH$.\\

We may assume $\delta$ is small and $r$ is large. Recall that $\chi$ is the cut-off function that is $1$ on $N_{\epsilon}$ and $0$ on $\IR^3 \backslash N_{2 \epsilon}$. And $\chi_r$ is the cut-off function that is $1$ on $B_r$ and $0$ on $\IR^3 \backslash B_{2r}$. Then we only need to show that $\chi w \psi$, $(1 - \chi) \chi_r w \psi$  and $(1 - \chi_r) w \psi$ are elements in $\IH$. \\

Since $w \psi$ is smooth on $B_{2r} \backslash N_{\delta}$, $(1 - \chi) \chi_r w \psi \in \IH$. From proposition \ref{proposition of formal adjoint and cokernel}, $\chi w \psi \in \IH$. From corollary \ref{corollary of cokernel in H near infinity}, $(1 - \chi_r) w \psi \in \IH$. Thus $w \psi \in \IH$.

\end{proof}

\begin{corollary}
The operator $\tilde{L}_{\Psi}$ has finite dimensional cokernel in $\IL_w$.
\end{corollary}

\begin{proof}
This is because each element in the cokernel of $\tilde{L}_{\Psi}$ in $\IL_w$ corresponds to an element in the kernel of $\tilde{L}^{\dagger}_{\Psi}$ in $\IH$. And we have shown that $\tilde{L}^{\dagger}_{\Psi}$ has finite dimensional kernel in $\IH$.
\end{proof}

Now we have our Fredholm theory for $\tilde{L}
_{\Psi}$.

\begin{theorem}\label{Theory of global Fredholmness}
The operator $\tilde{L}_{\Psi}$ is a Fredholm operator from $\IH$ to $\IL$.
\end{theorem}

\begin{proof}
We have shown that $\tilde{L}_{\Psi}$ has closed range in $\IL$ and has finite dimensional kernel in $\IH$. Moreover, it has finite dimensional cokernel in $\IL_w$. But $\IL$ and $\IL_w$ are equivalent. So $\tilde{L}_{\Psi}$ also has finite dimensional cokernel in $\IL$. Thus it is a Fredholm operator form $\IH$ to $\IL$.
\end{proof}

\subsection{The extra $1$-dimensional gauge transformation}\label{Subsection: The extra 1 dimensional gauge transformation}

Suppose $\Psi = A + \Phi$ is a solution to the Bogomolny equations with a knot singularity. Then there is $1$ dimensional gauge transformation given by

$$u = e^{t\Phi} $$

that is not detected by our gauge fixing condition. In fact, we assume $\Psi$ is an admissible configuration and a solution. An infinitesimal gauge transformation in this direction is give by an element $\psi_0$ in $\IH$ which can be written as  
$$\psi_0 = d_A \Phi = *F_A.$$

It can be checked directly that
$$\tilde{L}_{\Psi}(\psi_0) = *d_A\psi_0 - [\psi_0 \wedge \Phi] + *d_A *\psi_0 = *d_A d_A \Phi - [*F_A \wedge \Phi] + *d_A F_A = 0.$$

In order to study the moduli space of Bogomolny equations, we need one more extra $1$-dimensional gauge fixing condition on top of $\tilde{L}_{\Psi}$. The extra gauge fixing map send a configuration $\psi \in \IH$ to $<\psi, \psi_0>_{\IH}$, where the subscript $\IH$ indicates that the pair is given by $\IH$. One can thus extend $\tilde{L}_{\Psi}$ to a map from $\IH$ to $\IL \oplus \IR$, where the last component is given by $<\psi, \psi_0>_{\IH}$. We still denote it as $\tilde{L}_{\Psi}$ in this section when there is no ambiguity. And clearly the new version of $\tilde{L}_{\Psi}$ is still a Fredholm map but the index decreases by $1$.

\subsection{The moduli space as an analytical set}\label{Subsection: The moduli space as an analytical set}

The argument in this subsection is standard. Suppose $\Psi$ is a solution to the Bogomolny equations with knot singularity with monodromy $\gamma \neq \dfrac{1}{4}$ at the knot $K$, and has mass ``$M > 0$" and charge $k$ at infinity. Suppose for the moment that $\tilde{L}_{\Psi}$ is the operator from $\IH$ to $\IL \oplus \IR$ introduced in the last subsection. Then there is a direct corollary of the Fredholmness of $\tilde{L}_{\Psi}$:
\begin{theorem}
There exists a small open neighbourhood of $0$ in $ \IH$, namely $U$, and a real analytical map $f$ from $(\ker \tilde{L}_{\Psi}) \cap U$ to $\coker \tilde{L}_{\Psi} \subset \IL \oplus \IR$, such that solutions to the Bogomolny equation with gauge fixing conditon w.r.t. $\Psi$ in $U$ are 1-1 correspondent with $f^{-1}(0)$.
\end{theorem}
\begin{proof}
Since $\tilde{L}_{\Psi}$ is Fredholm, there exists a $C > 0$, such that for any $W \in \IL \oplus \IR$, letting $W^{\perp}$ be the projection of $W$ onto $(\coker \tilde{L}_{\Psi})^{\perp} \subset \IL \oplus \IR$, then there exists a unique $\psi \in (\ker \tilde{L}_{\Psi})^{\perp} \subset \IH$, such that $\tilde{L}_{\Psi}(\psi) = W^{\perp}$. Moreover, $\norm{\psi}_{\IH} \leq C \norm{W^{\perp}}_{\IL {\textstyle \oplus} {\scriptstyle \IR}}$.\\

It may also be assumed that $\norm{Q(\psi, \psi)}_{\IL} \leq C \norm{\psi}_{\IH}^2$ by proposition \ref{proposition of quadratic map is bounded}.\\

The proof has two steps:\\

\textbf{Step 1:} Assume $U$ is a small enough open neighbourhood of $0$ in $\IH$. For any $\psi_0 \in (\ker \tilde{L}_{\Psi}) \cap U$, there exists a unique $w(\psi_0) \in (\ker \tilde{L}_{\Psi})^{\perp}$, such that
$$\tilde{L}_{\Psi}(w(\psi_0)) + Q(\psi_0 + w(\psi_0), \psi_0 + w(\psi_0)) \in (\coker \tilde{L}_{\Psi}) \subset \IL \oplus \IR. $$

Choose any $\psi_0 \in (\ker \tilde{L}_{\Psi}) \cap U$. Consider $\psi_1, \psi_2, \cdots, \in (\ker \tilde{L}_{\Psi})^{\perp}$ such that
$$\tilde{L}_{\Psi}(\psi_1)  =  - Q(\psi_0, \psi_0)^{\perp};$$
$$\tilde{L}_{\Psi}(\psi_2)  = - 2Q(\psi_0, \psi_1)^{\perp} - Q(\psi_1, \psi_1)^{\perp};$$
$$\vdots$$
$$ \tilde{L}_{\Psi}(\psi_n) =  - 2 Q(\psi_0  + \cdots + \psi_{n-2}, \psi_{n-1})^{\perp} - Q(\psi_{n-1}, \psi_{n-1})^{\perp}.$$
$$\vdots $$
Here $(\cdot)^{\perp}$ is the projection to $(\coker \tilde{L}_{\Psi})^{\perp}$ in $\IL \oplus \IR$.\\

Since $U$ is a small neighborhood, it may be assumed that
$$C^2 \norm{\psi_0}^2_{\IH}+ \norm{\psi_0}_{\IH} \leq \dfrac{1}{2^4C^2}.$$ So $$\norm{\psi_1}_{\IH} \leq C ( \norm{Q(\psi_0, \psi_0)}_{\IL}) \leq C^2 \norm{\psi_0}_{\IH}^2 \leq \dfrac{1}{2^4C^2} ; $$
$$\norm{\psi_2}_{\IH} \leq C^2 \norm{\psi_1}_{\IH}(\norm{\psi_1}_{\IH} + 2\norm{\psi_0}_{\IH}) \leq \dfrac{1}{2}\norm{\psi_1}_{\IH}\leq \dfrac{1}{2^5 C^2};$$
$$\vdots$$
$\norm{\psi_n}_{\IH} \leq C^2\norm{\psi_{n-1}}_{\IH}(\norm{\psi_{n-1}}_{\IH} + 2\norm{\psi_{n-2}}_{\IH}  + \cdots + 2\norm{\psi_{0}}_{\IH}) \leq \dfrac{1}{2}\norm{\psi_{n-1}}_{\IH} \leq \dfrac{1}{2^{n + 3}C^2 }.$\\
Therefore letting $w(\psi_0) = \sum\limits_{n=1}^{+\infty}\psi_n$ in $\IH$ satisfies the requirement (except the uniqueness).\\

To prove the uniqueness, supposing that $\hat{w}(\psi_0)$ also satisfies the requirement, then
$$\norm{\hat{w}(\psi_0) - w(\psi_0)}_{\IH} \leq C\norm{Q(\hat{w}(\psi_0) + \psi_0, \hat{w}(\psi_0) - w(\psi_0))}_{\IL}$$
$$\leq C^2\norm{\hat{w}(\psi_0) + w(\psi_0)}_{\IH} \norm{\hat{w}(\psi_0) - w(\psi_0)}_{\IH}.$$

Since $U$ is small enough, it can be assumed that $C^2\norm{\hat{w}(\psi_0) + w(\psi_0)}_{\IH} < 1$, which implies that $\hat{w}(\psi_0) = w(\psi_0)$.\\

\textbf{Step 2:} Having $w(\psi_0)$ defined, it is obvious that, for any $\psi_0 \in (\ker \tilde{L}_{\Psi}) \cap U$, there exists an element $\psi_1$ such that $\tilde{L}_{\Psi}(\psi_0 + \psi_1) + Q(\psi_0 + \psi_1, \psi_0 + \psi_1) = 0$ if and only if $\psi_1 = w(\psi_0)$ and $\tilde{L}_{\Psi}(\psi_0 + w(\psi_0)) + Q(\psi_0 + w(\psi_0), \psi_0 + w(\psi_0)) = 0$, and if and only if the projection of $Q(\psi_0 + w(\psi_0), \psi_0 + w(\psi_0))$ onto $(\coker \tilde{L}_{\Psi})$ is $0$.\\

Let $f(\psi_0)$ be the projection of $Q(\psi_0 + w(\psi_0), \psi_0 + w(\psi_0))$ onto $(\coker \tilde{L}_{\Psi})$. Then $f$ is a map from $(\ker \tilde{L}_{\Psi}) \cap U$ to $(\coker \tilde{L}_{\Psi})$ which satisfies the requirements.
\end{proof}

\section{Existence of solutions from gluing}\label{Section: Existence of solutions from gluing}

For each fixed positive integer $k$, one may wonder whether a solution with charge $k$ and a knot singularity exists. The answer is yes when $\gamma \neq 
\dfrac{1}{4}$ for sufficiently large $M$. The solution can be constructed by a gluing method adapted from Taubes \cite{Jaffe1980VorticesTheories}. Here is a sketch of the proof:\\

We ``glue in" $k$ standard $1$-monopoles to  the model solution $\Psi_{\gamma, M, 0}$ centered far away from the knot. Then we get an approximate solution $\Psi$ to the Bogomolny equations. In another word, $V(\Psi)$ is small in certain sense. After a gauge transformation, we may assume $\Psi - \Psi_{\gamma, M, k} \in \IH$. So $\tilde{L}_{\Psi}$ is a Fredholm operator from $\IH$ to $\IL$.\\

In order to apply the implicit function theorem to get an actual solution to the Bogomolny equations, we require that the cokernel of $\tilde{L}_{\Psi}$ is $0$. This can be proved when $M$ is large enough and when the centers of the glued ``1"-monopoles are far away from the knot and from each other.\\

In this section, we always assume $\gamma \neq \dfrac{1}{4}$.

\subsection{Approximate solutions from the standard ``1"-monopole}\label{Subsection: Approximate solutions from the standard 1 monopole}

\begin{definition}
Choose a model point $P$ in $\IR^3$. Choose a Cartesian coordinates of $\IR^3$ centered at P, namely $\{w_1, w_2, w_3\}$. Let $w = M(w_1 + iw_2), t = Mw_3, r = \sqrt{|w|^2 + t^2}$. The following solution $\Psi = \Phi + A = \Phi + A_1dw_1 + A_2 dw_2 + A_3 dw_3$ is called \textbf{the 1-monopole}:
\begin{equation*} 
\left\{
\begin{array}{lr}
\Phi = \dfrac{M (r \coth r - 1)}{2r(r\cosh r - t\sinh r)} ((t\cosh r - r\sinh r)\begin{pmatrix}
-i& 0 \\
0 & i
\end{pmatrix} + \begin{pmatrix}
0& w \\
- \bar{w} &  0
\end{pmatrix}), &\\
A_3 = \dfrac{iM}{2r}(\dfrac{r \coth r - 1}{r \cosh r - t \sinh r})\begin{pmatrix}
0 & w \\
 \bar{w} & 0
\end{pmatrix},\\
A_1 + iA_2= \dfrac{M}{2(r \cosh r - t\sinh r)}((- \dfrac{ \cosh r}{r}  + \dfrac{ 1 }{ \sinh r }) \begin{pmatrix}
- w &  0\\
0 & w
\end{pmatrix}  + (\dfrac{t + r^2 - rt \coth r}{r})\begin{pmatrix}
0 &  0\\
2i&  0
\end{pmatrix}).
\end{array}
\right.
\end{equation*} 
\end{definition}

\begin{theorem}
The $1$-monopole satisfies the Bogomolny equations. Moreover,  consider the  half line $`` |w| = 0, t \geq 0"$. We call it the ``Dirac ray" $L$.  Then at a point whose distance to the Dirac ray $d \geq \dfrac{2}{M}$,

\begin{itemize}
     \item  $|\nabla_A \Psi| = |F_A| = O(\dfrac{1}{r^2})$;
     \item $\Phi =  - \dfrac{M}{2}\begin{pmatrix}
         -i & 0 \\
         0 & i
         \end{pmatrix} + O(\dfrac{1}{r})$; \item $|A| = O(\dfrac{1}{r})$;
     \item   $|[A, \Phi]| = O(e^{-r})$.
 \end{itemize}

\end{theorem}

\begin{proof}

 Suppose $$ U = \dfrac{(r+t)e^{-r}}{|w|}.$$ Then $$r\cosh r - t\sinh r  = \dfrac{|w|}{2}(U + U^{-1}), ~~~  t \cosh r - r\sinh r = \dfrac{|w|}{2} (U - U^{-1}). $$
\begin{equation*} 
\left\{
\begin{array}{lr}
\Phi = (\coth r - \dfrac{1}{r} )\dfrac{M}{ (U + U^{-1})} (\dfrac{(U - U^{-1})}{2}\begin{pmatrix}
-i& 0 \\
0 & i
\end{pmatrix} + \dfrac{1}{|w|}\begin{pmatrix}
0& w \\
- \bar{w} &  0
\end{pmatrix}), &\\ 
A_3 = iM(\coth r - \dfrac{1}{r})(\dfrac{1}{|w| (U + U^{-1})})\begin{pmatrix}
0 & w \\
 \bar{w} & 0
\end{pmatrix}, &\\
A_1 + iA_2 =  \dfrac{M}{|w| (U + U^{-1})}((- \dfrac{ \cosh r}{r}  + \dfrac{ 1 }{ \sinh r }) \begin{pmatrix}
- w &  0\\
0 & w
\end{pmatrix} - (\dfrac{-t - r^2 + rt \coth r}{r})\begin{pmatrix}
0 &  0\\
2i&  0
\end{pmatrix}).
\end{array}
\right.
\end{equation*}

Some calculations:
\begin{equation*}
\left\{
    \begin{array}{ll}
        F_{12} = \nabla_3 \Phi  = \partial_3 \Phi + [A_3, \Phi]  \\
          ~~~~ = M^2(\dfrac{t}{r^3} - \dfrac{t}{r\sinh^2 r})\dfrac{1}{ (U + U^{-1})} (\dfrac{(U - U^{-1})}{2}\begin{pmatrix}
-i& 0 \\
0 & i
\end{pmatrix} + \dfrac{1}{|w|}\begin{pmatrix}
0& w \\
- \bar{w} &  0
\end{pmatrix})\\
 ~~~~~~~~ - M^2(\coth r - \dfrac{1}{r})(\dfrac{|w|}{ r\sinh r( U + U^{-1}) })(  \dfrac{(U - U^{-1})}{2|w|}\begin{pmatrix}
0 & w \\
-\bar{w} & 0
\end{pmatrix} -  \begin{pmatrix}
-i & 0 \\
0 &  i
\end{pmatrix});\\
 F_{23} + iF_{31}  =(\nabla_1 + i\nabla_2) \Phi = (\partial_1 + i\partial_2) \Phi + [A_1 + iA_2, \Phi] \\
  ~~~~ = M^2 (\dfrac{w}{r^3} - \dfrac{w}{r\sinh^2 r})\dfrac{1}{ (U + U^{-1})} (\dfrac{(U - U^{-1})}{2}\begin{pmatrix}
-i& 0 \\
0 & i
\end{pmatrix} + \dfrac{1}{|w|}\begin{pmatrix}
0& w \\
- \bar{w} &  0
\end{pmatrix})&\\
~~~~ - (\coth r -\dfrac{1}{r})\dfrac{M^2w}{r \sinh r |w|(U+U^{-1})}(   t\begin{pmatrix}
-i& 0 \\
0 & i
\end{pmatrix} + \cosh r\begin{pmatrix}
0 & w \\
 -\bar{w} & 0 \\
\end{pmatrix} + \dfrac{r(U+U^{-1})}{|w|} \begin{pmatrix}
0 & 0 \\
\bar{w} & 0\\
\end{pmatrix}).
    \end{array}\right.
\end{equation*}

In particular, the $1$-monopole satisfies the Bogomolny equations.\\

Let $v = \dfrac{1}{ (U + U^{-1})} (\dfrac{(U - U^{-1})}{2}\begin{pmatrix}
-i& 0 \\
0 & i
\end{pmatrix} + \dfrac{1}{|w|}\begin{pmatrix}
0& z \\
- \bar{z} &  0
\end{pmatrix})$, which is a unit vector. Let $v^{\perp}$ denote any unit vector that is perpendicular to $v$. Then

\begin{equation*}
    \left\{
    \begin{array}{ll}
         \Phi = M(\coth r  - \dfrac{1}{r}) v,& \\
         F_{12} = M^2(\dfrac{t}{r^3} - \dfrac{t}{\sinh^2 r}) v - M^2(\coth r - \dfrac{1}{r})\dfrac{\rho}{r\sinh r} v^{\perp},&  \\
         F_{23} + iF_{31}  = M^2(\dfrac{z}{r^3} - \dfrac{z}{r\sinh^2 r}) v - M^2(\coth r -\dfrac{1}{r}) \dfrac{\sqrt{\rho^2 + 2t^2}}{r \sinh r} v^{\perp}.& 
    \end{array}
    \right.
\end{equation*}

Thus at a point whose distance to the Dirac ray $d \geq \dfrac{2}{M}$ so that $(r - t)^2 e^r \geq \rho$, \begin{equation*}
    \begin{array}{ll}
         \dfrac{(U - U^{-1})}{(U + U^{-1})}& = -1 + \dfrac{2\rho^2 }{\rho^2 + (r-t)^2 e^{2r}} = -1 + O(e^{-r}),
    \end{array}
\end{equation*}
then $$v = - \dfrac{1}{2}\begin{pmatrix}
         -i & 0 \\
         0 & i
         \end{pmatrix} + O(e^{-r}).$$

So all the properties listed in the theorem follow directly.

 \end{proof}

\textbf{Remark}
The 1-monopole defined here may be in an unfamliliar format to experts. In fact, under a gauge transformation, the 1-monopole is gauge equivalent to the following solution to the Bogomolny equation, which is better known (see for example Page 104 in \cite{Jaffe1980VorticesTheories}):\\
\begin{equation*}
\left\{
    \begin{array}{ll}
         \Phi ~~ = ~ (\dfrac{1}{r\tanh r} - \dfrac{1}{r^2}) (w_1 \sigma_1 + w_2 \sigma_2 + w_3 \sigma_3),  \\
         A = (\dfrac{1}{r \sinh r} - \dfrac{1}{r^2}) ((w_2\sigma_3 - w_3 \sigma_2) dw_1 + (w_3\sigma_1 - w_1 \sigma_3) dw_2 + (w_1\sigma_2 - w_2 \sigma_1) dw_3).
    \end{array}
    \right.
\end{equation*}
where $\sigma_1 = \begin{pmatrix}
 -i & 0\\
 0 & i
\end{pmatrix}, \sigma_2 = \begin{pmatrix}
 0 & -1\\
 1 & 0
\end{pmatrix}, \sigma_3 = \begin{pmatrix}
 0 & i\\
 i & 0
\end{pmatrix}$.\\

However, our format makes it much easier to create an approximate solution from gluing.\\

\textbf{Remark:} Alternatively, like what Taubes did in \cite{Jaffe1980VorticesTheories}, one may use a simpler format of the 1-monopole but has to work in a singular gauge to do the gluing. This is essentially equivalent to the approach here.

\begin{definition}\label{definition of adimissible approximate solutions}
If a configuraiton $\Psi = A + \Phi$ on $\IR^3 \backslash K$ satisfies the following constraints:
\begin{itemize}
    \item There exists a constant $M > 0$ such that $\lim\limits_{r \rightarrow \infty} |\Phi| = M$;
    \item For some small enough $\delta$, such that $\Psi = \Psi_{\gamma, M, 0}$ on $N_{100 \delta}$ (recall that $M$ is used to define the $\Phi_{\gamma, M, 0} = M \sigma$ portion in $\Psi_{\gamma, M, 0}$). In particular $V(\Psi) = *F_A = \nabla_A \Psi = 0$ on $N_{100 \delta}$;
    \item For some small enough $\mu > 0$ to be determined, the $L^2(\IR^3)$ norm of $V(\Psi)$ is less than $\mu$;
    \item After a gauge transformation supported on $\IR^3 \backslash B_r$, where $r$ is large enough, we may assume $\dfrac{\Psi}{|\Psi|} = e_0$, where $e_0$ is defined in subsection \ref{Subsection: The model configuration}.
\end{itemize}
Then $\Psi$ is called an \textbf{approximate solution}.
\end{definition}

Suppose $\Psi_1, \Psi_2, \Psi_3, \cdots \Psi_n$ are 1-monopole solutions centered at $P_1, P_2, \cdots, P_n$ and have Dirac rays $L_1, L_2, \cdots, L_n$ respectively. Suppose $M$ is large enough. Suppose $d > 0$ is a large enough constant and $R \gg d$ is even larger. Let $V_i(R)$ be the open set of points whose distance to the Dirac ray $L_i$ are less than $R$. Suppose moreover any two of the corresponded open neighbourhoods $V_1(R), V_2(R), \cdots, V_n(R)$ together with $B_{R}(O)$ (the ball of radius $R$ centered at the origin) do not intersect. Let $V_0(R)$ denote $B_R(O)$. Let $r$ be the distance function to the origin $O$. Then

\begin{proposition}
Let $\Psi = \sum\limits_{j = 0}^n \tilde{\chi}_j \Psi_j$, where each $\tilde{\chi}_j$  is a suitable cut-off function that is $1$ on $V_j(d)$ and $|\nabla \tilde{\chi}_j| \leq \dfrac{1}{dr} $ when $j \geq 1$. We may assume that the supports of $\tilde{\chi}_1, \tilde{\chi}_2, \cdots, \tilde{\chi}_n$ do not overlap. Moreover, let $ \tilde{\chi}_0 = 1 -\sum\limits_{j=1}^n \tilde{\chi}_j$. Given that  $d$ is large enough and $R \gg d$, then $\Psi$ is an approximate solution.
\end{proposition}

\begin{proof}
All the other requirements in definition \ref{definition of adimissible approximate solutions} except the third bullet ($V(\Psi)$ has a small $L^2$ norm) are straightforward. For the third bullet, consider
$$V(\Psi) = V(\sum\limits_{j = 0}^n \tilde{\chi}_j \Psi_j). $$
Note that for each $j$, $V(\Psi_j) = 0$. So we only need to estimate $|V(\Psi)|$ at a point $P$ where certain $\tilde{\chi}_j$ is neither $0$ nor $1$. We may assume $j = 1$ here. At the point $P$,

$$ \Psi = \tilde{\chi_1} \Psi_1 + (1 - \chi_1) \Psi_0.$$

One sees that $|V(\Psi)|$ is bounded above by a constant multiplication of
$$|\nabla \tilde{\chi}_1| |\Psi_0 - \Psi_1| + |Q(\Psi_0 - \Psi_1, \Psi_0)| + |\Psi_0 - \Psi_1|^2.$$

Assuming $R \gg d$ are both large, we have $|\nabla \tilde{\chi}_1| \leq \dfrac{c_0}{r}$, ~~$|\Psi_0 - \Psi_1| \leq \dfrac{c_0}{r}$, where $c_0$ is an arbitrarily small constant (depending on $R$ and $d$). Moreover, $|Q(\Psi_0 - \Psi_1, \Psi_0)| \leq c_0e^{-r}.$ So in fact 
$$|V(\Psi)| \leq \dfrac{3c_0}{r^2}, $$
where $c_0$ can be arbitrarily small. This is enough to say that the $L^2$ norm of $|V(\Psi)|$ is small.
\end{proof}

\subsection{$\tilde{L}_{\Psi}$ has $\{0\}$ cokernel}

\begin{proposition}\label{proposition that glued solution has 0 cokernel}
Let $\Psi$ be the approximate solution defined in subsection \ref{Subsection: Approximate solutions from the standard 1 monopole}. Then if $M$ is large enough, and if $\mu$ is small enough, then the cokernel of $\tilde{L}_{\Psi}$ is $\{0\}$.\\
\end{proposition}

The proof of this proposition is cumbersome, which occupies this subsection.\\

First of all, we can decompose $\psi$ into two parts, the part parallel to $\sigma$ and the part orthogonal to $\sigma$. To be more precise, let 
$$\psi = \psi^{\parallel} + \psi^{\perp}, ~~ \text{where} ~~ \psi^{\perp} = - [\sigma, [\sigma, \psi]], ~~ \psi^{\parallel} = <\psi, \sigma> \sigma.$$

Then for any fixed $C > 0$, there exists a $M_0 > 0$ and a $\epsilon_0 < \epsilon$, such that if $M \geq M_0$, then  

$$\epsilon_0 \int_{\IR^3 \backslash N_{\epsilon_0}} (|\nabla_A \psi|^2 + |[\Phi, \psi]|^2) d^3 x \geq 2C \int_{\partial N_{\epsilon_0}} |\psi^{\perp}|^2 d\theta ds.$$

The proof of the above statement is in subsection \ref{Subsection: Proof of some technical lemmas} (lemma \ref{technical lemma 2}). Given the above statement, then proposition \ref{proposition that glued solution has 0 cokernel} is a corollary of the following proposition:

\begin{proposition}\label{proposition middle tedious estimate}
Suppose $\psi \in \IH$ is a smooth configuration. Then there is a constant $C$ which doesn't depend on $\epsilon_0 < \epsilon$ such that,

$$C\int_{N_{\epsilon_0} \backslash K}\rho |\tilde{L}^{\dagger}_{\Psi} \psi|^2 d^3 x + C \int_{\partial N_{\epsilon_0}} |\psi^{\perp}|^2 d\theta ds + C \epsilon_0 \int_{\IR^3 \backslash N_{\epsilon_0}} |\tilde{L}^{\dagger}\psi|^2 d^3 x$$ $$\geq \epsilon_0 \int_{\IR^3 \backslash N_{\epsilon_0}} (|\nabla_A\psi|^2 + |[\Phi, \psi]|^2) d^3x.$$
\end{proposition}

\paragraph{Proof of proposition \ref{proposition that glued solution has 0 cokernel} using proposition \ref{proposition middle tedious estimate}}

If $\psi \in \IH$ is an element in the cokernel of $\tilde{L}^{\dagger}$, then by a standard elliptic regularity argument, $\psi$ is smooth on any compact subset of $\IR^3 \backslash K$. Moreover, by proposition \ref{proposition middle tedious estimate}
$$C \int_{\partial N_{\epsilon_0}} |\psi^{\perp}|^2 d\theta ds \geq \epsilon_0 \int_{\IR^3 \backslash N_{\epsilon_0}} (|\nabla_A \psi|^2 + |[\Phi, \psi]|^2) d^3x.$$

Since $M$ is large enough, we may assume $ \epsilon_0$ satisfies  
$$C \int_{\partial N_{\epsilon_0}} |\psi^{\perp}|^2 d\theta ds \leq \dfrac{\epsilon_0}{2} \int_{\IR^3 \backslash N_{\epsilon_0}} (|\nabla_A \psi|^2 + |[\Phi, \psi]|^2) d^3x. $$ 
So the only possibility is
$$\int_{\partial N_{\epsilon_0}} |\psi^{\perp}|^2 d\theta ds =  \int_{\IR^3 \backslash N_{\epsilon_0}} (|\nabla_A \psi|^2 + |[\Phi, \psi]|^2) d^3x = 0.$$

By theorem \ref{A Rellich theorem}, $$\int_{C} |\psi|^2 = 0$$
on any compact subset of $\IR^3 \backslash N_{\epsilon_0}$, which implies that $\psi = 0$ on $\IR^3 \backslash N_{\epsilon_0}$. So $\psi$ is supported in $N_{\epsilon_0}$. However, we may assume that $\epsilon_0$ is small enough such that proposition \ref{proposition of local Fredholmness} holds, which implies
$$||\psi ||_{\IH} \leq C ||\tilde{L}_{\Psi}^{\dagger} (\psi)||_{\IL} = 0.$$
So $\psi = 0$. Thus the kernel of $\tilde{L}_{\Psi}^{\dagger}$ in $\IH$ is $\{0\}$. By lemma \ref{lemma that equivalent cokernel in H and L}, $\tilde{L}_{\Psi}$ has a $\{0\}$ cokernel in $\IL$.\\

The remaining part of this subsection proves proposition \ref{proposition middle tedious estimate}.

\paragraph{Proof of proposition \ref{proposition middle tedious estimate}}

We need to introduce some notations here. Let $\rho d\Omega$ be the area form for $\partial N_{\epsilon_0}$. Then the following three measures on $N_{\epsilon_0} \backslash K$ are different, but  equivalent (bounded by each other):
$$\rho d\Omega, ~~ \rho d\rho d\theta, ~~ d^3x.$$

In particular, we may assume that
$$C \int_{N_{\epsilon_0} \backslash K} \rho |\tilde{L}^{\dagger}_{\Psi} \psi|^2 d^3 x \geq \int_{N_{\epsilon_0} \backslash K} \rho|\tilde{L}^{\dagger}_{\Psi} \psi|^2 \rho d\Omega. $$

Moreover, it suffices to prove the following modification of proposition \ref{proposition middle tedious estimate}:

\begin{equation}\label{modificaiton of the proposition}
    \begin{array}{ll}
         & ~~~~ \displaystyle \int_{N_{\epsilon_0} \backslash K}\rho |\tilde{L}^{\dagger}_{\Psi} \psi|^2 \rho d\Omega + C \int_{\partial N_{\epsilon_0}} |\psi^{\perp}|^2 d\Omega + 2 \epsilon_0 \int_{\IR^3 \backslash N_{\epsilon_0}} |\tilde{L}^{\dagger}\psi|^2 d^3 x \\
         & \geq \displaystyle \dfrac{1}{C} \int_{N_{\epsilon_0} \backslash K} \rho^2 |\nabla_{A_{\rho}} \psi|^2 d\rho d\Omega + \epsilon_0 \int_{\IR^3 \backslash N_{\epsilon_0}} (|\nabla_A\psi|^2 + |[\Phi, \psi]|^2) d^3x.
    \end{array}
\end{equation}

We need 5 steps to prove (\ref{modificaiton of the proposition}):\\

\textbf{Step 1:} We compute 
$$\epsilon_0 \int_{\IR^3 \backslash N_{\epsilon_0 }} |\tilde{L}^{\dagger}_{\Psi} \psi|^2 d^3x. $$
\begin{equation*}
    \begin{array}{ll}
         ~ & ~~~~\displaystyle \epsilon_0 \int_{\IR^3 \backslash N_{\epsilon_0 }} |\tilde{L}^{\dagger}_{\Psi} \psi|^2 d^3x =  \epsilon_0 \int_{\IR^3 \backslash N_{\epsilon_0 }} |\sum\limits_{j=1}^3 \tau_j \nabla_{A_j} \psi - [ \Phi, \psi]|^2 d^3x\\
         & = \epsilon_0  \displaystyle \int_{\IR^3 \backslash N_{\epsilon_0 }} (|\sum\limits_{j} \tau_j \nabla_{A_j} \psi|^2 - 2 \sum\limits_{j} <\tau_j \nabla_{A_j} \psi, [\Phi, \psi]>  + |[\Phi, \psi]|^2) d^3 x \\
         & = \displaystyle \epsilon_0  \int_{\IR^3 \backslash N_{\epsilon_0 }} (|\nabla_A \psi|^2 + \sum\limits_{j} (<\psi, \tau_j [(*F_A - d_A\Phi)_j, \psi]>)  + |[\Phi, \psi]|^2) d^3 x \\
         & ~~~~ + \displaystyle \epsilon_0  \int_{\partial N_{\epsilon_0 }}  (- \dfrac{1}{\epsilon_0 }< \tau_{\rho} \psi, \tau_{\theta} \nabla_{A_{\theta}} \psi> - < \tau_{\rho} \psi, \tau_s \nabla_{A_s} \psi> +  < \tau_{\rho} \psi, [\Phi, \psi]> ) \rho d\Omega.
    \end{array}
\end{equation*}

Here $\nabla_{A_{\theta}}$ is the covariant version of $\partial_{\theta}$. However, $\nabla_{A_s}$ is the covariant derivative in the direction that is orthogonal to $\partial_{\theta}$ and $\partial_{\rho}$. And $\tau_{\rho}, \tau_{\theta}, \tau_s$ is an orthonormal basis of the $1$-forms.  \\

\textbf{Step 2:} We compute 
$$\int_{N_{\epsilon} \backslash K} \rho^2 |\tilde{L}^{\dagger}_{\Psi} \psi|^2 d\Omega. $$

\begin{equation*}
    \begin{array}{ll}
         ~ & ~~~~\displaystyle \int_{N_{\epsilon} \backslash K} \rho^2 |\tilde{L}^{\dagger}_{\Psi} \psi|^2 d\Omega =  \int_{ N_{\epsilon} \backslash K} \rho^2|D_A \psi - [\Phi, \psi]|^2 d\Omega\\
         & \geq \displaystyle \int_{N_{\epsilon} \backslash K} \rho^2 (\sum\limits_{j}|\nabla_{A_j} \psi|^2  + |[\Phi, \psi]|^2 +  \sum\limits_{j} (<\psi, \tau_j [(*F_A - d_A\Phi)_j, \psi]>))d\Omega \\
         & ~~~~ + \displaystyle \int_{ N_{\epsilon} \backslash K}  <\psi, \dfrac{\tau_s}{\rho} \nabla_{A_{\theta}} \psi - \tau_{\theta}\nabla_{A_s}\psi + \tau_{\rho} [\Phi, \psi]> \rho d\Omega\\
         & ~~~~ -  \displaystyle \epsilon_0  \int_{\partial N_{\epsilon_0 }}  (- \dfrac{1}{\epsilon_0 }< \tau_{\rho} \psi, \tau_{\theta} \nabla_{A_{\theta}} \psi> - < \tau_{\rho} \psi, \tau_s \nabla_{A_s} \psi> +  < \tau_{\rho} \psi, [\Phi, \psi]> ) \rho d\Omega \\
         & ~~~~ - \displaystyle \liminf\limits_{\rho \rightarrow 0} \rho \int_{\partial N_{\rho}}|\psi|(\dfrac{1}{\rho} |\nabla_{A_{\theta}} \psi| + |\nabla_{A_s}\psi| +|[\Phi, \psi]|)\rho d\Omega.
         \end{array}
         \end{equation*}

The last term above vanishes, because
\begin{equation*}
    \begin{array}{ll}
        & ~~~~ \displaystyle \liminf\limits_{\rho \rightarrow 0} \rho \int_{\partial N_{\rho}}|\psi|(\dfrac{1}{\rho} |\nabla_{A_{\theta}} \psi| + |\nabla_{A_s}\psi| +|[\Phi, \psi]|)\rho d\Omega \\
         & \leq  \displaystyle \liminf\limits_{\rho \rightarrow 0}  \int_{\partial N_{\rho}}\rho (|\psi|^2 \cdot |\Phi| + |\psi||\nabla_A\psi|) \rho d\Omega  =  \displaystyle \liminf\limits_{\rho \rightarrow 0}  \int_{\partial N_{\rho}}\rho ( |\psi||\nabla_A\psi|) \rho d\Omega \\
         & \leq \displaystyle \liminf\limits_{\rho \rightarrow 0}  (\int_{\partial N_{\rho}} ( |\psi|^2) \rho d\Omega)^{\frac{1}{2}}(\int_{\partial N_{\rho}}\rho^3 (|\nabla_A\psi|^2) \rho d\Omega)^{\frac{1}{2}}.
    \end{array}
\end{equation*}
Here by lemma \ref{boundary convergent lemma} $$\displaystyle \lim\limits_{\rho \rightarrow 0}  \int_{\partial N_{\rho}} ( |\psi|^2) \rho d\Omega = 0. $$
$$\text{If} ~~~\displaystyle \liminf\limits_{\rho \rightarrow 0}\int_{\partial N_{\rho}}\rho^3 (|\nabla_A\psi|^2) \rho d\Omega > 0, ~\text{then}~~ \int_{N_{\epsilon}\backslash K}\rho |\nabla_A\psi|^2 \rho d\rho d\Omega = \infty,$$
which contradicts the fact that $\psi \in \IH.$\\

In addition, note that $*F_A = d_A\Phi = 0$ on $N_{\epsilon_0} \backslash K$. Thus

\begin{equation*}
    \begin{array}{ll}
         ~ & ~~~~\displaystyle \int_{N_{\epsilon_0} \backslash K} \rho |\tilde{L}^{\dagger}_{\Psi} \psi|^2 \rho d\Omega = \displaystyle \int_{N_{\epsilon_0} \backslash K} (\sum\limits_{j}|\nabla_{A_j} \psi|^2  + |[\Phi, \psi]|^2)\rho^2 d\Omega\\
         & ~~~~~~~~ + \displaystyle \int_{ N_{\epsilon_0} \backslash K}  <\psi, \dfrac{\tau_s}{\rho} \nabla_{A_{\theta}} \psi - \tau_{\theta}\nabla_{A_s}\psi + \tau_{\rho} [\Phi, \psi]> \rho d\Omega\\
         & ~~~~~~~~ -  \displaystyle \epsilon_0  \int_{\partial N_{\epsilon_0 }}  (- \dfrac{1}{\epsilon_0 }< \tau_{\rho} \psi, \tau_{\theta} \nabla_{A_{\theta}} \psi> - < \tau_{\rho} \psi, \tau_s \nabla_{A_s} \psi> +  < \tau_{\rho} \psi, [\Phi, \psi]> ) \rho d\Omega.
         \end{array}
         \end{equation*}
Combine the above expression with the one in step 1, we get

\begin{equation*}
    \begin{array}{ll}
         & ~~~~ \displaystyle \int_{N_{\epsilon_0} \backslash K} \rho |\tilde{L}^{\dagger}_{\Psi} \psi|^2 \rho d\Omega + \epsilon_0 \int_{\IR^3 \backslash N_{\epsilon}} |\tilde{L}^{\dagger}_{\Psi}\psi|^2 d^3x  \\
         & = \displaystyle \int_{N_{\epsilon_0 \backslash K}}(|\nabla_A \psi|^2 + |[\Phi, \psi]|^2) \rho^2 d\Omega + \displaystyle \int_{ N_{\epsilon_0} \backslash K}  <\psi, \dfrac{\tau_s}{\rho} \nabla_{A_{\theta}} \psi - \tau_{\theta}\nabla_{A_s}\psi + \tau_{\rho} [\Phi, \psi]> \rho d\Omega\\
         & ~~~~ + \displaystyle  \epsilon_0  \int_{\IR^3 \backslash N_{\epsilon_0 }} (|\nabla_A \psi|^2 + \sum\limits_{j} (<\psi, \tau_j [(*F_A - d_A\Phi)_j, \psi]>)  + |[\Phi, \psi]|^2) d^3 x.
    \end{array}
\end{equation*}

And we have a more general version
\begin{equation*}
    \begin{array}{ll}
         & ~~~~\displaystyle \int_{N_{\epsilon_0} \backslash K} \rho |\tilde{L}^{\dagger}_{\Psi} \psi|^2 \rho d\Omega + \epsilon_0 \int_{\IR^3 \backslash N_{\epsilon}} |\tilde{L}^{\dagger}_{\Psi}\psi|^2 d^3x  \geq  \displaystyle \int_{N_{\epsilon_0} \backslash K} \rho |\tilde{L}^{\dagger}_{\Psi} \psi|^2 \rho d\Omega + t \epsilon_0 \int_{\IR^3 \backslash N_{\epsilon}} |\tilde{L}^{\dagger}_{\Psi}\psi|^2 d^3x  \\
         & = \displaystyle \int_{N_{\epsilon_0 \backslash K}}(|\nabla_A \psi|^2 + |[\Phi, \psi]|^2) \rho^2 d\Omega + \displaystyle \int_{ N_{\epsilon_0} \backslash K}  <\psi, \dfrac{\tau_s}{\rho} \nabla_{A_{\theta}} \psi - \tau_{\theta}\nabla_{A_s}\psi + \tau_{\rho} [\Phi, \psi]> \rho d\Omega\\
         & ~~~~ + \displaystyle t \epsilon_0  \int_{\IR^3 \backslash N_{\epsilon_0 }} (|\nabla_A \psi|^2 + \sum\limits_{j} (<\psi, \tau_j [(*F_A - d_A\Phi)_j, \psi]>)  + |[\Phi, \psi]|^2) d^3 x \\
         & ~~~~ - (1 - t)\displaystyle \epsilon_0  \int_{\partial N_{\epsilon_0 }}  (- \dfrac{1}{\epsilon_0 }< \tau_{\rho} \psi, \tau_{\theta} \nabla_{A_{\theta}} \psi> - < \tau_{\rho} \psi, \tau_s \nabla_{A_s} \psi> +  < \tau_{\rho} \psi, [\Phi, \psi]> ) \rho d\Omega,
    \end{array}
\end{equation*}
where $t$ is any constant such that $\dfrac{1}{4} \leq t \leq 1$.\\

\textbf{Step 3:} One can check directly that
$$|<\psi, \tau_j[*F_A - d_A\Phi_j, \psi]>| \leq |V(\Psi)| |Q(\psi, \psi)|,$$
where $Q$ is the quadratic form defined in section \ref{Section: Introduction and preliminary set-ups}. Since $\mu$ is small enough, and from subsection \ref{Subsection: The quadratic operator Q} we know that the $L^2(\IR^3 \backslash N_{\epsilon_0})$ norm of $|Q(\psi, \psi)|$ is bounded above by
$$C \int_{\IR^3\backslash N_{\epsilon_0}} (|\nabla_A \psi|^2 + |[\Phi, \psi]|^2) d^3x, $$
we may assume that

$$\displaystyle   \int_{\IR^3 \backslash N_{\epsilon_0 }} (|\nabla_A \psi|^2 + 2\sum\limits_{j} (<\psi, \tau_j [(*F_A - d_A\Phi)_j, \psi]>)  + |[\Phi, \psi]|^2) d^3 x \geq 0. $$

Thus we have

\begin{equation*}
    \begin{array}{ll}
         & ~~~~ \displaystyle \int_{N_{\epsilon_0} \backslash K} \rho |\tilde{L}^{\dagger}_{\Psi} \psi|^2 \rho d\Omega + \epsilon_0 \int_{\IR^3 \backslash N_{\epsilon}} |\tilde{L}^{\dagger}_{\Psi}\psi|^2 d^3x  \\
         & \geq \displaystyle \int_{N_{\epsilon_0 \backslash K}}(|\nabla_A \psi|^2 + |[\Phi, \psi]|^2) \rho^2 d\Omega + \displaystyle \int_{ N_{\epsilon_0} \backslash K}  <\psi, \dfrac{\tau_s}{\rho} \nabla_{A_{\theta}} \psi - \tau_{\theta}\nabla_{A_s}\psi + \tau_{\rho} [\Phi, \psi]> \rho d\Omega\\
         & ~~~~ + \displaystyle \dfrac{1}{8} \epsilon_0  \int_{\IR^3 \backslash N_{\epsilon_0 }} (|\nabla_A \psi|^2 + |[\Phi, \psi]|^2) d^3 x \\
         & ~~~~ - \displaystyle (1 - t)\epsilon_0  \int_{\partial N_{\epsilon_0 }}  (- \dfrac{1}{\epsilon_0 }< \tau_{\rho} \psi, \tau_{\theta} \nabla_{A_{\theta}} \psi> - < \tau_{\rho} \psi, \tau_s \nabla_{A_s} \psi> +  < \tau_{\rho} \psi, [\Phi, \psi]> ) \rho d\Omega.
    \end{array}
\end{equation*}

So far, in order to prove (\ref{modificaiton of the proposition}), it suffices to prove that for any $\psi \in \IH$ such that $\tilde{L}_{\Psi}^{\dagger} \psi = 0$, for some $\dfrac{1}{4} \leq t \leq 1$,
\begin{equation}\label{modification of the proposition 2}
\begin{array}{ll}
     &  \displaystyle C  \int_{\partial N_{\epsilon_0}} |\psi^{\perp}|^2 \rho d\Omega + \int_{N_{\epsilon_0} \backslash K} ((|\nabla_{A} \psi|^2 + |[\Phi, \psi]|^2) \rho^2 +  <\psi, \tau_s \nabla_{A_{\theta}} \psi - \rho \tau_{\theta}\nabla_{A_s}\psi + \rho \tau_{\rho} [\Phi, \psi]>) d\Omega\\
         & ~~~~ - \displaystyle   (1 - t)\int_{\partial N_{\epsilon_0 }}  (- \epsilon_0 < \tau_{\rho} \psi, \tau_{\theta} \nabla_{A_{\theta}} \psi> - \epsilon_0^2 < \tau_{\rho} \psi, \tau_s \nabla_{A_s} \psi> + \epsilon_0^2 < \tau_{\rho} \psi, [\Phi, \psi]> )  d\Omega \geq 0.
\end{array}
\end{equation}

We call the above inequality the ``$t$" version of (\ref{modification of the proposition 2}). It is not necessary to prove it for all $t$. We only need to prove it for some choice of $t$ in each situation.\\

\textbf{Step 4:} Since $\Psi = \Psi_{\gamma, M, 0}$ on $N_{\epsilon_0} \backslash K$, both $[\Phi, \cdot]$ and $\nabla_A$  preserves the decomposition $\psi = \psi^{\parallel} + \psi^{\perp}$. Thus we may prove (\ref{modification of the proposition 2}) in two cases $\psi = \psi^{\parallel}$ and $\psi = \psi^{\perp}$ respectively. This step assumes $\psi = \psi^{\parallel}$.\\

Note that $\Phi = M \sigma$ on $N_{\epsilon_0} \backslash K$. So $[\Phi, \psi^{\parallel}] = 0$. And the $t = 1$ version  (\ref{modification of the proposition 2}) becomes

$$ \int_{N_{\epsilon_0} \backslash K} ((|\nabla_{A} \psi|^2) \rho^2 + ( <\psi, \tau_s \nabla_{A_{\theta}} \psi - \rho \tau_{\theta}\nabla_{A_s}\psi>)) d\Omega \geq 0.$$

The $t = \dfrac{1}{2}$ version of (\ref{modification of the proposition 2}) becomes

\begin{equation}
\begin{array}{ll}
     &  \displaystyle  \int_{N_{\epsilon_0} \backslash K} ((|\nabla_{A} \psi|^2) \rho^2 + ( <\psi, \tau_s \nabla_{A_{\theta}} \psi - \rho \tau_{\theta}\nabla_{A_s}\psi>)) d\Omega\\
         & ~~~~ - \dfrac{1}{2} \displaystyle   \int_{\partial N_{\epsilon_0 }}  (- \epsilon_0 < \tau_{\rho} \psi, \tau_{\theta} \nabla_{A_{\theta}} \psi> - \epsilon_0^2 < \tau_{\rho} \psi, \tau_s \nabla_{A_s} \psi> )  d\Omega \geq 0.
\end{array}
\end{equation}

Note that since $F_A = 0$, $\nabla_{A_{\rho}}, \nabla_{A_{\theta}}, \nabla_{A_s}$ commutes. Thus like what we did in subsection \ref{Subsection: Compare L with the operator N}, we may consider a Fourier component of $\psi$ by assuming $\psi = \begin{pmatrix}
 u(\rho) \\ v(\rho)
\end{pmatrix} e^{- im\theta - iks}$, where $m$ is an integer, $k$ takes values in a discrete subset of $\IR$.  Then we want to prove the larger one of the following two expressions is non-negative:
$$\int_0^{\epsilon} \rho^2(|\partial_{\rho}u|^2 + |\partial_{\rho} v|^2) + (m^2 + m + k^2\rho^2) (|u|^2 + |v|^2) + 2\rho <iku, v> d\rho$$
and
$$\int_0^{\epsilon} \rho^2(|\partial_{\rho}u|^2 + |\partial_{\rho} v|^2) + (m^2 + m + k^2\rho^2) (|u|^2 + |v|^2) + 2\rho <iku, v> d\rho$$
$$ - \dfrac{1}{2} m\epsilon_0(|u(\epsilon_0)|^2 + |v(\epsilon_0)|^2) - \epsilon_0^2 <iku, v>. $$

If $m^2 + m \geq 1$, then the first expression above is greater or equal to

$$\int_0^{\epsilon}  (1 + k^2\rho^2) (|u|^2 + |v|^2) + 2\rho <iku, v> d\rho \geq 0.$$

If $m^2 + m = 0$, then $m = 0$ or $m = -1$. In either case, the secon expression is not less than
$$\int_0^{\epsilon} \rho^2(|\partial_{\rho}u|^2 + |\partial_{\rho} v|^2) +  k^2\rho^2 (|u|^2 + |v|^2) + 2\rho <iku, v> d\rho - \epsilon_0^2 <iku, v> $$
$$= \int_0^{\epsilon} \rho^2 (|\partial_{\rho}u|^2 + |\partial_{\rho}v|^2) +  k^2 \rho^2 (|u|^2 + |v|^2) -\rho^2 \partial_{\rho}(<iku, v>) d\rho \geq 0.$$

\textbf{Step 5} In this step we assume $\psi = \psi^{\perp}$. We still assume  $\psi = \begin{pmatrix}
 u(\rho) \\ v(\rho)
\end{pmatrix} e^{- im\theta - iks}$, where $m$ is an integer. Like what we did in subsection \ref{Subsection: Compare L with the operator N}, we assume $\lambda = m + 2\gamma j$.  Then we want to prove the following expression is non-negative for some $\dfrac{1}{4} \leq t \leq 1$.
$$\int_0^{\epsilon} \rho^2(|\partial_{\rho}u|^2 + |\partial_{\rho} v|^2) + ({\lambda}^2 + \lambda + (k^2 + 4M^2)\rho^2) (|u|^2 + |v|^2) + 2\rho <iku \pm 2Mu, v> d\rho$$
$$ (C - (1 - t)\lambda) \epsilon_0 (|u(\epsilon_0)|^2 + |v(\epsilon_0)|^2) - 2(1 - t)\epsilon_0^2 <iku \pm 2Mu, v>. $$
If $\lambda^2 + \lambda \geq 0$, then the $t = 1$ version of the above expression is no less than
$$\int_0^{\epsilon} \rho^2(|\partial_{\rho}u|^2 + |\partial_{\rho} v|^2) + ( k^2 + 4M^2)\rho^2 (|u|^2 + |v|^2) + 2\rho <(ik \pm 2M)u, v> d\rho + C \epsilon_0 (|u(\epsilon)|^2 + |v(\epsilon)|^2).$$
By lemma \ref{lemma of a preliminary inequality on N_epsilon}, and assuming $C$ is large enough, then
$$\int_0^{\epsilon} \rho^2(|\partial_{\rho}u|^2 d\rho + C\epsilon_0 (|u(\epsilon)|^2 + |v(\epsilon)|^2) \geq \dfrac{1}{4} \int_0^{\epsilon} (|u|^2 + |v|^2) d\rho.$$
So the first expression is no less than
$$\int_0^{\epsilon} (\dfrac{1}{4} + (k^2 + 4M^2) \rho^2)(|u|^2 + |v|^2) + 2\rho <(ik \pm 2M)u, v> d\rho \geq 0.$$

If $\lambda^2 + \lambda < 0$, then either $\lambda = - 2\gamma$ or $\lambda = 2\gamma - 1$. In either case, assuming $\alpha = \min\{2\gamma, 1 - 2\gamma \}$, then the $t$ version of the expression is no less than
$$\int_0^{\epsilon} \rho^2(|\partial_{\rho}u|^2 + |\partial_{\rho} v|^2) + (- \dfrac{1}{4} + (k^2 + 4M^2)\rho^2) (|u|^2 + |v|^2) + 2\rho <(ik \pm 2M)u, v> d\rho$$
$$ + C \epsilon_0 (|u(\epsilon_0)|^2 + |v(\epsilon_0)|^2) - 2(1 - t)\epsilon_0^2 <(ik \pm 2M)u, v>(\epsilon_0) $$
$$ = \int_0^{\epsilon} \rho^2(|\partial_{\rho}u|^2 + |\partial_{\rho} v|^2) + (- \dfrac{1}{4} + (k^2 + 4M^2)\rho^2) (|u|^2 + |v|^2) + \rho^2 \partial_{\rho} (<(ik \pm 2M)u, v>) d\rho$$
$$ + C \epsilon_0 (|u(\epsilon_0)|^2 + |v(\epsilon_0)|^2) - (1 - 2t)\epsilon_0^2 <(ik \pm 2M)u, v>(\epsilon_0) $$
$$ \geq \int_0^{\epsilon} (\rho \partial_{\rho}u + \dfrac{1}{2} u + (-ik \pm 2M)\rho v)^2 + (\rho \partial_{\rho}v + \dfrac{1}{2} v + (ik \pm 2M)\rho u)^2  d\rho$$
$$ + \int_0^{\epsilon}-  \dfrac{1}{2}(|u|^2 + |v|^2) - \dfrac{1}{2}\rho \partial_{\rho}(|u|^2 + |v|^2)  - \partial_{\rho}(\rho^2 <u, (-ik \pm 2M)v>)  d\rho   $$
$$ + C \epsilon_0 (|u(\epsilon_0)|^2 + |v(\epsilon_0)|^2) - (1 - 2t)\epsilon_0^2 <(ik \pm 2M)u, v> (\epsilon_0)$$
$$\geq (C - \dfrac{1}{2}) \epsilon_0 (|u(\epsilon_0)|^2 + |v(\epsilon_0)|^2) - (2 - 2t) \epsilon_0^2 <(ik \pm 2M)u(\epsilon_0), v(\epsilon_0)>.$$

We may assume $t = 1$ and $C \geq \dfrac{1}{2}$, then the above expression is non-negative. And the proof is finished.

\subsection{Proof of some technical lemmas}\label{Subsection: Proof of some technical lemmas}

We prove some technical lemmas used in the last subsection.\\

Suppose $\Psi = A + \Phi = \Psi_{\gamma, M, 0}$ on $N_{\epsilon} \backslash K$. For each configuration $\psi$, it can be decomposed into the parts that parallel to $\sigma$ and perpendicular to $\sigma$. To be more precise,
$$\psi = \psi^{\parallel} + \psi^{\perp}, $$
where $\psi^{\parallel} = <\psi, \sigma> \sigma$, $\psi^{\perp} = [\sigma, [\psi, \sigma]].$ Note that $\nabla_A$ acts on $\psi^{\parallel}$ and $\psi^{\perp}$ separately on $N_{\epsilon} \backslash K$.

\begin{lemma}\label{lemma of theta bound}
Suppose $\psi$ is any smooth configuration $N_{\epsilon} \backslash K$. Fix  $s$ and $\rho$, then $$\min\{4\gamma^2, \dfrac{(1-2\gamma)^2}{2}\} \int_0^{2\pi} (|\partial_{\theta}\psi|^2 + |\psi^{\perp}|^2) d\theta  \leq \int_0^{2\pi} |\nabla_{A_{\theta}} \psi|^2 d\theta \leq (1 + 4\gamma^2)\int_0^{2\pi} (|\partial_{\theta}\psi|^2 + |\psi^{\perp}|^2) d\theta,$$
$$\min\{4\gamma^2, (1-2\gamma)^2\} \int_0^{2\pi}  |\psi^{\perp}|^2 d\theta  \leq \int_0^{2\pi} |\nabla_{A_{\theta}} \psi|^2 d\theta. $$
\end{lemma}

\begin{proof}
Recall that as we discussed in subsection \ref{Subsection: Compare L with the operator N}, each Fourier component $e^{-im\theta - iks}h_j$ is also an eigenvector of the  operator $ \nabla_{\theta} = \partial_{\theta} + \gamma[\sigma, \cdot]$ with eigenvalue $-i(m + 2\gamma j)$.\\

Suppose $\psi_0 = w(\rho) e^{-im\theta - iks}h_j$ is one Fourier component of $\psi$. Then
$$ \int_0^{2\pi} |\nabla_{A_{\theta}} \psi_0|^2 d\theta = \int_0^{2\pi} |\partial_{\theta} \psi_0 + \gamma [\sigma, \psi_0]|^2 d\theta = (m + 2\gamma j)^2  \int_0^{2\pi} |\psi_0|^2 d\theta,$$
$$ \int_0^{2\pi} |\partial_{\theta}\psi_0|^2  d\theta = m^2  \int_0^{2\pi}|\psi_0|^2 d\theta, ~~~~  \int_0^{2\pi}  |\psi_0^{\perp}|^2 d\theta = j^2 \int_0^{2\pi}|\psi_0|^2 d\theta.$$
Recall that $m$ is an integer and $j = 0$ or $\pm 1$, it is easy to see 
$$\min\{4\gamma^2, \dfrac{(1-2\gamma)^2}{2}\}(m^2 + j^2) \leq (m + 2\gamma j)^2 \leq (1 + 4\gamma^2)(m^2 + j^2),$$
$$\min\{4\gamma^2, (1-2\gamma)^2\}j^2 \leq (m + 2\gamma j)^2.$$
Hence the lemma is true.
\end{proof}

\begin{lemma}\label{technical lemma 1}

For any $c > 0$, there exists an $\epsilon_1 \leq \delta$, for any $\epsilon \leq \epsilon_1$, if $\psi \in \IH$ is smooth, then
$$ \min{\{ 2\gamma - c, (1-2\gamma) - c\}} \int_{\partial N_{\epsilon}} |\psi^{\perp}|^2 \rho d\theta ds \leq  \epsilon \int_{\IR^3 \backslash N_{\epsilon}} |\nabla_A \psi^{\perp}|^2 d^3 x. $$
\end{lemma}

\begin{proof}
It may be assumed that $\psi = \psi^{\perp}$. We only prove the case $\gamma \leq \dfrac{1}{4}$, in which $ \min{\{ 2\gamma - c, (1-2\gamma) - c\}} = 2\gamma -  c$ for convenience. The other case ($\dfrac{1}{4} < \gamma < \dfrac{1}{2}$) follows by the same reason.\\

The proof has two steps:\\

\textbf{Step 1} proves the following statement: there exists an $\epsilon_1 \leq \delta$ such that, if $\psi \in \IH$ is smooth, then
$$  (2\gamma -  c) \int_{\partial N_{\epsilon_1}} |\psi^{\perp}|^2 \rho d\theta ds \leq  \epsilon_1 \int_{\IR^3 \backslash N_{\epsilon_1}} |\nabla_A \psi^{\perp}|^2 d^3 x. $$

Choose any $0 < 2\epsilon_1 \leq \delta$, then from theorem \ref{Appendix1} and theorem \ref{Appendix2}, there exists some constant $C(2\epsilon_1)$, such that 
$$2 \epsilon_1 \int_{\IR^3 \backslash N_{2\epsilon_1}} |\nabla_A \psi|^2  d^3 x  \geq C(2\epsilon_1) \int_{\partial N_{2\epsilon_1}} |\psi|^2 \rho d\theta ds. $$
If $C(2\epsilon_1) \geq 2\gamma - \gamma^3$, step 1 is done. Otherwise, $C(2\epsilon_1) < 2\gamma - \gamma^3$.\\

By lemma \ref{lemma of theta bound},
$$ 4\gamma^2\int_0^{2\pi} |\psi^{\perp}|^2 d\theta  \leq \int_0^{2\pi} |\partial_{\theta} \psi + \gamma [\sigma, \psi]|^2 d\theta, $$
so
$$ \int_0^{2\pi} |\nabla_A \psi|^2 d\theta \geq \int_0^{2\pi} (|\partial_{\rho}\psi|^2 +  \dfrac{4\gamma^2}{\rho^2}|\psi^{\perp}|^2) d\theta. $$

Choose $\mu > 0$ small enough such that $$4\gamma^2 - 2^{\mu} C(2 \epsilon_1)\mu - 2^{2\mu} C(2 \epsilon_1)^2 \geq 0,$$
\begin{equation*}
    \begin{array}{ll}
         ~ & ~~~~ \displaystyle \int_{\IR^3 \backslash N_{\epsilon_1}} |\nabla_A \psi|^2 d^3 x \geq \displaystyle \int_{ N_{2\epsilon_1} \backslash N_{\epsilon_1}} |\nabla_A \psi|^2 d^3 x  + \dfrac{C(2\epsilon_1)}{2\epsilon_1} \int_{\partial N_{2\epsilon_1}} |\psi|^2 \rho d\theta ds  \\
         & \geq \displaystyle \int_0^l\int_0^{2\pi}\int_{\epsilon_1}^{2\epsilon_1} (\rho |\partial_{\rho} \psi|^2  + \dfrac{4\gamma^2}{\rho} |\psi|^2) d\rho d\theta ds + C(2\epsilon_1) \int_0^l \int_0^{2\pi} |\psi(2\epsilon_1)|^2 d\theta ds\\
         &= \displaystyle \int_0^l\int_0^{2\pi}\int_{\epsilon_1}^{2\epsilon_1} (\rho |\partial_{\rho} \psi|^2  + \dfrac{4\gamma^2}{\rho} |\psi|^2) d\rho d\theta ds + (2 \epsilon_1)^{\mu} C(2\epsilon_1) \int_0^l \int_0^{2\pi} \int_{\epsilon_1}^{2\epsilon_1} \partial_{\rho}( \rho^{-\mu}  |\psi|^2)d\rho d\theta ds\\
         & ~~~ + \displaystyle 2^{\mu}C(2\epsilon_1) \int_0^l \int_0^{2\pi} |\psi(\epsilon_1)|^2 d\theta ds\\
         & \geq \displaystyle \int_0^l\int_0^{2\pi}\int_{\epsilon_1}^{2\epsilon_1} (\rho |\partial_{\rho} \psi|^2  + \dfrac{4\gamma^2}{\rho} |\psi|^2) d\rho d\theta ds -   (2\epsilon_1)^{\mu} C(2\epsilon_1) \mu\int_0^l \int_0^{2\pi} \int_{\epsilon_1}^{2\epsilon_1}  \rho^{-(1 + \mu)} |\psi|^2 d\rho d\theta ds\\
         & ~~~- \displaystyle  \int_0^l \int_0^{2\pi} \int_{\epsilon_1}^{2\epsilon_1}  ( \rho|\partial_{\rho} \psi|^2 + (2\epsilon_1)^{2\mu} C(2\epsilon_1)^2\rho^{-(2\mu+1)}|\psi|^2 ) d\rho d\theta ds  + \displaystyle \dfrac{2^{\mu}C(2\epsilon_1)}{\epsilon_1} \int_{\partial N_{\epsilon_1}} |\psi(\epsilon_1)|^2 d^2 x\\
         & = \displaystyle \int_0^l \int_0^{2\pi}\int_{\epsilon_1}^{2\epsilon_1} (\dfrac{4\gamma^2}{\rho} -  \dfrac{(2\epsilon_1)^{\mu} C(2 \epsilon_1)\mu}{\rho^{1 + \mu}} -  \dfrac{(2\epsilon_1)^{2\mu} C(2\epsilon_1)^2}{\rho^{2\mu +1}} ) |\psi|^2) d\rho d\theta ds\\
         & ~~~    + \displaystyle \dfrac{2^{\mu}C(2\epsilon_1)}{\epsilon_1} \int_{\partial N_{\epsilon_1}} |\psi(\epsilon_1)|^2 d^2 x\\
         & \geq \displaystyle \int_0^l \int_0^{2\pi}\int_{\epsilon_1}^{2\epsilon_1}  (4\gamma^2 - 2^{\mu} C(2 \epsilon_1)\mu - 2^{2\mu} C(2 \epsilon_1)^2)\dfrac{1}{\rho} |\psi|^2) d\rho d\theta ds  + \displaystyle \dfrac{2^{\mu}C(2\epsilon_1)}{\epsilon_1} \int_{\partial N_{\epsilon_1}} |\psi(\epsilon_1)|^2 d^2 x\\
         & \geq \displaystyle  \displaystyle \dfrac{2^{\mu}C(2\epsilon_1)}{\epsilon_1} \int_{\partial N_{\epsilon_1}} |\psi(\epsilon_1)|^2 d^2 x.
    \end{array}
\end{equation*}
 Let $C(\epsilon_1) = 2^{\mu}C(2\epsilon_1)$. Again if $C(\epsilon_1) \geq 2\gamma - c$, then step 1 is done. Otherwise, substitute $2\epsilon_1$ with $\epsilon_1$, repeat the above process until the constant be greater or equal than $2\gamma - c$, which finishes the first step.\\

\textbf{Step 2} proves that, for any $\epsilon \leq \epsilon_1$, if $\psi \in \IH_{ \epsilon}$ is smooth, then
$$    ( 2\gamma - c) \int_{\partial N_{\epsilon}} |\psi^{\perp}|^2 \rho d\theta ds \leq  \epsilon \int_{\IR^3 \backslash N_{\epsilon}} |\nabla_A \psi^{\perp}|^2 d^3 x. $$

In fact,
\begin{equation*}
    \begin{array}{ll}
         ~ & ~~~~ \displaystyle \int_{\IR^3 \backslash N_{\epsilon}} |\nabla_A \psi|^2 d^3 x \geq \displaystyle \int_{ N_{\epsilon} \backslash N_{\epsilon_1}} |\nabla_A \psi|^2 d^3 x  + \dfrac{1}{\epsilon_1} \int_{\partial N_{2\epsilon_1}} |\psi|^2 \rho d\theta ds  \\
         & \geq \displaystyle \int_0^l\int_0^{2\pi}\int_{\epsilon}^{\epsilon_1} (\rho |\partial_{\rho} \psi|^2  + \dfrac{4\gamma^2}{\rho} |\psi|^2) d\rho d\theta ds +  (2\gamma -  c)\int_0^l \int_0^{2\pi} |\psi(\epsilon_1)|^2 d\theta ds\\
         &= \displaystyle \int_0^l\int_0^{2\pi}\int_{\epsilon}^{\epsilon_1} (\rho |\partial_{\rho} \psi|^2  + \dfrac{4\gamma^2}{\rho} |\psi|^2) d\rho d\theta ds +  (2\gamma -  c) \int_0^l \int_0^{2\pi} \int_{\epsilon}^{\epsilon_1} \partial_{\rho}(  |\psi|^2)d\rho d\theta ds\\
         & ~~~ + \displaystyle  (2\gamma - c) \int_0^l \int_0^{2\pi} |\psi(\epsilon)|^2 d\theta ds\\
         & \geq  \displaystyle  (2\gamma -  c) \int_0^l \int_0^{2\pi} |\psi(\epsilon)|^2 d\theta ds,
    \end{array}
\end{equation*}
so the proposition is proved.
\end{proof}

\begin{lemma}\label{technical lemma 2}

Fix any large constant $C$. There exists an $\epsilon_2 \leq \epsilon$. Suppose $|M|$ is large enough (which depends on the knot, $\epsilon_2$ and $C$),  then if $\psi \in \IH$ is smooth. Then
$$  C \int_{\partial N_{\epsilon_2}} |\psi^{\perp}|^2 \rho d\Omega \leq  \epsilon \int_{\IR^3 \backslash N_{\epsilon_2}} (|\nabla_A  \psi^{\perp}|^2 + |[\Phi, \psi]|^2)d^3 x. $$
\end{lemma}

\begin{proof}
The proof is similar with the proof of lemma \ref{technical lemma 1}. From the start, there exists a constant $C(2\epsilon_1)$, such that
$$2 \epsilon_1 \int_{\IR^3 \backslash N_{2\epsilon_1}} (|\nabla_A \psi|^2 + 4M^2 |\psi^{\perp}|^2) d^3 x  \geq C(2\epsilon_1) \int_{\partial N_{2\epsilon_1}} |\psi^{\perp}|^2 \rho d\Omega. $$
If $C(2\epsilon_1) \geq C$, then the proof is done. Otherwise, $C(2\epsilon_1) < C$.\\

Since $|M|$ is large enough, it can be assumed that
$$4(M\epsilon_1)^2 - 2 C(2 \epsilon_1) - 4 C(2 \epsilon_1)^2 \geq 0.$$
\begin{equation*}
    \begin{array}{ll}
         ~ &  ~~~~ \displaystyle \int_{\IR^3 \backslash N_{\epsilon_1}} (|\nabla_A \psi|^2 + 4M^2 |\psi^{\perp}|^2) d^3 x\\
         & \geq \displaystyle \int_{ N_{2\epsilon_1} \backslash N_{\epsilon_1}} (|\nabla_A \psi|^2 + 4M^2 |\psi^{\perp}|^2) d^3 x  + \dfrac{C(2\epsilon_1)}{2\epsilon_1} \int_{\partial N_{2\epsilon_1}} |\psi^{\perp}|^2 \rho d\Omega  \\
         & \geq \displaystyle \int_0^l\int_0^{2\pi}\int_{\epsilon_1}^{2\epsilon_1} (\rho |\partial_{\rho} \psi^{\perp}|^2  + 4M^2 |\psi^{\perp}|^2) d\rho d\theta ds + C(2\epsilon_1) \int_0^l \int_0^{2\pi} |\psi^{\perp}(2\epsilon_1)|^2 d\theta ds\\
         &= \displaystyle \int_0^l\int_0^{2\pi}\int_{\epsilon_1}^{2\epsilon_1} (\rho |\partial_{\rho} \psi^{\perp}|^2  + 4M^2 |\psi^{\perp}|^2) d\rho d\theta ds + (2 \epsilon_1) C(2\epsilon_1) \int_0^l \int_0^{2\pi} \int_{\epsilon_1}^{2\epsilon_1} \partial_{\rho}( \rho^{-1}  |\psi^{\perp}|^2)d\rho d\theta ds\\
         & ~~~ + \displaystyle 2C(2\epsilon_1) \int_0^l \int_0^{2\pi} |\psi^{\perp}(\epsilon_1)|^2 d\theta ds\\
         & \geq \displaystyle \int_0^l\int_0^{2\pi}\int_{\epsilon_1}^{2\epsilon_1} (\rho |\partial_{\rho} \psi^{\perp}|^2  + 4M^2 |\psi^{\perp}|^2) d\rho d\theta ds -   (2\epsilon_1) C(2\epsilon_1)\int_0^l \int_0^{2\pi} \int_{\epsilon_1}^{2\epsilon_1}  \rho^{-2} |\psi^{\perp}|^2 d\rho d\theta ds\\
         & ~~~- \displaystyle  \int_0^l \int_0^{2\pi} \int_{\epsilon_1}^{2\epsilon_1}  ( \rho|\partial_{\rho} \psi^{\perp}|^2 + (2\epsilon_1)^{2} C(2\epsilon_1)^2\rho^{-3}|\psi^{\perp}|^2 ) d\rho d\theta ds  + \displaystyle \dfrac{2C(2\epsilon_1)}{\epsilon_1} \int_{\partial N_{\epsilon_1}} |\psi^{\perp}(\epsilon_1)|^2 d^2 x\\
         & = \displaystyle \int_0^l \int_0^{2\pi}\int_{\epsilon_1}^{2\epsilon_1}  (4M^2\rho -  \dfrac{(2\epsilon_1) C(2 \epsilon_1)}{\rho^2} -  \dfrac{(2\epsilon_1)^2 C(2 \epsilon_1)^2}{\rho^3} ) |\psi^{\perp}|^2 d\rho d\theta ds\\
         & ~~~~ + \displaystyle \dfrac{2C(2\epsilon_1)}{\epsilon_1} \int_{\partial N_{\epsilon_1}} |\psi^{\perp}(\epsilon_1)|^2 d^2 x\\
         & \geq \displaystyle \int_0^l \int_0^{2\pi}\int_{\epsilon_1}^{2\epsilon_1} (4(M\epsilon_1)^2 - 2 C(2 \epsilon_1) - 4 C(2 \epsilon_1)^2) \dfrac{1}{\rho} |\psi^{\perp}|^2 d\rho d\theta ds + \displaystyle \dfrac{2C(2\epsilon_1)}{\epsilon_1} \int_{\partial N_{\epsilon_1}} |\psi^{\perp}(\epsilon_1)|^2 d^2 x\\
         & \geq \displaystyle  \displaystyle \dfrac{2C(2\epsilon_1)}{\epsilon_1} \int_{\partial N_{\epsilon_1}} |\psi^{\perp}(\epsilon_1)|^2 d^2 x.
    \end{array}
\end{equation*}

Let $C(\epsilon_1) = 2C(2\epsilon_1)$. Again if $C(\epsilon_1) \geq C$, then the proof is finished. Otherwise, substitute $2\epsilon_1$ with $\epsilon_1$, repeat the above process until the constant be greater or equal than $C$, which is valid since $M$ is assumed to be large enough. This finished the proof.
\end{proof}

\subsection{The existence of the solutions}

Suppose $\Psi$ is the approximate solution defined in subsection \ref{Subsection: Approximate solutions from the standard 1 monopole}. In this subsection, we need to view $\tilde{L}_{\Psi}$ as an operator from $\IH$ to $\IL \oplus \IR$, where the additional component $\IR$ is already introduced in subsection \ref{Subsection: The moduli space as an analytical set}.

\begin{theorem}\label{local neighborhood manifold structure theorem}
There exists a small open neighbourhood of $0$ in $ \IH$, namely $U$,  such that, each point in $(\ker \tilde{L}_{\Psi}) \cap U$ corresponds to an actual solution to the Bogomolny equations which differs $\Psi$ by an element in $\IH$. In particular, such a solution exists.
\end{theorem}

\begin{proof}
Since $\tilde{L}_{\Psi}$ is Fredholm with 0 cokernel, there exists an $C > 0$, such that for any $W \in \IL \oplus \IR$, there exists a unique $\psi \in (\ker \tilde{L}_{\Psi})^{\perp} \subset \IH$, such that $\tilde{L}_{\Psi}(\psi) = W$. Moreover, $\norm{\psi}_{\IH} \leq C \norm{W}_{\IL \oplus \IR}$.\\

Choose any $\psi_0 \in (\ker \tilde{L}_{\Psi}) \cap U$. Consider $\psi_1, \psi_2, \cdots, \in (\ker \tilde{L}_{\Psi})^{\perp}$ such that  $$\tilde{L}_{\Psi}(\psi_1) = - V(\Psi + \psi_0) = -V(\Psi) - Q(\psi_0, \psi_0);$$
$$\tilde{L}_{\Psi}(\psi_2) = - V(\Phi + \psi_0 + \psi_1) = - 2Q(\psi_0, \psi_1) - Q(\psi_1, \psi_1);$$
$$\vdots$$
$$ \tilde{L}_{\Psi}(\psi_n) = - V(\Phi + \psi_0 + \cdots \psi_{n-1})  =  - 2 Q(\psi_0  + \cdots + \psi_{n-2}, \psi_{n-1}) - Q(\psi_{n-1}, \psi_{n-1}).$$
$$\vdots $$

Since $U$ is a small neighborhood and $\mu$ is also small enough, it may be assumed that $C^2 C_{\Psi}(\mu + C_{\Psi} \norm{\psi_0}^2_{\IH} + 2\norm{\psi_0}_{\IH}) \leq \dfrac{1}{2^4}.$ So
$$\norm{\psi_1}_{\IH} \leq C (\norm{V(\Psi)}_{\IL \oplus \IR} + \norm{Q(\psi_0, \psi_0)}_{\IL \oplus \IR}) \leq C(\mu + C_{\Psi}\norm{\psi_0}_{\IH}^2) \leq \dfrac{1}{2^4CC_{\Psi}} ; $$
$$\norm{\psi_2}_{\IH} \leq CC_{\Psi}\norm{\psi_1}_{\IH}(\norm{\psi_1}_{\IH} + 2\norm{\psi_0}_{\IH}) \leq \dfrac{1}{2}\norm{\psi_1}_{\IH}\leq \dfrac{1}{2^{5} C C_{\Psi} };$$
$$\norm{\psi_3}_{\IH} \leq CC_{\Psi}\norm{\psi_2}_{\IH}(\norm{\psi_2}_{\IH} + 2\norm{\psi_1}_{\IH} + 2\norm{\psi_0}_{\IH}) \leq \dfrac{1}{2}\norm{\psi_2}_{\IH}\leq \dfrac{1}{2^{6} C C_{\Psi} };$$
$$\vdots$$
$$\norm{\psi_n}_{\IH} \leq CC_{\Psi}\norm{\psi_{n-1}}_{\IH}(\norm{\psi_{n-1}}_{\IH} + 2\norm{\psi_{n-2}}_{\IH}  + \cdots + 2\norm{\psi_{0}}_{\IH}) \leq \dfrac{1}{2}\norm{\psi_{n-1}}_{\IH} \leq \dfrac{1}{2^{n+3} C C_{\Psi}}.$$
Hence $\psi = \sum\limits_{n=0}^{+\infty}\psi_n$ is in $\IH$ which is clearly a solution to the Bogonolmy equations.
\end{proof}

Here is an immediate corollary.

\begin{corollary}
For each $\gamma \neq \dfrac{1}{4}$ and non-negative integer $k$, there exists a large enough $M$ and a solution to the Bogomolny equations $\Psi$ such that
$$\Psi - \Psi_{\gamma, M, k} \in \IH. $$
\end{corollary}

\section{Behaviors of a solution when approaching infinity}\label{Section:Behaviors of a solution when approaching infinity}

In this section, we assume $\Psi = A + \Phi$ is a smooth solution to the Bogomolny equations on $\IR^3 \backslash N_{\epsilon}$. The only additional assumption (besides being smooth)  is:

$$\int_{\IR^3 \backslash N_{\epsilon}} |F_A|^2 d^3x < +\infty.$$

The goal of this section is to prove that there is a non-negative real number $M$ and (if $M > 0$) a non-negative integer $k$ such that $\Psi$ is a solution with mass $M$ and charge $k$ in the sense that we defined in section \ref{Section: Introduction and preliminary set-ups}. Moreover $\Psi$ satisfies all the good asymptotic behaviors approaching infinity.\\

The case $M = 0$ is possible but not studied in this paper. The author hopes to fill this gap in a future paper. Note that if there is not a knot singularity, then $M = 0$ would imply $\Psi = 0$ after a gauge transformation. So the solution is just a trivial one. It is likely that even if there is a knot singularity, $M = 0$ still implies at least $\Phi = 0$. But the author doesn't have a completed proof when writing this paper.\\

All the arguments in this section follow the spirits of Taubes' arguments in \cite{Jaffe1980VorticesTheories} (but details are slightly different). We need to be more careful in details because of the existence of a knot singularity.\\

\subsection{Assymptotic behaviors of $|\Phi|$ and $|\nabla_{A}^k\Phi|$ at infinity}\label{Subsection: Assymptotic behaviors of Phi and nabla Phi at infinity}

Throughout this section, assume $\chi$ is a cut-off function that equals $1$ on $\IR^3 \backslash N_{2\epsilon}$ and supported in $\IR^3 \backslash N_{\epsilon}$. Assume $\chi_R$ is another cut-off function which is $1$ in $B_R$ (the ball of radius $R$ centered at the origin) and equals $0$ in $\IR^3 \backslash B_{2R}$. We may always assume that  $R$ is large enough and $|\nabla^k (\chi_R)|$ is bounded above independently in $R$ (but the bound can depend on $k$).\\

\begin{lemma}\label{first lemma of Phi}
There are two non-negative functions $v_1$ and $v_2$ on $\IR^3$ such that

$$\chi |\Phi|^2 = v_1 + v_2,$$
and $v_1$ is bounded, $v_2$ is in $L^2(\IR^3)$.

\end{lemma}

\begin{proof}

 Suppose $\nabla_{A_i}$ be the covariant derivative $\nabla_A$ in the $x_i$ direction, where $x_1, x_2, x_3$ are Euclidean coordinates of $\IR^3$. Suppose $$*F_A = F_1dx_1 + F_2dx_2 + F_3dx_3.$$ 
 
 We have the Bogomolny equation
$$ \nabla_{A_i} \Phi = F_i ~~~~ \text{for} ~~~ i = 1, 2, 3.$$

And the Bianchi identity
 
 $$\sum\limits_{i=1}^3 \nabla_{A_i}^2 \Phi = \sum\limits_{i=1} \nabla_{A_i}F_i = 0.$$
 
So

$$\Delta (|\Phi|^2) = \sum\limits_{i = 1}^3 (2|\nabla_{A_i} \Phi|^2 + 2<\nabla_{A_i}^2 \Phi, \Phi>) = 2|\nabla_A \Phi|^2 = 2|F_A|^2.$$

$$\Delta (\chi|\Phi|^2) = 2 \chi |F_A|^2 + (\text{compactly supported terms}).$$

Thus $f = \Delta(\chi |\Phi|^2)$ is in $L^1(\IR^3)$. Let $G$ be the Green's function. Suppose $g = G(f)$. Then by Hardy-Littlewood-Sobolev inequality, we have

$$m \{|g| > \lambda\} \leq C \lambda^{-3}, $$
where $m\{|g| > \lambda\}$ is the measure of the set $\{|g| > \lambda\}$, $\lambda$ is any positive real number.\\

Let $h = g - \chi|\Phi|^2$. Then 
$m \{h > \lambda\} \leq C \lambda^{-3}.$ Which implies that
$$\dfrac{3}{4R^3}\int_{B_R(x_0)} h d^3x \leq \dfrac{C}{R},$$
where $B_R(x_0)$ is the ball of radius $R$ centered at $x_0$, and $C$ is a constant that doesn't depend on $R$ or $x_0$.\\

On the other hand, $h$ is a harmonic function. So by the mean value property 
$$h(x_0) = \dfrac{3}{4R^3} \int_{B_R(x_0)} u(x) d^3x.$$

Let $R \rightarrow +\infty$, we know that $h(x_0) \leq 0$. But a non-negative harmonic function on $\IR^3$ can only be a non-negative constant. So we know that there is a constant $M > 0$ such that $h = - M^2$.\\

Now we have $g = \chi|\Phi|^2 - M^2.$. Let $$v_1 = \min\{M^2 + 1, \chi|\Phi|^2\},  ~~ v_2 = \chi|\Phi|^2 - v_1.$$ 

Clearly $v_1$ is bounded. We have $v_2 \geq 0$ and
$$m \{v_2 > \lambda \} \leq C(\lambda + 1)^3. $$

So the 
$$\int_{\IR^3} |v_2|^2 d^3 x = \int_0^{+\infty} \int_{v_2 > \sqrt{\lambda}} m\{v_2 > \sqrt{\lambda}\} d\lambda  \leq C \int_0^{+\infty} (\sqrt{\lambda} + 1)^{-3}d\lambda < +\infty. $$

Thus $v_1$ and $v_2$ are what we want.
\end{proof}

\textbf{Remark:} Note that $M \geq 0$ is well-defined in the above proof. As what will be shown later, this number is the ``mass" of $\Psi$. We will continue to use the notation $M$ to represent the same number as in the above proof.

\begin{lemma}\label{lemma for gradient of F}
We have
$$\int_{\IR^3 \backslash N_{2\epsilon}} |\nabla_A F_A|^2 d^3 x < +\infty. $$
\end{lemma}

\begin{proof}
We assume $\chi_1 = \chi \cdot \chi_R$. In this proof, $C$ is a constant that doesn't depend on $R$. (But it may have different values from line by line.)
Then
$$\int_{\IR^3} |\nabla_A (\chi_1 F_A)|^2 d^3 x = - \int_{\IR^3\backslash N_{\epsilon}} \chi_1 \sum\limits_{i=1}^3<F_A, \nabla^2_{A_i}(\chi_1 F_A)> d^3x $$
$$\leq C\int_{\IR^3\backslash N_{\epsilon}} |F_A|^2 d^3x + \dfrac{1}{2}\int_{\IR^3\backslash N_{\epsilon}} |\nabla_A(\chi_1 F_A)|^2 d^3x + \int_{\IR^3 \backslash N_{\epsilon}} \chi_1^2 <F_A, \sum\limits_{i=1}^3\nabla_{A_i}^2 F_A> d^3 x.  $$
So
$$\int_{\IR^3 \backslash N_{\epsilon}} |\nabla_A (\chi_1 F_A)|^2 d^3 x  \leq C + 2\int_{\IR^3 \backslash N_{\epsilon}} \chi_1^2 <F_A, \sum\limits_{i=1}^3\nabla_{A_i}^2 F_A> d^3 x. $$

Suppose the Bogomolny equation is written as $$*F_A = d_A\Phi = F_1dx_1 + F_2dx_2 + F_3dx_3.$$ And we have the Bianchi identity
$$\sum\limits_{i=1}^{3}\nabla_{A_i} F_i = 0. $$

Thus we have
$$\sum\limits_{i=1}^{3} \nabla^2_{A_i} F_1 =- [[F_1, \Phi], \Phi] + 4[F_2, F_3].$$
So $$|\sum\limits_{i=1}^{3} \nabla^2_{A_i} F_1| \leq |F_A||\Phi|^2 + 4|F_A|^2.$$
Similar identity holds for $F_2, F_3$ as well. So 
$$|<F_A, \sum\limits_{i=1}^3\nabla_{A_i}^2 F_A>| \leq 3|F_A|^2|\Phi|^2 + 12|F_A|^3. $$
And

$$\int_{\IR^3 \backslash N_{\epsilon}} |\nabla_A (\chi_1 F_A)|^2 d^3 x  \leq C + 6\int_{\IR^3 \backslash N_{\epsilon}}\chi_1^2 (|F_A|^2|\Phi|^2 + 4|F_A|^3)d^3 x.$$

Note that $\chi|\Phi|^2 = v_1 + v_2$, where $v_1, v_2$ are the functions defined in the the proof of lemma \ref{first lemma of Phi}. In fact, we modify it a little bit to assume $|\Phi|^2 = v_1 + v_2$ here, where $v_1$ is bounded and $v_2$ is in $L^2(\IR^3 \backslash N_{\epsilon})$. So

$$\int_{\IR^3 \backslash N_{\epsilon}} |\nabla_A (\chi_1 F_A)|^2 d^3 x  \leq C + 6\int_{\IR^3 \backslash N_{\epsilon}}\chi_1^2( |F_A|^2v_2 + 4 |F_A|^3)d^3 x$$
$$\leq C + 6(\int_{\IR^3 \backslash N_{\epsilon}} |\chi_1 F_A|^6 d^3x)^{\frac{1}{4}}(\int_{\IR^3 \backslash N_{\epsilon}} |\chi_1 F_A|^2 d^3x)^{\frac{1}{4}}(\int_{\IR^3 \backslash N_{\epsilon}} (v_2^2 + |F_A|^2) d^3x)^{\frac{1}{2}}$$ $$ \leq C + C (\int_{\IR^3 \backslash N_{\epsilon}} |\chi_1F_A|^6 d^3x)^{\frac{1}{4}}. $$

By Sobolev embedding $L^{1,2}(\IR^3) \subset L^6$, and assuming $C$ is large enough, we have

$$\int_{\IR^3 \backslash N_{\epsilon}} |\nabla_A (\chi_1 F_A)|^2 d^3 x  \leq C + C (\int_{\IR^3 \backslash N_{\epsilon}} |\nabla_A(\chi F_A)|^2 d^3x)^{\frac{3}{4}}. $$

So 

$$\int_{B_R \backslash N_{2 \epsilon}}|\nabla_A(F_A)|^2 d^3 x \leq \int_{\IR^3 \backslash N_{\epsilon}} |\nabla_A(\chi F_A)|^2 d^3x \leq C. $$

Since $C$ doesn't depend on $R$, letting $R \rightarrow +\infty$ we know that

$$ \int_{\IR^3 \backslash N_{2\epsilon}} |\nabla_{A} F_A|^2 d^3 x \leq C.$$

\end{proof}

\begin{corollary}
\begin{itemize}
    \item The followings are in $L^2 (\IR^3 \backslash N_{2\epsilon})$:
$$|\nabla |\nabla |\Phi|||\leq   |\nabla (|F_A|)| = |\nabla |\nabla_A \Phi|| \leq |\nabla_A \nabla_A \Phi| = |\nabla_A F_A|.$$
\item The followings are in $L^p(\IR^3 \backslash N_{2\epsilon})$, where $2 \leq p \leq 6$:
$$|\nabla |\Phi|| \leq |\nabla_A \Phi| = |F_A|. $$
\item The following is in $L^{p}(\IR^3 \backslash N_{2\epsilon})$, where $6 \leq p \leq \infty$: $M - |\Phi|.$
\item Moreover, $\lim\limits_{|x| \rightarrow +\infty} (M - |\Phi(x)|) = 0$ uniformly in $|x|$.
\end{itemize}

\end{corollary}

\begin{proof}
The first bullet is obvious. In the second bullet, $L^2$ is given and $L^6$ is a standard Sobolev embedding theorem. Other $L^p$ with $2 < p < 6$ can be obtained by interpolation.\\

Once we know that $\lim\limits_{|x| \rightarrow +\infty}(M - |\Phi(x)|) = 0$ uniformly in $|x|$, then the third bullet follows by the following argument:\\

$L^{\infty}$ is obvious. $L^6$ is obtained by the fact that $\nabla (M - |\Phi|) = - \nabla |\Phi|$ is in $L^2$ and Sobolev inequality. The remaining cases are obtained by interpolation.\\

So it suffices to prove the fourth bullet only.

Recall that
$$\Delta (|\Phi|)^2 = 2|\nabla_{A}\Phi|^2.  $$

We know that $|\nabla_A \Phi|^2$ is in $L^p(\IR^3 \backslash N_{2\epsilon})$, where $1 \leq p \leq 3$. On the other hand, assuming $\chi$ is the cut-off function defined in the proof of lemma \ref{lemma for gradient of F}, then

$$\Delta(\chi |\Phi|^2) - (\Delta |\Phi|^2) $$

is supported in $N_{2\epsilon} \backslash N_{\epsilon}$, so $\Delta(\chi |\Phi|^2)$ has bounded $L^p(\IR^3)$ norm for any $1 \leq p \leq 3$.\\

Assume $G$ is the Green's function of $\Delta$. Let $g = G(f)$, where $f = \Delta (\chi |\Phi|^2)$. Then based the properties of $G$, for any $R > 0$

$$|l(y)| \leq \int_{|x - y| \leq R} \dfrac{C}{|x - y|} |f(y)| d^3x + \dfrac{C}{R} \int_{|x - y| \geq R} |f(y)| d^3x$$
$$\leq CR^{\frac{1}{2}}(\int_{|x - y|\leq R} |f(x)|^2 d^3x)^{\frac{1}{2}} + \dfrac{C}{R}\int_{\IR^3} |f| d^3x < +\infty,$$
where $C$ changes from line by line but doesn't depend on $y$. This  indicates that $\lim\limits_{|y| \rightarrow +\infty}|g(y)| = 0 $ uniformly in $y$.\\

We have already proved $g = \chi|\Phi|^2 - M^2$ in the proof of lemma \ref{first lemma of Phi}. This finishes the proof.
\end{proof}

\begin{corollary}
Suppose $M > 0$. When $R$ is large enough, $\dfrac{\Phi}{|\Phi|}$ gives a well-defined unit section on $\partial B_R$. Thus the charge ``$k$" is well-defined. 
\end{corollary}

\begin{lemma}
In fact, we can do the same estimate and get that 

$$\int_{\IR^3 \backslash N_{2\epsilon}} |\nabla_A^m F_A|^2 d^3 x < +\infty$$
for any non-negative integer $m$.
\end{lemma}

\begin{proof}
We prove by induction. We know that when $m = 0, 1$ this is already proved. Assume it is proved for $1, 2, \cdots, m-1$. In particular, since $F_A$ is smooth, we actually have

$$\int_{\IR^3\backslash N_{\epsilon}} | \nabla_{A}^j F_A|^2 d^3x < +\infty,$$
for $j = 1, 2, \cdots, m-1$.\\

If $ 1 \leq j \leq m - 2$, then  $\chi |\nabla_A^j F_A|$ is in $L^6(\IR^3)$ by Sobolev embedding, and hence $|\nabla_{A}^j F_A|$ in $L^p(\IR^3 \backslash N_{\epsilon})$ for any $2 \leq p \leq 6$.\\

In addition, if $1 \leq j \leq m - 4$

$$\lim\limits_{|x| \rightarrow +\infty} |\nabla^j_A F_A| = 0 $$
uniformly in $|x|$ for the following reason:\\

Since $\chi |\nabla^j_A F_A|$ is in $L^{2,p}(
\IR^3)$ for any $2 \leq p \leq 6$, 
$$\Delta(\chi |\nabla_A^j F_A|^2)| \leq C + C|\nabla_A^{j+1}F_A|^2 + C|\nabla_{A}^jF_A||\nabla_A^{j+2}F_A| \in L^1(\IR^3 \backslash N_{\epsilon}) \cap L^2(\IR^3 \backslash N_{\epsilon}).$$

Then use a similar argument with the proof of lemma \ref{first lemma of Phi}, we get $$\lim\limits_{|x| \rightarrow +\infty} |\nabla^j_A F_A| = 0 $$
uniformly in $|x|$.\\

Let $\chi_1 = \chi \chi_R$. Suppose $C$ is a large enough constant that doesn't depend on $R$, but may change from line to line. Consider 

$$\sum\limits_{i=1}^3\int_{\IR^3 \backslash N_{\epsilon}} |\nabla_{A_i}\nabla_{A}^{m-1} (\chi_1 F_A)|^2 d^3 x \leq  \int_{\IR^3 \backslash N_{\epsilon}}|<\nabla_{A}^{m-1}(\chi_1 F_A), \sum\limits_{i=1}^3\nabla_{A_i}^2\nabla_A^{m-1}(\chi_1 F_A) > | d^3 x$$
$$\leq C \int_{\IR^3 \backslash N_{\epsilon}} ( ( \chi_1|\nabla_A^{m-1} F_A| ) ( \sum\limits_{j_1 + j_2 = m-2}\chi_1|\nabla_A^{j_1}F_A||\nabla_A^{j_2}F_A| ) + \chi_1^2\sum\limits_{i=1}^3 |\nabla_{A}^{m-1}( F_A)|  |\nabla_A^{m-1}\nabla_{A_i}^2 F_A|) d^3 x + C.$$

 Note that $$|\nabla_A^{m-1} \sum\limits_{i=1}^3(\nabla_{A_i}^2 F_A)| \leq C |\nabla_{A}^{m-1}(|F_A||\Phi|^2 + |F_A|^2)|.$$ So all the terms are bounded above by
 $$C\int_{\IR^3 \backslash N_{\epsilon}}\chi_1^2|\nabla_A^{m-1}F_A|^2(1 + |F_A|) d^3x + C. $$

By the same argument as the proof of lemma \ref{first lemma of Phi}, using $L^6 \subset L^{1,2}$, we conclude that

$$\sum\limits_{i=1}^3\int_{\IR^3 \backslash N_{\epsilon}} |\nabla_{A_i}\nabla_{A}^{m-1} (\chi_1 F_A)|^2 d^3 x \leq  < C.$$

Let $R \rightarrow +\infty$, then the proof is finished.

\end{proof}

\begin{corollary}
From the above proof we know that 
$$\lim\limits_{|x| \rightarrow +\infty} |\nabla_A^m F_A| = 0$$
uniformly in $|x|$ (but not in $m$). In particular, $|\nabla_A^m F_A|$ is in $L^p(\IR^3 \backslash N_{\epsilon})$ for all $2 \leq p \leq \infty$. 
\end{corollary}

\subsection{More on assymptotic behaviors}

Suppose $\Psi$ is a solution to the Bogomolny equations on $\IR^3 \backslash B_R$. We have already seen that the ``mass" $M$ is well-defined from $\Psi$. In this subsection, we assume $M > 0$.

\begin{lemma}
Let $\eta = [F_A, \Phi]$. Then for any $ 0 < \delta < M$, there is a constant $C$ which depends on $\delta$ such that $$|\eta| \leq C e^{- \delta |x|}.$$
\end{lemma}
\begin{proof}

$$\Delta |\eta|^2  = 2\sum\limits_{i=1}^3 (|\nabla_{A_i}\eta|^2 + <\eta, \nabla^2_{A_i}[F_A, \Phi]>) $$
$$ = 2\sum\limits_{i=1}^3 (|\nabla_{A_i} \eta|^2 + <\eta, [\nabla^2_{A_i}F_A, \Phi] + 2[\nabla_{A_i}F_A, F_i]>).  $$

We have $$\sum\limits_{i=1}^3[\nabla_{A_i}^2F_1, \Phi]  = - [[[F_1, \Phi], \Phi], \Phi] + 4[[F_2, F_3], \Phi] = |\Phi|^2[F_1, \Phi] + 4[[F_2, F_3], \Phi]. $$

Similar identity holds for $F_2, F_3$. So

$$ \sum\limits_{i, j=1}^3 <[F_j, \Phi], [\nabla_{A_i}^2 F_j, \Phi]> =  |\eta|^2|\Phi|^2 + \sum\limits_{i, j=1}^3 4<[F_1, \Phi], [[F_2, F_3], \Phi]>. $$

Note that $$|[F_1, \Phi]| \cdot |[[F_2, F_3], \Phi]| \leq |\eta| \cdot (|[[F_2, \Phi], F_3]| + |[[F_3, \Phi], F_2]|) \leq 4|\eta|^2 |F_A| = |\eta|^2 \cdot o(1),$$
where $o(1)$ is a term that approaches $0$ uniformly when $|x| \rightarrow +\infty$.\\

On the other hand, 

$$[\nabla_{A_i} F_j, F_i] = \dfrac{1}{|\Phi|^2} ([[\nabla_{A_i} F_j, \Phi], [F_i, \Phi]] ~ + <\nabla_{A_i}F_j, \Phi>[\Phi, F_i]~ + <F_i, \Phi> [\nabla_{A_i} F_j, \Phi]).$$
The first term above is perpendicular with $\eta$. The second term is bounded above by $|\eta| \cdot o(1)$. For the last term, we have

$$[\nabla_{A_i} F_j, \Phi] = \nabla_{A_i}[F_j, \Phi] - [F_j, F_i]$$ $$= \nabla_{A_i}[F_j, \Phi] - \dfrac{1}{|\Phi|^2}([[F_j, \Phi], [F_i, \Phi]] + <F_i, \Phi> [F_j, \Phi] + <F_j, \Phi> [\Phi, F_i]),$$
where $[[F_j, \Phi], [F_i, \Phi]]$ is perpendicular with $\eta$, and the remaining terms are bounded above by $$|\nabla_A \eta| ~ + |\eta| \cdot o(1). $$
Thus

$$<\eta, \sum\limits_{i=1}^{3}[\nabla_{A_i} F_A, F_i]> ~ \leq ~ |\eta|^2 \cdot o(1) + |\eta| |\nabla_A \eta| \cdot o(1) = (|\eta|^2 + |\nabla_A \eta|^2) \cdot o(1). $$

After all

$$\Delta |\eta|^2  = 2|\nabla_A \eta|^2 + |\eta|^2 \cdot |\Phi|^2 +  (|\eta|^2 + |\nabla_A \eta|^2) \cdot o(1)$$ $$= 2|\nabla_A \eta|^2 + M^2 |\eta|^2 + (|\eta|^2 + |\nabla_A \eta|^2) \cdot o(1). $$

Now consider $s(x) = C e^{-\delta |x|}$. We have

$$\Delta s = (\delta^2 - \dfrac{2\delta}{|x|}) C e^{- \delta |x|} \leq \delta^2 s(x) \leq (M^2 - o(1))s(x).$$

So $\Delta (|\eta|^2 - s) \geq M^2(|\eta|^2 -s)$ for large enough $|x|$ uniformly. In particular, we may assume that $R$ is large enough such that

$$\Delta(|\eta|^2 - s) \geq M^2 (|\eta|^2 - s) ~~ \text{on} ~~ \IR^3 \backslash B_R.$$

On the other hand, we may assume $C$ is large enough such that
$$|\eta|^2 - s \leq 0 ~~ \text{on} ~ \IR^3 \backslash B_R.$$

Since $$\lim\limits_{|x| \rightarrow +\infty} |\eta|^2 - s = 0,$$
by a standard maximal principle, we know that
$$|\eta|^2 - s \leq 0 ~~ \text{on} ~ \IR^3 \backslash B_R.$$

This finishes the proof.

\end{proof}

In fact, by a similar but more cumbersome argument, we can also show that, $$|\nabla_A \eta| = O(e^{-\delta |x|}).$$ 

Here it is

\begin{lemma}

\end{lemma}

\begin{proof}
Note that $$\nabla_A^k \eta  = [\nabla_{A}^k F_A, \Phi] + (\cdots),$$
where $(\cdots)$ are terms that look like $[\nabla_A^j F_A, \nabla_A^{k - 1 - j} F_A]$ for $j \leq k - 1$. We may use induction to assume that this term is bounded above by a constant times
$$|\dfrac{1}{|\Phi|^2}[[[\nabla_A^jF_A, \Phi], \Phi], F_{A}^{k - 1 - j}F_A]| + |<\nabla_A^kF_A, \Phi> [\Phi, \nabla_a^{K - 1 - J}F_A]|,$$
which further bounded above by 
$Ce^{-\delta |x|}.$ Thus we only need to show that 
$$[\nabla_A^k F_A, \Phi] = O(e^{-\delta |x|}).$$

In fact, let $\xi = [\nabla_A^k F_A, \Phi]$, then
$$\Delta(|\xi|^2) = |\nabla_A \xi|^2 + \sum\limits_{i = 1}^3<\xi, \nabla_{A_i}^2[\nabla_A^k F_A, \Phi]>. $$

Note that $<\xi, \sum\limits_{i = 1}^3 \nabla_{A_i}^2[\nabla_A^k F_A, \Phi]>$ can be written as 
$ |\xi|^2 |\Phi|^2 + (\cdots)$, where by induction $(\cdots)$ is a sum of terms that either bounded above by $Ce^{- \delta |x|}$ , or bounded above by $o(1) \cdot (|\xi|^2 + |\nabla_A\psi|^2). $ So after all,

$$\Delta |\xi|^2 \geq M^2 |\xi|^2 - C(e^{- \delta |x|}).  $$

By the same maximal principle argument as in the previous lemma, we get for a slightly smaller $\delta$,
$$|\xi| = O(e^{- \delta |x|}).  $$

\end{proof}

\subsection{$\Psi$ differs the model configuration by an element in $\IH$ after a gauge transformation near infinity}

The same assumption as in the last subsection: Suppose $\Psi$ is a solution to the Bogomolny equations on $\IR^3 \backslash B_R$. Suppose the ``mass" $M$ is positive.

\begin{theorem}
After a gauge transformation, for some $k$, $\Psi$ satisfies 
$$(1 - \chi_{R})(\Psi - \Psi_{M, k}) \in \IH.$$
\end{theorem}

\begin{proof}
Recall that $\Psi = A + \Phi$. We may always use a gauge transformation to make $\Phi$ parallel to $e_0$. Now we assume $\Phi = |\Phi| e_0$. We have already proved (in subsection \ref{Subsection: Assymptotic behaviors of Phi and nabla Phi at infinity}) that $|\Phi| = M + o(1)$ as $|x| \rightarrow +\infty$.\\

The Bogomolny equations say
$$F_A = *d_A \Phi = *d(|\Phi|) e_0 + |\Phi| (*d_A e_0).$$

However,
$$[F_A, \Phi] = [F_A, |\Phi| e_0] = |\Phi|^2 [*d_Ae_0, e_0]. $$

And we know from the previous subsection that

$$|[F_A, \Phi]| \leq C e^{- \delta |x|}.$$

On the other hand, $*d_A e_0$ is perpendicular with $e_0$. So $|[*d_A e_0, e_0]| = |*d_A e_0|.$ So we have

$$|d_A e_0| \leq C e^{- \delta |x|}.  $$

So in fact $e_0$ is almost a covariantly invariant section under $A$ when $|x|$ is large. In particular, let $A_0$ be the connection similar with $A$ but
$d_{A_0} e_0 = 0$. To be more precise,
$$A_0 = [de_0, e_0] + <A, e_0>e_0. $$Then
$$O(e^{-\delta |x|}) = (d_A - d_{A_0})e_0 = [A - [de_0, e_0], e_0].  $$
This implies that
$$|A - A_0| = O(e^{-\delta |x|}).$$

So we may write $A$ as

$$A = [de_0, e_0] + <A, e_0> e_0 + O(e^{- \delta |x|}).$$
And thus
$$F_A = - \dfrac{1}{2} [de_0, de_0] + (d<A, e_0>) e_0 + O(e^{-\delta |x|}). $$

On the other hand,
$$*F_A = d_A \Phi = d(|\Phi|) e_0 + O(e^{-\delta |x|}). $$

One can check directly that

$$- \dfrac{1}{2}[de_0, de_0]=  *\dfrac{k}{|x|^2}d(|x|) e_0.$$ And recall $r = |x|$ on $\IR^3 \backslash B_R$, 

$$*d|\Phi| =  *\dfrac{k}{r^2} dr+ d<A, e_0> + O(e^{-\delta |x|}). $$
So 
$$\Delta |\Phi| = d*d|\Phi| = O(e^{-\delta |x|}) ~~~ (\text{as}~~ |x| \rightarrow +\infty).$$

Suppose $$\Delta ((1 - \chi_R)|\Phi|) = f = O(e^{-\delta |x|}).$$ Let $G$ be the Green's function of the Laplacian. Then $|\Phi| - G(f)$ is a harmonic function on $\IR^3$. Note that one can check directly that
$$G(f) = \dfrac{c}{|x|} + O(\dfrac{1}{|x|^2}), ~~~ \text{as} ~~ |x| \rightarrow +\infty,$$
where $c$ is a constant. In fact,
$$G(f)(x) = \int_{\IR^3} \dfrac{1}{4\pi |x - y|} f(y) d^3y = \dfrac{1}{4\pi|x|} (\int_{\IR^3} f(y)d^3y) + \dfrac{1}{4\pi}\int_{\IR^3} (\dfrac{1}{|x - y|} - \dfrac{1}{|x|}) f(y)d^3y,$$
where $$|\dfrac{1}{4\pi}\int_{\IR^3} (\dfrac{1}{|x - y|} - \dfrac{1}{|x|}) f(y)d^3y| \leq \dfrac{1}{4\pi}\int_{\IR^3} \dfrac{|y||f(y)|}{|x - y| |x|} d^3y = O(\dfrac{1}{|x|^2}). $$
Together with the fact that $\lim\limits_{|x| \rightarrow +\infty} |\Phi| = M$, we get for large $|x|$,
$$|\Phi| = M + G(f) =  M + \dfrac{c}{|x|} + O(\dfrac{1}{|x|^2}).$$

In fact, the expression of $G(f)$ has an expansion for large $|x|$
$$G(f) = \dfrac{c}{|x|} + \sum\limits_{j=2}^{+\infty} \dfrac{G(f)_{j}(S)}{|x|^j},$$
where each $G(f)_j(S)$ is a smooth function on the sphere. In particular, recall that $r = |x|$ on $\IR^3 \backslash B_R$. Then for large $r$,
$$d (|\Phi|) = - \dfrac{c}{r^2} dr + O(\dfrac{1}{r^3}).$$

Note that this already implies that $(1 - \chi_R)(\Phi - \Phi_{M, k}) \in \IH$.\\

On the other hand, 

$$d|\Phi| = -\dfrac{k}{r^2} dr + *d<A, e_0> + O(e^{-\delta |x|}). $$

So
$$(\dfrac{k - c}{r^2}) * dr = d<A, e_0> + O(\dfrac{1}{r^3}).$$

Since $d<A, e_0>$ is an exact form on $\partial B_r$. Taking the integral along $\partial B_r$ as $r \rightarrow +\infty$ implies that $k = c$. So

$$d<A, e_0> = O(\dfrac{1}{r^3}). $$

It is standard that after a $U(1)$ gauge transformation that doesn't affect $\Phi$ (which changes $<A, e_0>$ by an closed $1$-form on $\IR^3 \backslash B_R$), we may assume 
$$d(*<A, e_0>) = 0.$$
This implies that $$\nabla (<A, e_0>) = O(\dfrac{1}{r^3}).$$ In particular, $$\partial_r (|<A, e_0>|) = O(\dfrac{1}{r^3}).$$ Integrating along the $r$ direction, we get
$$<A, e_0> = O(\dfrac{1}{r^2}).$$ This implies

$$|A_0 - [de_0, e_0]| = O(\dfrac{1}{r^2}), ~~~ |\nabla(A_0 - [de_0, e_0])| = O(\dfrac{1}{r^3}). $$

Together with the fact that 
$$|A - A_0| = O(e^{-\delta r}), ~~~ |\nabla(A - A_0)| = O(e^{-\delta r}), $$
it implies that $(1 - \chi_R)(A - A_{M, k}) \in \IH$. So finally,
$$(1 - \chi_R)(\Psi - \Psi_{M, k}) \in \IH.$$
This is what we want to prove.
\end{proof}

\section{A Ulenbeck type of regularity theory near the knot} \label{Section: A Ulenbeck type of regularity theory near the knot}
The goal of this section is to prove the following theorem:

\begin{theorem}\label{middle theorem}
Suppose $\Psi = A + \Phi$ is a smooth solution to the Bogomolny equations on $N_{2\epsilon} \backslash K$ with 

$$\int_{N_{2\epsilon} \backslash K} \rho^2 |F_A|^2 d\rho d\theta ds < +\infty. $$

Then either $\Psi$ can be extended to a smooth solution on $N_{2\epsilon}$, or there exists a constant $\gamma \in (0, \dfrac{1}{2})$ such that $\Psi$ has a knot singularity of monodromy $\gamma$.
\end{theorem}

The idea of the proof is essentially borrowed from L. M. Sibner and  R. J. Sibner's arguments in \cite{Sibner1992ClassificationHolonomy}.
\subsection{A remark on the domain of $\gamma$}

We have always assumed that $\gamma \in (0, \dfrac{1}{2})$. In principle, one can also define a knot singularity with monodromy $\gamma$ such that $\gamma \notin (0, \dfrac{1}{2})$. However, the following propositions indicate that it suffices to study the case $\gamma \in (0, \dfrac{1}{2})$.

\begin{proposition}
Suppose $\Psi$ is a smooth configuration on $N_{2\epsilon} \backslash K$ with such that
$$\chi (\Psi - \Psi_{\gamma, M, 0}) \in \IH $$ for some $\gamma \notin (0, \dfrac{1}{2})$. Then there is an $SO(3)$ gauge transformation that sends $\Psi$ to $u(\Psi)$ such that
$$\chi (u(\Psi) - \Psi_{\gamma', M, 0}) \in \IH $$
for some $\gamma' \in [0, \dfrac{1}{2})$.
\end{proposition}

\begin{proof}Let $u = e^{\frac{1}{2}i\theta}\sigma$ near the knot. Note that $u$ has a sign ambiguity. (The function $e^{\frac{1}{2}i\theta}$ is only well-defined up to a sign.) So it represents an $SO(3) = SU(2)/\{\pm I\}$ gauge transformation. But clearly $u$ sends $\Psi_{\gamma, M, 0}$ to $\Psi_{\gamma + \frac{1}{2}, M, 0}$. Thus $\Psi_{\gamma, M, 0}$ and $\Psi_{\gamma', M, 0}$ differ by a gauge transformation as long as $\gamma - \gamma'$ is an integer or a half integer. So the remaining task is to prove that the gauge transformation $u$ preserves $\IH$ near the knot.\\

 Suppose $\psi \in \IH$. We may assume $\psi$ is supported in $N_{\epsilon} \backslash K$. Then clearly $|u| = 1$ and suppose $u(\psi) = u \psi u^{-1}$, then $u(\psi)$ differes $\psi$ on $N_{\epsilon} \backslash K$ only and
$$|u(\psi)| = |\psi|, ~~~ |\nabla(u(\psi))| \leq 4 |d u| |\psi| + |d\psi| \leq \dfrac{C}{\rho} |\psi| + |d\psi|.$$

Clearly $$ \rho^2|\nabla(u(\psi))|^2 +  |u(\psi)|^2 \leq C(\rho^2 |\nabla \psi|^2 + |\psi|^2). $$
So $u(\psi) \in \IH$. 

\end{proof}

\begin{proposition}
If $\Psi$ is a smooth solution to the Bogomolny equations on $N_{2\epsilon} \backslash K$ such that
$$\chi(\Psi - \Psi_{0, M, 0})\in \IH. $$
Then $\Psi$ can be extended to a smooth configuration on $N_{2\epsilon}$.
\end{proposition}

The proof of this proposition defers to section \ref{Section: Limit behaviors near the knot}.

Given the above two propositions, theorem \ref{middle theorem} is equivalent with the following proposition:

\begin{proposition}\label{final regularity proposition to prove}
Suppose $\Psi = A + \Phi$ is a smooth solution to the Bogomolny equations on $N_{2\epsilon} \backslash K$ with 

$$\int_{N_{2\epsilon} \backslash K} \rho^2 |F_A|^2 d\rho d\theta ds < +\infty. $$

Then there exists a $\gamma \in [0, \dfrac{1}{2})$ such that,
after a gauge transformation, $\Psi$ satisfies
$$\chi(\Psi - \Psi_{\gamma, M, 0}) \in \IH.$$
\end{proposition}

The remaining of this section proves the above proposition.

\subsection{A seemingly weaker statement}

\begin{lemma}\label{lemma for the size of Phi near the knot}
Suppose $\Psi$ is a smooth solution to the Bogomolny equations on $N_{\epsilon} \backslash K$ such that
$$\int_{N_{\epsilon} \backslash K} \rho^2|F_A|^2 d\rho d\theta ds < +\infty.$$
Then 
$$\int_{N_{\epsilon} \backslash K} |\Phi|^2 d\rho d\theta ds < +\infty.$$
\end{lemma}

\begin{proof}
Note that $$|\nabla|\Phi|| \leq |\nabla_A \Phi| = |F_A|.$$
So the lemma \ref{lemma for the size of Phi near the knot} is a direct corollary of lemma \ref{lemma of a preliminary inequality on N_epsilon}.
\end{proof}

The upcoming lemma \ref{local regularity lemma for connection} seems to be weaker than proposition \ref{final regularity proposition to prove} in some sense. But we'll show that actually it implies proposition \ref{final regularity proposition to prove}.

\begin{lemma}\label{local regularity lemma for connection}

Suppose $A$ is an $su(2)$-valued 1-form on $N_{2\epsilon} \backslash K$, representing a connection on the trivial $SU(2)$-bundle. Suppose the curvature is $F_A$. Suppose $\chi F_A$ has finite $\IL$ norm. By shrinking $\epsilon$, it may be assumed that the $\IL$ norm of $\chi F_A$ is small enough. Then there exists a constant $C$ which doesn't depend on $A$, and a flat fiducial connection namely $A^{f} = \gamma \sigma d\theta + \tilde{\gamma} \sigma ds$ for some constants $\gamma$ and $\tilde{\gamma}$, and a smooth gauge transformation $u$, such that after the gauge transformation (so that $A$ becomes u(A)),
$$\int_{N_{\epsilon} \backslash K} \rho (|\nabla(u(A) - A^{f})|^2 + |[\sigma, (u(A) - A^f)]|^2) d^3x \leq C \norm{\chi F_A}^2_{\IL}.$$
\end{lemma}

\paragraph{Proof that lemma \ref{local regularity lemma for connection} implies theorem \ref{final regularity proposition to prove}}~\\

We assume lemma \ref{local regularity lemma for connection} is true here. Suppose $\Psi = A + \Phi$ is a smooth solution to the Bogomolny equation on $N_{2\epsilon} \backslash K$ such that
$$\int_{N_{2\epsilon} \backslash K} \rho^2 |F_A|^2 d\rho d\theta ds < +\infty.$$ Suppose $u$ is the gauge transformation given in lemma \ref{local regularity lemma for connection}. Let $a = \chi(u(A) - A^{f})$. Then from the assumption,

$$\int_{N_{\epsilon} \backslash K} \rho^2 |\nabla a|^2 d\rho d\theta ds < +\infty.$$

By lemma \ref{lemma of a preliminary inequality on N_epsilon}, we know that $a \in \IH.$ On the other hand, one can check directly that 

$$\chi(A^f - \Psi_{\gamma, M, 0}) \in \IH.$$

So the remaining task is to prove that

$$\chi u(\Phi) \in \IH.$$

By Bogomolny equations,

$$* u(F_A) = du(\Phi) + [u(A), u(\Phi)]. $$
Assume $C$ is a large constant that may change from line to line which doesn't depend on $\delta$. Assume $\delta \ll \epsilon$. Suppose $\chi_{\delta}$ is a smooth cut-off function that is $1$ on $N_{\delta}$ and $0$ on $\IR^3 \backslash N_{2 \delta}$ and $|\nabla \chi_{\delta}| \leq \dfrac{2}{\delta}.$
$$\int_{N_{2\epsilon} \backslash K} \rho^2 |\nabla (( 1 - \chi_{\delta}) u(\Phi))|^2 d\rho d\theta ds$$ $$\leq  \int_{N_{2\epsilon} \backslash K} ((C^2 |F_A|^2 +  C^2|a|^2|\Phi|^2 + (1 + \dfrac{1}{C}) |A^f|^2 |\Phi|^2)(1 - \chi_{\delta})^2\rho^2 + C|\nabla((1 - \chi_{\delta})u(\Phi))||\Phi| \rho + C|\Phi|^2) d\rho d\theta ds $$
$$\leq C^2 + \int_{N_{2\delta} \backslash K} ( C^2|a|^2 (1 - \chi_{\delta})^2 |\Phi|^2 + \dfrac{1}{C} |\nabla((1 - \chi_{\delta})u(\Phi))|^2) \rho^2 d\rho d\theta ds.$$

So

$$\int_{N_{2\epsilon} \backslash K} \rho^2 |\nabla ((1 - \chi_{\delta})u(\Phi))|^2 d\rho d\theta ds \leq C + C \int_{N_{2\delta} \backslash K} \rho^2 (1 - \chi_{\delta})^2 |a|^2 |\Phi|^2 d\rho d\theta ds$$
$$\leq C + C ((\int_{N_{2\delta} \backslash K}\rho^2 |\nabla a|^2)^{\frac{1}{2}} + C) ((\int_{N_{2\epsilon} \backslash K} \rho^2 |\nabla ((1 - \chi_{\delta})u(\Phi))|^2 d\rho d\theta ds )^{\frac{1}{2}} + C).$$

The last step is by lemma \ref{lemma of a preliminary inequality on N_epsilon} and lemma \ref{weight version of Sobolev embedding lemma}. And the above inequality is enough to get

$$\int_{N_{2\epsilon} \backslash K} \rho^2 |\nabla ((1 - \chi_{\delta})u(\Phi))|^2 d\rho d\theta ds \leq C.$$

Since $C$ doesn't depend on $\delta$, we have
$$\int_{N_{2\epsilon} \backslash K} \rho^2 |\nabla (u(\Phi))|^2 d\rho d\theta ds < +\infty.$$

And by lemma \ref{lemma of a preliminary inequality on N_epsilon} again and the fact that $u(\Phi)$ is smooth, we know that $$\chi u(\Phi) \in \IH.$$

The remaining of this section is to prove lemma \ref{local regularity lemma for connection}.

\subsection{Some preliminary Uhlenbeck type analysis}

\begin{lemma}(Uhlenbeck lemma) \label{Uhlenbeck lemma}
There exists a large enough constant $C$ such that, supposing $B_r$ is a ball of radius $r$, and supposing
$$\int_{B_r} |F_A|^2 d^3x \leq \dfrac{1}{Cr}, $$
then after some gauge transformation $u$, letting $u(A) = d^*(uAu^{-1} - (du)u^{-1})$,
\begin{equation*}
    \left\{\begin{array}{ll}
        d^*(u(A)) = 0, ~~~ *u(A)|_{\partial B_r} = 0;&  \\
        \displaystyle \int_{B_r}(|\nabla A|^2 + \dfrac{1}{r^2}|A|^2) d^3x \leq  C  \int_{B_r} |F_A|^2 d^3x.& 
    \end{array}
    \right.
\end{equation*}
Moreover, the above gauge transformation $u$ is unique up to a constant gauge transformation.
\end{lemma}

\begin{proof}
$$\int_{B_r}|F_A|^{\frac{3}{2}}d^3x \leq (\int_{B_r}|F_A|^2 d^3x )^{\frac{3}{4}} (\int_{B_r} 1 d^3 x)^{\frac{1}{4}} \leq (\dfrac{4\pi}{3C^3})^{\frac{1}{4}}. $$
Since $C>0$ is large enough, if may be assumed that $(\dfrac{4\pi}{3C^3})^{\frac{1}{4}}$ is smaller than the Ulenbeck's constant. Hence the lemma can be proved by a standard Ulenbeck's argument.
\end{proof}

\begin{lemma}\label{Uhlenbeck transition estimate lemma}
Suppose $B_1, B_2$ are two open balls in $\IR^3$ that overlap. Let $V$ be an open subset such that $V \subset \bar{V} \subset B_1 \cap B_2$. Suppose $A$ is an $su(2)$-connection on $B_1 \cup B_2$ such that $\norm{F_A}_{\IL^2(B_1 \cup B_2)}$ is small enough, where $F_A$ is the curvature. Suppose $a, b$ are the $su(2)$-valued 1-forms that represent $A$ in the Uhlenbeck gauges on $B_1, B_2$ respectively. On $B_1 \cap B_2$, they are related by a gauge transformation $g$, i.e., $a = gbg^{-1} - (dg)g^{-1}$ on $B_1 \cap B_2$. Then \\

(1) there exists a constant $C>0$ which doesn't depend on $A$, such that any two points $S_1, S_2 \in V$,
$$|g(S_1) - g(S_2)| \leq C \norm{F_A}^2_{\IL^2(B_1 \cup B_2)};$$

(2) for any point $S \in V$, 
$$\int_{B_1 \cap B_2} (|g - g(V)|^2 + |\nabla g|^2 + |\nabla \nabla g|^2) d^3x \leq C \norm{F_A}^2_{\IL^2(B_1 \cup B_2)}.$$
\end{lemma}

\begin{proof}
Consider the identity $d^*dg = a\cdot ag - 2ag \cdot b + gb\cdot b$.\\

Let $L = a\cdot ag - 2ag\cdot b + gb\cdot b$, then for some constant $C_0 > 0$, 
$$\norm{L}_{\IL^2(B_1 \cup B_2)} \leq C_0 (\norm{a}_{\IL^4_{B_1 \cup B_2}} + \norm{b}_{\IL^4(B_1 \cup B_2)})^2 \leq C_0^2 \norm{F_A}^2_{\IL^2(B_1 \cup B_2)}. $$
The second inequality above is from the Uhlenbeck gauge condition and Sobolev embedding.\\

On $B_1 \cap B_2$,
$$d^*dg = L, ~~~ \text{and} \dfrac{\partial g}{\partial \Vec{n}} = 0,$$
where $\Vec{n}$ is the outside normal vector of $\partial (B_1 \cap B_2)$.\\
So (2) is a direct corollary of the standard elliptic regularity theorem.\\

To prove (1), supposing $G(x, y)$ is the Green's function on $B_1 \cap B_2$ with Neumann boundary condition, then for any point $x$ on $V$,
$$g(x) = \int_{B_1 \cap B_2} L(y)G(x, y)dy + c, ~~\text{where $c$ is a constant that doesn't depend on $x$.}$$
$$|g(S_1) - g(S_2)| \leq \int_{B_1 \cap B_2} |L(y)||G(S_1, y) - G(S_2, y)|d^3y.$$

Note that for any fixed $x$, the function $G(x,y)$ has only one singularity at $y = x$ which $\sim \dfrac{1}{\text{dist}(x, \cdot)}$, where dist$(x, \cdot)$ is the distance to the point $x$, and $|L|$ can be bounded by $2(|a|^2 + |b|^2)$. Hence applying a Hardy's inequality (see theorem \ref{Appendix8} in the appendix) it can be shown that 
$$ |g(S_1) - g(S_2)| \leq C_1 \int_{B_1 \cap B_2} (|\nabla |a||^2 + |\nabla |b||^2 + |a|^2 + |b|^2)d^3x \leq C \norm{F_A}_{\IL^2(B_1 \cup B_2)}^2, $$
where $C$ and $C_1$ are constants. 
\end{proof}

\begin{remark}
When rescaling $B_1 \cup B_2$ by a factor $R$, the above inequalities should be\\

(1) $$|g(S_1) - g(S_2)| \leq C R \norm{F_A}^2_{\IL^2(B_1 \cup B_2)};$$

(2)  $$\int_{B_1 \cap B_2} (\dfrac{1}{R^3}|g - g(V)|^2 + \dfrac{1}{R}|\nabla g|^2 + R|\nabla \nabla g|^2) d^3x \leq C R\norm{F_A}^2_{\IL^2(B_1 \cup B_2)},$$ where $C$ is independent with $R$.
\end{remark}

\begin{lemma}\label{Uhlenbeck gluing lemma}
(gluing Uhlenbeck gauges) Let $U_1$ and $U_2$ be two bounded open sets in $\IR^3$ with smooth boundary that have a nonempty connected overlap. Suppose the closure of $U_1 \backslash U_2$ and the closure of $U_2 \backslash U_1$ do not overlap. Suppose $\chi$ is a smooth cut-off function which is $0$ on $U_1 \backslash U_2$ and $1$ on $U_2 \backslash U_1$ and $|\nabla \chi|$, $|\nabla \nabla \chi|$ are bounded.  Suppose $A$ is an $su(2)$-valued 1-form  on $U_1 \cup U_2$, viewed as a connection on the trivial $SU(2)$-bundle. Suppose $U$ is an open set that covers $U_1 \cup U_2$ and suppose the curvature $F_A$ has small enough $\IL^2$ norm on $U$. Let $u_1$ and $u_2$ be $SU(2)$-valued functions on $U_1$ and $U_2$ respectively, viewed as gauge transformations. Suppose there exists a constant $C > 0$ such that $$|u_1^{-1}u_2 - 1| \leq C 
\int_{U} |F_A|^2 d^3x ~\text{on}~ U_1 \cap U_2 ~ \text{pointwise}.$$ 
$ \displaystyle \int_{U_1 \cap U_2} (|u_1^{-1}u_2 - 1|^2 + |\nabla (u_1^{-1}u_2)|^2 + |\nabla \nabla (u_1^{-1} u_2)|^2) d^3x \leq C \int_{U} |F_A|^2 d^3x$.\\

Let $a_1 = u_1^{-1}Au_1 + u_1^{-1} du_1, a_2 = u_2^{-1}Au_2 + u_2^{-1}du_2$. Assume
$$ \int_{U_i} (|\nabla a_i|^2 + |a_i|^2) d^3x \leq C \int_{U} |F_A|^2 d^3x, ~~ \text{for $i=1, 2$}. $$
Then there exists an $SU(2)$-valued function $u$ on $U_1 \cup U_2$ (still viewed as a gauge transformation) such that letting $a = u^{-1}Au + u^{-1}du$, then\\

(1) $u = u_1$ on $U_1 \backslash U_2$ and $u = u_2$ on $U_2 \backslash U_1$;\\

(2) $ \displaystyle |u^{-1}u_i - 1| \leq C \int_{U} |F_A|^2 d^3x$ on $U_1 \cap U_2$ for $i = 1,2$;\\

(3) and there exists a constant $C'$ such that
$$ \displaystyle \int_{U_1 \cup U_2} (|\nabla a|^2 + |a|^2) d^3x \leq  C' \displaystyle \int_{U} |F_A|^2 d^3x. $$
\end{lemma}

\begin{proof}
Suppose $e^v = u_1^{-1}u_2$ on $U_1 \cap U_2$.  Letting $u = u_1 e^{\chi v}$ on $U_1$ and $u = u_2$ on $U_2 \backslash U_1$, then obviously $u$ satisfies (1) (2). To prove (3), it is enough to show that 

\begin{equation*}
\displaystyle \int_{U_1 \cap U_2} (|\nabla a|^2 + |a|^2) d^3x \leq C' \displaystyle \int_{U} |F_A|^2 d^3x.
\end{equation*}
\begin{equation*}
    \begin{array}{ll}
         ~ & ~~~~ \displaystyle \int_{U_1 \cap U_2} (|\nabla a|^2 + |a|^2) d^3x\\
         & = \displaystyle \int_{U_1 \cap U_2} (|\nabla (e^{-\chi v}a_1 e^{\chi v} + e^{-\chi v} d(e^{\chi v}))|^2 + |e^{-\chi v}a_1 e^{\chi v} + e^{-\chi v}d(e^{\chi v})|^2) d^3x \\ & \leq  \displaystyle 5\int_{U_1 \cap U_2} (|\nabla (e^{-\chi v})|^2|a_1|^2 + |\nabla a_1|^2 + |a_1|^2 |\nabla (e^{\chi v})|^2 + |\nabla \nabla(e^{\chi v})|^2 + |\nabla(e^{-\chi v})|^2 |\nabla (e^{\chi v})|^2) d^3x \\
         & ~~~~ +\displaystyle 2\int_{U_1 \cap U_2} (|a_1|^2 + |\nabla(e^{\chi v})|^2) d^3x \\
         & \leq  \displaystyle C_1\int_{U_1 \cap U_2} ( |\nabla a_1|^2 + (|a_1|^2 + |\nabla \chi|^2|v|^2 + |\nabla v|^2) (|\nabla \chi|^2|v|^2 + |\nabla v|^2) + |\nabla \nabla \chi|^2|v|^2 + |\nabla \nabla v|^2  \\
         & ~~~~ + |a_1|^2 + |\nabla\chi|^2 |v|^2 + |\nabla v|^2) d^3x \\
          & \leq  \displaystyle C_2\int_{U_1 \cap U_2} ( |\nabla a_1|^2 + (|a_1|^2 + |v|^2 + |\nabla v|^2) (|v|^2 + |\nabla v|^2) +|v|^2 + |\nabla \nabla v|^2 + |a_1|^2 +  |\nabla v|^2) d^3x \\
         & \leq  C' \displaystyle \int_{U} |F_A|^2 d^3x,
    \end{array}
\end{equation*}
where $C_1, C_2, C'$ are constants.
\end{proof}

\begin{remark}
When re-scaling $U_1 \cup U_2$ and $U$ by a factor $R$, the above theorem should be stated as following : Assume $$|u_1^{-1}u_2 - 1| \leq C R
\int_{U} |F_A|^2 d^3x ~\text{on}~ U_1 \cap U_2 ~ \text{pointwise}.$$
$ \displaystyle \int_{U_1 \cap U_2} (\dfrac{1}{R^3}|u^{-1}u_i - 1|^2 + \dfrac{1}{R}|\nabla (u_1^{-1}u_2)|^2 + R|\nabla \nabla (u_1^{-1} u_2)|^2) d^3x \leq CR \int_{U} |F_A|^2 d^3x$.\\
Let $a_1 = u_1^{-1}Au_1 + u_1^{-1} du_1, a_2 = u_2^{-1}Au_2 + u_2^{-1}du_2$. Assume
$$ \int_{U_i} (R|\nabla a_i|^2 + \dfrac{1}{R}|a_i|^2) d^3x \leq C R\int_{U} |F_A|^2 d^3x, ~~ \text{for $i=1, 2$}. $$
Then there exists an $SU(2)$-valued function $u$ on $U_1 \cup U_2$ (still viewed as a gauge transformation) such that letting $a = u^{-1}Au + u^{-1}du$, then\\

(1) $u = u_1$ on $U_1 \backslash U_2$ and $u = u_2$ on $U_2 \backslash U_1$;\\

(2) $ \displaystyle |u^{-1}u_i - 1| \leq CR \int_{U} |F_A|^2 d^3x$ on $U_1 \cap U_2$ for $i = 1,2$;\\
then
$$ \displaystyle \int_{U_1 \cup U_2} (R|\nabla a|^2 + \dfrac{1}{R}|a|^2) d^3x \leq  C' \displaystyle R\int_{U} |F_A|^2 d^3x. $$
\end{remark}

\begin{proof}
\begin{equation*}
    \begin{array}{ll}
         ~ & ~~~~ \displaystyle \int_{U_1 \cap U_2} (R|\nabla a|^2 + \dfrac{1}{R}|a|^2) d^3x\\
         & = \displaystyle \int_{U_1 \cap U_2} (R|\nabla (e^{-\chi v}a_1 e^{\chi v} + e^{-\chi v} d(e^{\chi v}))|^2 + \dfrac{1}{R}|e^{-\chi v}a_1 e^{\chi v} + e^{-\chi v}d(e^{\chi v})|^2) d^3x \\
         & \leq  \displaystyle 100\int_{U_1 \cap U_2} ( R|\nabla a_1|^2 + R|a_1|^2 |\nabla (\chi v)|^2 + R|\nabla \nabla(\chi v)|^2 + R|\nabla(\chi v)|^4) d^3x \\
         & ~~~~ +\displaystyle 100\int_{U_1 \cap U_2} (\dfrac{1}{R}|a_1|^2 + \dfrac{1}{R}|\nabla(\chi v)|^2) d^3x \\
         & \leq  \displaystyle C_1\int_{U_1 \cap U_2} ( R|\nabla a_1|^2 + R(|a_1|^2 + \dfrac{1}{R^2}|v|^2 + |\nabla v|^2) (\dfrac{1}{R^2}|v|^2 + |\nabla v|^2) + \dfrac{1}{R^3}|v|^2 + R|\nabla \nabla v|^2 ) d^3x \\
         & ~~~~ +\displaystyle C_1\int_{U_1 \cap U_2} (\dfrac{1}{R}|a_1|^2 + \dfrac{1}{R^3} |v|^2 + \dfrac{1}{R}|\nabla v|^2) d^3x \\
         & \leq  C' R\displaystyle \int_{U} |F_A|^2 d^3x,
    \end{array}
\end{equation*}
where $C_1, C'$ are constants. 
\end{proof}

\subsection{Open covering}
Suppose $F_A \in \IL^2_{ \epsilon}$,
$$\int_{\IR^3 \backslash K} \rho_{\epsilon} |F_A|^2 d^3x < +\infty. $$

Consider an open covering $\{ U_{\alpha} \}$ of $N_{\epsilon} \backslash K$, where each $\alpha$ represents a triple $(n_{\rho}, n_{\theta}, n_s)$. To be more precise, let $t$ be a large enough integer such that $2^t \epsilon \gg l$ and $t \gg 1000$, then $n_{\rho}, n_{\theta}, n_s$ are nonnegative integers such that $n_{\rho} \in [2, +\infty), n_{\theta} \in [0, t-1], n_s \in [0, 2^{t+n_{\rho}}-1]$.  Define $U_{\alpha} = U_{n_{\rho}, n_{\theta}, n_s}$ as the following: $U_{n_{\rho}, n_{\theta}, n_s} = \{ (r, \theta, s) \in N_{\epsilon} |$
$$~ \rho \in (\dfrac{\epsilon}{2^{(3n_{\rho} + 1)}}, \dfrac{\epsilon}{2^{(3n_{\rho}-4)}}), \theta \in (\dfrac{2\pi (3n_{\theta} -4)}{3t + 2}, \dfrac{2\pi (3n_{\theta} + 1)}{3t + 2}), s \in (\dfrac{(3n_s -4)l}{3 \times 2^{(t+n_{\rho} )}+2}, \dfrac{(3n_s + 1)l}{3 \times 2^{(t+n_{\rho})}+2}) \}. $$
where the value of $\theta$ is defined modulo $2\pi$ and the value of $s$ modulo $l$, the subscripts $n_{\theta}$ and $n_s$ are also defined in the sence of modulo $3t + 2$ and modulo $3 \times 2^{t + n_{\rho}} + 2$ respectively.\\

 Let the point $P_{\alpha}$ be
$$P_{\alpha} =P_{n_{\rho}, n_{\theta}, n_s} = \{ r = \dfrac{\epsilon}{2^{3n_{\rho}}}, \theta = \dfrac{6\pi n_{\theta}}{3t+2}, s = \dfrac{3n_s l}{3 \times 2^{(t+n_{\rho} )}+2} \}.$$

Two open sets or points are called ``neighbors" if exactly one of the three subscripts differ by $1$.\\

The open cover $U_{\alpha}$ is a sub open cover of open balls $B_{\alpha} = B_{n_{\rho}, n_{\theta}, n_s}$, such that  the ball $B_{n_{\rho},n_{\theta},n_s}$ has radius $2^{-n_{\rho}} \eta$, where $\eta$ is a constant. Each ball $B_{n_{\rho},n_{\theta},n_s}$ is a subset of $N_{2^{-(n_{\rho} - 10)}\epsilon} \backslash N_{2^{-(n_{\rho} + 10)} \epsilon}$ but covers $U_{n_{\rho}, n_{\theta}, n_s}$. \\

From lemma \ref{Uhlenbeck lemma}, assuming $\epsilon$ is small enough, on each ball $B_{\alpha}$, there exists an $SU(2)$-valued function $u_{\alpha}$ such that letting $A_{\alpha} = u_{\alpha}(A)$ on $B_{\alpha}$, then
\begin{equation*}
    \left\{\begin{array}{ll}
         d^*(A_{\alpha}) = 0, ~~~ *A_{\alpha}|_{\partial B_{\alpha}} = 0;&\\
        \displaystyle \int_{B_{{\alpha}}} |\nabla A_{\alpha}|^2 + \dfrac{1}{\rho^2}|A_{\alpha}|^2 d^3x \leq C \int_{B_{\alpha}} |F_A|^2 d^3 x. & 
    \end{array}
    \right.
\end{equation*}
 Using the freedom of a constant gauge transformation, it may be assumed that
 \begin{itemize}
     \item $u_{{n_{\rho},n_{\theta},n_s}}(P_{{n_{\rho},n_{\theta},n_s}}) = u_{{n_{\rho},(n_{\theta} + 1),n_s}}(P_{{n_{\rho},n_{\theta},n_s}})$ for all $n_{\theta} = 1, 2, \cdots, t-1$ (note that $0$ is excluded);
     \item $u_{{n_{\rho},1,n_s}}(P_{{n_{\rho},0,n_s}}) = u_{{n_{\rho} ,1 , (n_s + 1)}}(P_{{n_{\rho},0 ,n_s}})$ for all $n_s = 1, 2, 3, \cdots, 2^{t+n_{\rho}}-1$ (note that $0$ is excluded again);
     \item $u_{{n_{\rho},1,1}}(P_{{n_{\rho},0,0}}) = u_{{(n_{\rho} + 1) ,1, 1}}(P_{{n_{\rho},0 ,0}})$ for all $n_{\rho} = 1, 2, 3, \cdots;$
     \item Letting $g_{n_{\rho}, n_s} = u_{n_{\rho}, 0, n_s}^{-1}u_{n_{\rho}, 1, n_s}$, then $g_{n_{\rho}, n_s} (P_{n_{\rho}, 1 ,n_s})$ is called the \textbf{ $\theta$ monodromy of $A$ at the level $(n_{\rho}, n_s)$}, denoted as $\Gamma^{\theta}_{n_{\rho}, n_s}$;
     \item Letting $g_{n_{\rho}} = u_{n_{\rho}, 1, 0}^{-1}u_{n_{\rho}, 1, 1}$, then $g_{n_{\rho}} (P_{n_{\rho}, 1 , 1})$ is called the \textbf{$s$ monodromy of $A$ at the level $n_{\rho}$}, denoted as $\Gamma^{s}_{n_{\rho}}$.
 \end{itemize}
 
Note that even all the above conditions are satisfied, there is still one more freedom to choose the Uhlenbeck gauges (change all the above $u$ by a same constant gauge transformation), which will be used later.
 
\begin{corollary}
Suppose $U_{n_{\rho}, n_s} = \union\limits_{n_{\theta}} U_{n_{\rho}, n_{\theta}, n_s}$, $U_{n_{\rho}} = \union\limits_{n_s}U_{n_{\rho}, n_s}$. There exists a constant $C> 0$. 
\begin{itemize}
    \item Suppose $\alpha = (n_{\rho}, n_{\theta}, n_s)$, $\alpha'$ is either $(n_{\rho} + 1, n_{\theta}, n_s)$, $(n_{\rho}, n_{\theta}+1, n_s)$ (in wich case $n_{\theta}\neq 0$) or $(n_{\rho}, n_{\theta}, n_s)$ (in which case $n_s \neq 0$). Suppose $P \in U_{\alpha} \cap U_{\alpha'}$ is any point in the overlap. Then 
$$|u_{\alpha}(P) - u_{\alpha'}(P)|  \leq C\int_{U_{n_{\rho}} \cup U_{n_{\rho} + 1}} \rho |F_A|^2 d^3x.$$ 
    \item Suppose $\alpha = (n_{\rho}, 0, n_s)$, $\alpha' = (n_{\rho}, 1, n_s)$. Suppose $P \in U_{\alpha} \cap U_{\alpha'}$ is any point in the overlap. Then
    $$|u_{\alpha}(P) \Gamma^{\theta}_{n_{\rho}, n_s} - u_{\alpha'}(P)|  \leq C\int_{U_{n_{\rho}} \cup U_{n_{\rho} + 1}} \rho |F_A|^2 d^3x.$$ 
    \item Suppose $\alpha = (n_{\rho}, n_{\theta}, 0)$, $\alpha' = (n_{\rho}, n_{\theta}, 1)$. Suppose $P \in U_{\alpha} \cap U_{\alpha'}$ is any point in the overlap. Then
    $$|u_{\alpha}(P) \Gamma^s_{n_{\rho}} - u_{\alpha'}(P)|  \leq C\int_{U_{n_{\rho}} \cup U_{n_{\rho} + 1}} \rho |F_A|^2 d^3x.$$ 
\end{itemize} 
\end{corollary}
 
\begin{proof}
This is a direct corollary of lemma \ref{Uhlenbeck transition estimate lemma}.
\end{proof} 

\begin{corollary}
Assuming $\displaystyle \int_{N_{\epsilon}}\rho |F_A|^2 d^3x \leq +\infty$, then there exist a $\Gamma^{\theta} \in SU(2)$ and a $\Gamma^s \in SU(2)$, such that
$$\lim\limits_{n_{\rho} \rightarrow +\infty} \Gamma^{\theta}_{n_{\rho}, n_s} = \Gamma^{\theta} ~ \text{uniformly in $n_s$}, ~~ \lim\limits_{n_{\rho} \rightarrow +\infty} \Gamma^s_{n_{\rho}} = \Gamma^s. $$
Moreover, $\Gamma^{\theta}$ and $\Gamma^s$ commute with each other.
\end{corollary}

\begin{proof}
Taking $\Gamma^{\theta}_{n_{\rho}, n_s}$ as an example, choosing any point $P_1 \in U_{(n_{\rho} + 1), 0, n_s} \cap U_{n_{(n_{\rho} + 1), 1, n_s}} \cap U_{n_{\rho}, 0,  n_s} \cap U_{n_{\rho}, 1, n_s}$ (note that this set is not empty). Then 
$$| \Gamma^{\theta}_{n_{\rho}, n_s} - u^{-1}_{n_{\rho}, 0,  n_s}(P_1)u_{n_{\rho}, 1, n_s}(P_1)|  \leq C\int_{U_{n_{\rho}} \cup U_{n_{\rho} + 1}} \rho |F_A|^2 d^3x.$$ 
$$| \Gamma^{\theta}_{(n_{\rho}+1), n_s} - u^{-1}_{(n_{\rho} + 1), 0, n_s}(P_1)u_{(n_{\rho} + 1), 1, n_s}(P_1)|  \leq C\int_{U_{n_{\rho + 1}} \cup U_{n_{\rho} + 2}} \rho |F_A|^2 d^3x.$$ 
\begin{equation*}
    \begin{array}{ll}
         ~& ~~~~ |u^{-1}_{n_{\rho}, 0,  n_s}(P_1)u_{n_{\rho}, 1, n_s}(P_1) - u^{-1}_{(n_{\rho} + 1), 0, n_s}(P_1)u_{(n_{\rho} + 1), 1, n_s}(P_1)| \\
         & \leq  |u_{(n_{\rho} + 1), 0, n_s}(P_1)u^{-1}_{n_{\rho}, 0,  n_s}(P_1)- u_{(n_{\rho} + 1), 1, n_s}(P_1)u^{-1}_{n_{\rho}, 1, n_s}(P_1) |\\
         & \leq |u_{(n_{\rho} + 1), 0, n_s}(P_1)u^{-1}_{n_{\rho}, 0,  n_s}(P_1)-1| + | u_{(n_{\rho} + 1), 1, n_s}(P_1)u^{-1}_{n_{\rho}, 1, n_s}(P_1) -1 |  \displaystyle \leq 2C\int_{U_{n_{\rho}} \cup U_{n_{\rho} + 1}} \rho |F_A|^2 d^3x.
    \end{array}
\end{equation*}
Hence 
$$|\Gamma^{\theta}_{n_{\rho}, n_s} - \Gamma^{\theta}_{(n_{\rho} + 1), n_s}| \leq 4 C \int_{U_{n_{\rho}} \cup U_{n_{\rho} + 1} \cup U_{n_{\rho} + 2} } \rho |F_A|^2 d^3x. $$

Since $\displaystyle \int_{N_{\epsilon} \backslash K} \rho |F_A|^2 d^3x < +\infty$, for any $n_s$ the sequence $\Gamma^{\theta}_{n_{\rho}, n_s}$ is a Cauchy sequence , hence there exists an $\Gamma^{\theta}_{n_s} = \lim\limits_{n_{\rho} \rightarrow +\infty} \Gamma^{\theta}_{n_{\rho}, n_s} $. (The convergence is uniform in $u_s$.) On the other hand,
$$| \Gamma^{\theta}_{n_{\rho}, n_s} - u^{-1}_{n_{\rho}, 0,  n_s}(P_1)u_{n_{\rho}, 1, n_s}(P_1)|  \leq C\int_{U_{n_{\rho}} \cup U_{n_{\rho} + 1}} \rho |F_A|^2 d^3x.$$ 
$$| \Gamma^{\theta}_{n_{\rho}, (n_s + 1)} - u^{-1}_{n_{\rho}, 0, (n_s + 1)}(P_1)u_{n_{\rho}, 1, (n_s + 1)}(P_1)|  \leq C\int_{U_{n_{\rho}} \cup U_{n_{\rho} + 1}} \rho |F_A|^2 d^3x.$$ 
\begin{equation*}
    \begin{array}{ll}
         ~& ~~~~ |u^{-1}_{n_{\rho}, 0,  n_s}(P_1)u_{n_{\rho}, 1, n_s}(P_1) - u^{-1}_{n_{\rho} , 0, (n_s + 1)}(P_1)u_{n_{\rho}, 1, (n_s + 1)}(P_1)| \\
         & \leq  |u_{n_{\rho}, 0, (n_s + 1)}(P_1)u^{-1}_{n_{\rho}, 0,  n_s}(P_1)- u_{n_{\rho}, 1, (n_s + 1)}(P_1)u^{-1}_{n_{\rho}, 1, n_s}(P_1) |\\
         & \leq |u_{n_{\rho}, 0, (n_s + 1)}(P_1)u^{-1}_{n_{\rho}, 0,  n_s}(P_1)- 1 | + |u_{n_{\rho}, 1, (n_s + 1)}(P_1)u^{-1}_{n_{\rho}, 1, n_s}(P_1) -1 |  \displaystyle \leq 2C\int_{U_{n_{\rho}} \cup U_{n_{\rho} + 1}} \rho |F_A|^2 d^3x.
    \end{array}
\end{equation*}
So
$$|\Gamma^{\theta}_{n_{\rho}, n_s} - \Gamma^{\theta}_{n_{\rho}, (n_s + 1)}| \leq 4 C \int_{U_{n_{\rho}} \cup U_{n_{\rho} + 1} } \rho |F_A|^2 d^3x. $$
Hence $\Gamma^{\theta}_{n_s}$ are the same for different $n_s$, i.e.,  there exists an $\Gamma^{\theta} = \lim\limits_{n_{\rho} \rightarrow +\infty} \Gamma^{\theta}_{n_{\rho}, n_s}$ which convergent uniformly in $u_s$.\\

The existence of $\Gamma^s$ can be proved by the similar reason and skipped.\\

Finally, to prove $\Gamma^{\theta}$ and $\Gamma^s$ commute, note that for any $P_1 \in U_{n_{\rho}, 0,1} \cup U_{n_{\rho}, 1, 1} \cup U_{n_{\rho},1,0}$, 
$$| \Gamma^{\theta}_{n_{\rho}, 1} - u^{-1}_{n_{\rho}, 0,  1}(P_1)u_{n_{\rho}, 1, 1}(P_1)|  \leq C\int_{U_{n_{\rho}} \cup U_{n_{\rho} + 1}} \rho |F_A|^2 d^3x.$$ 
$$| \Gamma^s_{n_{\rho}} - u^{-1}_{n_{\rho},1, 0}(P_1)u_{n_{\rho}, 1, 1}(P_1)|  \leq C\int_{U_{n_{\rho}} \cup U_{n_{\rho} + 1}} \rho |F_A|^2 d^3x.$$ 
$$| \Gamma^{\theta}_{n_{\rho}, 0} - u^{-1}_{n_{\rho}, 0,  0}(P_1)u_{n_{\rho}, 1, 0}(P_1)|  \leq C\int_{U_{n_{\rho}} \cup U_{n_{\rho} + 1}} \rho |F_A|^2 d^3x.$$ 
$$| \Gamma^s_{n_{\rho}} - u^{-1}_{n_{\rho},0, 0}(P_1)u_{n_{\rho}, 0, 1}(P_1)|  \leq C\int_{U_{n_{\rho}} \cup U_{n_{\rho} + 1}} \rho |F_A|^2 d^3x.$$ 
So
\begin{equation*}
    \begin{array}{ll}
         |\Gamma^{\theta} \Gamma^s - \Gamma^s \Gamma^{\theta}| & = \lim\limits_{n_{\rho} \rightarrow +\infty}|\Gamma^{\theta}_{n_{\rho}, 0} \Gamma^s_{n_{\rho}} - \Gamma^s_{n_{\rho}} \Gamma^{\theta}_{n_{\rho},1}|\\
         & = \lim\limits_{n_{\rho} \rightarrow +\infty}  | u^{-1}_{n_{\rho},0, 0}(P_1)u_{n_{\rho}, 1, 1}(P_1)- u^{-1}_{n_{\rho},0, 0}(P_1)u_{n_{\rho}, 1, 1}(P_1)| = 0. 
    \end{array}
\end{equation*}
\end{proof}

Recall that there is an extra freedom of a constant gauge transformation not be used in the definitions of local Uhlenbeck gauges. This freedom can be used to change $\Gamma^{\theta}$ and $\Gamma^s$ up to conjugacy. In particular, it can be assumed that $\Gamma^{\theta} = e^{2\pi\gamma \sigma}$ and $\Gamma^s = e^{2\pi \tilde{\gamma} \sigma}$ for some $\gamma, \tilde{\gamma} \in (0,1)$.

\subsection{Proof of lemma \ref{local regularity lemma for connection}}

\begin{theorem}
One can adapt all the $u_{\alpha}$ to $\tilde{u}_{\alpha}$ and let $\tilde{A}_{\alpha} = \tilde{u}_{\alpha}(A)$, such that for some possibly larger $\tilde{C}$,
\begin{itemize}
    \item $ \displaystyle \int_{U_{{\alpha}}} |\nabla \tilde{A}_{\alpha}|^2 + \dfrac{1}{\rho^2}|\tilde{A}_{\alpha}|^2 d^3x \leq \tilde{C} \sum\limits_{\alpha' ~\text{is either}~ \alpha~ or~ \text{a neighbor of}~ \alpha}\int_{U_{\alpha'}} |F_A|^2 d^3 x.$
    \item Suppose $\alpha = (n_{\rho}, n_{\theta}, n_s)$, $\alpha'$ is either $(n_{\rho} + 1, n_{\theta}, n_s)$, $(n_{\rho}, n_{\theta}+1, n_s)$ (in wich case $n_{\theta}\neq 0$) or $(n_{\rho}, n_{\theta}, n_s)$ (in which case $n_s \neq 0$), then $\tilde{u}_{\alpha} = \tilde{u}_{\alpha'}$ on $ U_{\alpha} \cap U_{\alpha'}$.
    \item Suppose $\alpha = (n_{\rho}, 0, n_s)$, $\alpha' = (n_{\rho}, 1, n_s)$, then $\tilde{u}_{\alpha} \Gamma^{\theta} = \tilde{u}_{\alpha'} $ on  $ U_{\alpha} \cap U_{\alpha'}$.
    \item Suppose $\alpha = (n_{\rho}, n_{\theta}, 0)$, $\alpha' = (n_{\rho}, n_{\theta}, 1)$, then $ \tilde{u}_{\alpha} \Gamma^s = \tilde{u}_{\alpha'}$ on $ U_{\alpha} \cap U_{\alpha'}$.
\end{itemize}
\end{theorem}

\begin{proof}

Suppose $\chi(t)$ is a smooth cut-off function which is $0$ when $t \leq 0 $ and $1$ when $t \geq 1$. It may also be assumed that $\chi(t) = \chi(1-t)$. Let $\chi(a, b, t) = \chi(\dfrac{t-a}{b-a})$.\\

Let $\widehat{U}_{\alpha} = \union\limits_{\alpha' \text{is a neighbor of} \alpha} U_{\alpha}$. Define the following functions on $\widehat{U}_{n_{\rho}, n_{\theta}, n_s}$ which depend only one of the parameters $\rho, \theta$ and $s$ respectively: \\

1. Let $ \chi^{\rho}_{\alpha} = \chi^{\rho}_{n_{\rho}, n_{\theta}, n_s}(\rho) = \chi(\dfrac{\epsilon}{2^{3n_{\rho} - 1}}, \dfrac{\epsilon}{2^{3n_{\rho} - 4}}, \rho)$ be a function which depends only on $\rho$. For convenience, let $\chi^{\rho + 1}_{\alpha}$ denote $\chi^{\rho}_{(n_{\rho} + 1), n_{\theta}, n_s}$. \\

2. Let  $\chi^{\theta}_{\alpha} = \chi^{\theta}_{n_{\rho}, n_{\theta}, n_s} (\theta) = \chi(\dfrac{2\pi(3n_{\theta} - 4)}{3t+2} , \dfrac{2\pi(3n_{\theta}-2)}{3t+2} ,\theta)$ be a function which depends only on $\theta$. For convenience, let $\chi^{\theta + 1}_{\alpha}$ denote $\chi^{\theta}_{n_{\rho}, (n_{\theta} + 1), n_s}$. \\

3. Let $\chi^{s}_{\alpha} = \chi^{s}_{n_{\rho}, n_{\theta}, n_s}(s) = \chi(\dfrac{(3n_{s} - 4)l}{3 \times 2^{t+n_{\rho} + 2}}, \dfrac{(3n_{s} - 2)l}{3 \times 2^{t+n_{\rho} + 2}},s)$ be a function which depends only on $s$. For convenience, let $\chi^{s + 1}_{\alpha}$ denote $\chi^{s}_{n_{\rho}, n_{\theta}, (n_s + 1)}$. \\

There are three steps to construct $\tilde{u}_{\alpha}$:\\

Step 1: Let $v^{\theta}_{\alpha} = v^{\theta}_{n_{\rho}, n_{\theta}, n_s}$ be the $su(2)$-valued function on $U_{n_{\rho}, (n_{\theta} - 1), n_s} \cap U_{n_{\rho}, n_{\theta}, n_s} $ such that $e^{v^{\theta}_{\alpha}} =  \left\{
\begin{array}{ll}
    u^{-1}_{n_{\rho}, (n_{\theta} - 1), n_s} u_{n_{\rho}, n_{\theta}, n_s} & \text{when} ~~~~ n_{\theta} \neq 1 \\
     (u_{n_{\rho}, (n_{\theta} - 1), n_s} \Gamma^{\theta})^{-1} u_{n_{\rho}, n_{\theta}, n_s} & \text{when} ~~~~ n_{\theta} = 1 
\end{array}
\right.$\\
and let $v^{\theta + 1}_{n_{\rho}, n_{\theta}, n_s}$ be the $su(2)$-valued function on $U_{n_{\rho}, n_{\theta}, n_s} \cap U_{n_{\rho}, (n_{\theta} + 1), n_s}$ such that $e^{v^{\theta + 1}_{\alpha}} = \left\{
\begin{array}{ll}
    u^{-1}_{n_{\rho}, n_{\theta}, n_s} u_{n_{\rho}, (n_{\theta} + 1), n_s} & \text{when} ~~~~ n_{\theta} \neq 0 \\
    (u_{n_{\rho}, n_{\theta}, n_s} \Gamma^{\theta})^{-1} u_{n_{\rho}, (n_{\theta} + 1), n_s} & \text{when} ~~~~ n_{\theta} = 0
\end{array}
\right.$.\\

Let $$\hat{u}_{\alpha} = \left\{
\begin{array}{ll}
      u_{\alpha} e^{-\chi^{\theta}_{\alpha} v^{\theta}_{\alpha}} & \text{on} ~~~~U_{n_{\rho}, (n_{\theta} - 1), n_s} \cap U_{n_{\rho}, n_{\theta}, n_s}; \\
      e^{\chi^{\theta + 1}_{\alpha} v^{\theta + 1}_{\alpha}} u_{\alpha} & \text{on}~~~~  U_{n_{\rho}, n_{\theta}, n_s} \cap U_{n_{\rho}, (n_{\theta} + 1), n_s}; \\ 
      u_{\alpha} & \text{elsewhere on} ~~~~ U_{n_{\rho}, n_{\theta}, n_s}.
\end{array}
\right.$$

2. Let $v^s_{\alpha} = v^s_{n_{\rho}, n_{\theta}, n_s}$ be the $su(2)$-valued function on $U_{n_{\rho}, n_{\theta}, (n_s - 1)} \cap U_{n_{\rho}, n_{\theta}, n_s} $ such that $e^{v^{s}_{\alpha}} = \left\{ \begin{array}{ll}
    u^{-1}_{n_{\rho}, n_{\theta}, (n_s - 1)} u_{n_{\rho}, n_{\theta}, n_s} & \text{when} ~~~~ n_s \neq 1 \\
    (u_{n_{\rho}, n_{\theta}, (n_s - 1)} \Gamma^s)^{-1} u_{n_{\rho}, n_{\theta}, n_s} & \text{when} ~~~~ n_s = 1
\end{array}
\right.$\\
and let $v^{s + 1}_{n_{\rho}, n_{\theta}, n_s}$ be the $su(2)$-valued function on $U_{n_{\rho}, n_{\theta}, n_s} \cap U_{n_{\rho}, n_{\theta}, (n_s + 1)}$ such that $e^{v^{s + 1}_{\alpha}} = \left\{ \begin{array}{ll}
    u^{-1}_{n_{\rho}, n_{\theta}, n_s } u_{n_{\rho}, n_{\theta}, (n_s + 1)} & \text{when} ~~~~ n_s \neq 0 \\
     (u_{n_{\rho}, n_{\theta}, n_s }\Gamma^s)^{-1} u_{n_{\rho}, n_{\theta}, (n_s + 1)} & \text{when} ~~~~ n_s = 0
\end{array}
\right.$.\\

Let $$\hat{\hat{u}}_{\alpha} = \left\{
\begin{array}{ll}
      \hat{u}_{\alpha} e^{-\chi^{s}_{\alpha} v^{s}_{\alpha}} & \text{on} ~~~~U_{n_{\rho}, n_{\theta}, (n_s - 1)} \cap U_{n_{\rho}, n_{\theta}, n_s}; \\ 
        e^{\chi^{s + 1}_{\alpha} v^{s + 1}_{\alpha}} \hat{u}_{\alpha} & \text{on}~~~~  U_{n_{\rho}, n_{\theta}, n_s} \cap U_{n_{\rho}, n_{\theta}, (n_s + 1)}; \\ 
      \hat{u}_{\alpha} & \text{elsewhere on} ~~~~ U_{n_{\rho}, n_{\theta}, n_s}.
\end{array}
\right.$$

3. Let $v^{\rho}_{\alpha} = v^{\rho}_{n_{\rho}, n_{\theta}, n_s}$ be the $su(2)$-valued function on $U_{(n_{\rho} - 1), n_{\theta}, n_s} \cap U_{n_{\rho}, n_{\theta}, n_s}$ such that $e^{v^{\rho}_{\alpha}} = u^{-1}_{(n_{\rho} - 1), n_{\theta}, n_s} u_{n_{\rho}, n_{\theta}, n_s}$ and let $v^{\rho + 1}_{n_{\rho}, n_{\theta}, n_s}$ be the $su(2)$-valued function on $U_{n_{\rho}, n_{\theta}, n_s} \cap U_{(n_{\rho} + 1), n_{\theta}, n_s}$ such that $e^{v^{\rho + 1}_{\alpha}} = u^{-1}_{n_{\rho}, n_{\theta}, n_s} u_{(n_{\rho} + 1),n_{\theta}, n_s}$. \\

Let $\tilde{u}_{\alpha} = \left\{
\begin{array}{ll}
     \hat{\hat{u}}_{\alpha} e^{-\chi^{\rho}_{\alpha} v^{\rho}_{\alpha}} & \text{on} ~~~~U_{(n_{\rho} - 1), n_{\theta}, n_s} \cap U_{n_{\rho}, n_{\theta}, n_s}; \\ 
     e^{\chi^{\rho + 1}_{\alpha} v^{\rho + 1}_{\alpha}} \hat{\hat{u}}_{\alpha} & \text{on}~~~~  U_{n_{\rho}, n_{\theta}, n_s} \cap U_{(n_{\rho} + 1), n_{\theta}, n_s}; \\ 
     \hat{\hat{u}}_{\alpha} & \text{elsewhere on} ~~~~ U_{n_{\rho}, n_{\theta}, n_s}.
\end{array}
\right.$\\

Applying the argument of lemma \ref{Uhlenbeck gluing lemma} twice on each of the above steps, it is then proved that $\tilde{u}_{\alpha}$ constructed above satisfies all the requirements.

\end{proof}

Let $A^f = \gamma \sigma d\theta + \dfrac{\tilde{\gamma}}{2\pi l} \sigma ds$. It is flat and trivial on each open ball $B_{n_{\rho}, n_{\theta}, n_s}$. However $A^f$ is not trivial globally. One way to describe $A^f$ is to set $A^f$ to be the trivial connection on each ball $B_{n_{\rho}, n_{\theta}, n_s}$ but with a nontrivial constant gauge transformation on $B_{n_{\rho}, 0, n_s} \cap B_{n_{\rho}, 1, n_s}$, namely $\widehat{g}_{n_{\rho}, n_s} = e^{2\pi\gamma \sigma}$ and a nontrivial constant gauge transformation on $B_{n_{\rho}, 1, 0} \cap B_{n_{\rho}, 1, 1}$.\\

Clearly in the above gauge, on each $U_{n_{\rho}, n_{\theta}, n_s}$, letting $A^f_{\alpha} = 0$, then
$$  \int_{U_{{\alpha}}} |\nabla (\tilde{A}_{\alpha} - A^f_{\alpha})|^2 + \dfrac{1}{\rho^2}|\tilde{A}_{\alpha} - A^f_{\alpha}|^2 d^3x \leq \tilde{C} \sum\limits_{\alpha' ~\text{is either}~ \alpha~ or~ \text{a neighbor of}~ \alpha}\int_{U_{\alpha'}} |F_A|^2 d^3 x.$$
Or equivalently, for some possibly larger $\tilde{C}$,
$$  \int_{U_{{\alpha}}} (\rho|\nabla^{A^f_{\alpha}} (\tilde{A}_{\alpha} - A^f_{\alpha})|^2 + \dfrac{1}{\rho}|\tilde{A}_{\alpha} - A^f_{\alpha}|^2 )d^3x \leq \tilde{C} \sum\limits_{\alpha' ~\text{is either}~ \alpha~ or~ \text{a neighbor of}~ \alpha}\int_{U_{\alpha'}} \rho |F_A|^2 d^3 x.$$
The above expression is gauge invariant. In particular, on each $U_{\alpha}$, there exists a gauge tranformation $v_{\alpha}$ that transforms $A^f_{\alpha}$ back to $A^f =  \gamma \sigma d\theta + \dfrac{\tilde{\gamma}}{2\pi l} \sigma ds$ and also make $v_{\alpha} \circ \tilde{u}_{\alpha}$ a global gauge transformation (i.e., they agree on all intersections of different $U_{\alpha}$). Letting $\tilde{\tilde{A}} = v_{\alpha} (\tilde{A}_{\alpha})$ (it is defined globally, hence the subscript is not needed), then on each $U_{\alpha}$,
$$  \int_{U_{{\alpha}}} (\rho|\nabla^{A^f} (\tilde{\tilde{A}} - A^f)|^2 + \dfrac{1}{\rho}|\tilde{\tilde{A}} - A^f|^2 )d^3x \leq \tilde{C} \sum\limits_{\alpha' ~\text{is either}~ \alpha~ or~ \text{a neighbor of}~ \alpha}\int_{U_{\alpha'}} \rho |F_A|^2 d^3 x.$$

Summing over all the $U_{\alpha}$, for some possibly larger $\tilde{C}$,
$$  \int_{N_{\epsilon}} (\rho|\nabla^{A^f} (\tilde{\tilde{A}} - A^f)|^2 + \dfrac{1}{\rho}|\tilde{\tilde{A}} - A^f|^2 )d^3x \leq \tilde{C} \int_{N_{\epsilon}} \rho |F_A|^2 d^3 x.$$

\section{Limit behaviors near the knot}\label{Section: Limit behaviors near the knot}

In this section, we assume $\Psi = A + \Phi$ is a smooth solution to the Bogomonly equations on $\IR^3 \backslash K$, such that the curvature $F_A \in \IL$. Recall that $\chi$ is the cut-off function which is $1$ on $N_{\epsilon}$ and $0$ on $\IR^3
\backslash N_{2\epsilon}$. From the previous section, we may assume that after a gauge transformation, $\chi(\Psi - \Psi_{\gamma, M, 0}) \in \IH$. In this section, we assume $\gamma \neq \dfrac{1}{4}$. And we derive more detailed limit behaviors of $\Psi$ near the knot.

\subsection{The poly-homogeneous expansion}
Let $\alpha = \min\{2\gamma, 1 - 2\gamma\} \in (0, \dfrac{1}{2})$. 

\begin{conjecture}\label{the conjecture}

If $\Psi$ is a solution to the Bogomolny equations on $\IR^3 \backslash K$ such that 
$$\Psi - \Psi_{\gamma, M, k} \in \IH $$
for some $\gamma \neq \dfrac{1}{4}$. Then
$$\Psi \sim \Psi_{\gamma, M, k} + \sum\limits_{\lambda} \sum\limits_{p = 0}^{n_{\lambda}} \psi_{\lambda, p}\rho^{\lambda}(\log \rho)^p .$$
Here $\lambda$ is taken from a discrete subset of $\IR$ starting from $- \alpha$. Each $n_{\lambda}$ is a non-negative integer. In addition, when $\lambda = -\alpha$, $n_{\lambda} = 0$. 
All $\psi_{\lambda, p}$ are smooth in $\theta$ and $s$.  In particular for some small enough $\delta > 0$, $$\psi = \rho^{-\alpha}\psi_{\alpha} + O(\rho^{- \alpha + \delta}).$$ 
\end{conjecture}

Here the notation $\sim$ means that the difference between $\Psi$ and any partial sum from the right hand side is $O(\rho^{\lambda}(\log)^p)$, where  $\psi_{\lambda, p}\rho^{\lambda}(\log \rho)^p$ is the next smaller term in the summation. All the derivatives of $\Psi$ has also similar expansions.\\

It is likely that conjecture \ref{the conjecture} is true based on the theory of micro-local analysis. See for example \cite{Mazzeo2013TheCondition}. However, the author can only prove a weaker statement in this paper:

\begin{proposition}\label{proposition of polyhomogeneous expansion}
\begin{itemize}
    \item There is a configuration $\psi_{-\alpha} $in $\theta, s$ which is smooth in $\theta$ direction and lies in $L^2(d\theta ds)$, such that in a neighbourhood of the knot,
    $$\Psi - \Psi_{\gamma, M, 0} = \psi_{-\alpha}\rho^{-\alpha} + \psi', $$
    where $\psi'$ is a configuration such that 
    $$\lim\limits_{\rho \rightarrow 0} \rho^{-\alpha - \delta}\int_0^{l}\int_0^{2\pi}|\psi'|^2 d\theta ds = 0,$$
    where $\delta > 0$ is a small enough number.
    \item When $0 \leq \alpha < \dfrac{1}{6}$, $\psi_{-\alpha}$ is also smooth in $s$ direction.
\end{itemize}
\end{proposition}

In another word, the author can only get the leading term in the poly-homogeneous expansion. It seems that when $0 \leq \alpha < \dfrac{1}{6}$, one can use the same method as in this paper to prove the full poly-homogeneous expansion. However, this part is not rigorously written down by the author. And when $\dfrac{1}{6} < \alpha < \dfrac{1}{2}$ it doesn't seem to work.

\begin{corollary}\label{Corollary of gamma is 0 implies smooth}
Assume $\gamma = 0$.  If $\Psi$ is a smooth solution to the Bogomolny equations on $\IR^3 \backslash K$ such that 
$$\Psi - \Psi_{0, M, k} \in \IH.$$
Then $\Psi$ extends to a global smooth solution to the Bogomolny equations on $\IR^3$.
\end{corollary}

\begin{proof}
From the previous proposition, when $\gamma = 0$, $\alpha = 0$, then near the knot 
$$\Psi - \Psi_{0, M, k} = \psi_0 + O(\rho^{\delta}), ~~~ \nabla(\Psi - \Psi_{0, M, k}) = \nabla \psi_0 +  O(\rho^{-1 + \delta}),$$
where $\psi_0$ is smooth term which is independent with $\rho$, $\delta$ is a small enough positive constant.\\

Using the equation $\tilde{N}(\psi_0) = O(\rho^{-\delta})$ for some small enough positive constant $\delta$, we know that $\psi_0$ is also independent with $\theta$. Thus
$$\nabla \psi_0 = (\partial_s \psi_0) ds = O(1). $$
Note that $\Psi_{0, M, k}$ is a well-defined smooth configuration on the entire $\IR^3$. So, $|\nabla \Psi| = O(\rho^{-1 + \delta}) + O(1) = O(\rho^{-1 + \delta})$, which implies the curvature $|F_A| = O(\rho^{-1 + \delta})$. So in fact the curvature $|F_A|$ is in the ordinary $L^2$ norm $L^2(\IR^3)$. Thus classical theory on the Bogomolny equations say, see for example \cite{Jaffe1980VorticesTheories}, that $\Psi$ is gauge equivalent to a smooth solution to the Bogomolny equations on $\IR^3$.

\end{proof}

The remaining part of this subsection is only parenthetical. It explains why this poly-homogeneous expansion is likely to be true from the viewpoint of micro-local analysis, or more precisely, Mazzeo's theory in \cite{Mazzeo1991EllipticI}. Readers may skip it. Starting from  the next subsection, our proof of proposition \ref{proposition of polyhomogeneous expansion} doesn't rely on any background knowledge in micro-local analysis.\\

The operator $\tilde{L}_{\Psi}$ can be viewed as a elliptic edge calculator which is elliptic in the interior of $\IR^3 \backslash K$ and has an ``edge elliptic" type near two boundaries: approaching the knot and approaching infinity.

\paragraph{Near the knot}~\\

Treat the knot as a boundary of the manifold $\IR^3 \backslash K$. Then the normal operator of $\tilde{L}_{\Psi}$ near the knot is represented by $\tilde{N}$. The general definition of the normal operator can be found in much literature, see \cite{Mazzeo1991EllipticI} for example. But in our situation, the operator $\tilde{N}$ is similar with the operator $N$ introduced in section \ref{Section: The Fredholm theory near the knot}. (They only differ by a bounded zero order term, which is not important for the sake of doing analysis.)\\

To be more precise, the normal operator $\tilde{N}$ acts on a trivial bundle of dimension $6$. The six dimensions are labeled by a pair $(I, j)$ where  $I = u, v$ and $j = 0, \pm 1$. That is to say, $E = \bigoplus\limits_{I = u,v, j = 0, \pm 1} l_{I, j}$, where each $l_{I,j}$ is a trivial complex line bundle near the knot. Let $E_j = l_{u, j} \oplus l_{v, j}$, $j = 0, \pm 1$. Then $N$ acts on each $E_j$ separately, denoted as $N_j$ respectively.\\

Suppose $ \begin{pmatrix} \alpha \\ \beta \end{pmatrix} = \begin{pmatrix}
ue^{-i(m+1)\theta - iks} \\ v e^{-im\theta - iks}
\end{pmatrix}$ is a Fourier component of a section of $E_j$, where both $m$ and $\dfrac{2kl}{2\pi}$ are integers. Then for each $j$ the normal operator acting on the Fourier component is given by
$$\tilde{N}_{j}\begin{pmatrix}
ue^{-i(m+1)\theta - iks} \\ v e^{-im\theta - iks}
\end{pmatrix}= \begin{pmatrix}
     (( \partial_{\rho} - \dfrac{(m  + 2\gamma j)}{\rho}) v +iku) e^{-i(m+1)\theta - iks}\\
 ((\partial_{\rho} + \dfrac{(m+1 + 2\gamma j)}{\rho})u - ik v)e^{-im\theta - iks}\end{pmatrix}.$$
For each fixed $j$, the indicial roots are
$$ \lambda_{u, j} = -(m+1 +2\gamma j), ~~~ \lambda_{v, j} = (m + 2\gamma j).$$

 Then locally $N_j$ is a bounded operator from $\rho^{\delta + \frac{1}{2}}H_e(E_j)$ to $\rho^{\delta - \frac{1}{2}}L^2(E_j)$ for any $\delta$, where
locally (for a configuration $\psi$ supported in $N_{\epsilon} \backslash K$), the $H_e$ norm is given by the square root of 

$$\int_{N_{\epsilon} \backslash K}\rho^2 (|\partial_{\rho} \psi|^2 + |\partial_s \psi|^2) +  (|\partial_{\theta}\psi|^2 + |\psi|^2) d\rho d\theta ds.$$

The $L^2$ norm is given by the square root of 

$$\int_{N_{\epsilon}} |\psi|^2 d\rho d\theta ds.$$

One can apply Mazzeo's theory in \cite{Mazzeo1991EllipticI} and easily check that the weights $\delta$ that make $\tilde{L}_{\Psi}$ Fredholm are all numbers such that
$$\min\{- 2\gamma |j| ,  2\gamma |j| - 1\} < \delta < \max\{- 2\gamma |j|, 2\gamma |j| - 1\}.$$

In particular, if $\gamma \neq \dfrac{1}{4}$, then $\delta = \dfrac{1}{2}$ is the weight that we chose to use in section \ref{Section: Introduction and preliminary set-ups}.

\paragraph{Near infinity}~\\

Approaching infinity, $\tilde{L}_{\Psi}$ is just the twisted Dirac operator plus a zero order term. It is standard that the indicial roots are a discrete subset of $\IR$. And the Fredholm weights are all real numbers except the indicial roots. \\

To be more precise, near infinity (say, in $\IR^3 \backslash B_R$) $\tilde{L}_{\Psi}$ is a bounded map from $r^{\delta} H_{\infty}$ to $r^{\delta - 1}L^2$, where the $H_{\infty}$ norm is given by the square root of
$$\int_{\IR^3 \backslash B_R} (r^2 |\nabla \psi|^2 + |\psi|^2) r^2 dr dS.$$
And the $L^2$ norm is just the ordinary $L^2$ norm under the measure $r^2 dr dS$.\\

Note that the $\IH$ defined in section \ref{Section: Introduction and preliminary set-ups} is not the same as any version of $r^{\delta} H_{\infty}$ near infinity. We'll not further study this aspect in this paper.

\paragraph{Asymptotic behavior near the knot}~\\

We assume $\gamma \neq \dfrac{1}{4}$. According to Mazzeo's elliptic edge theory, every element $\psi$ in the kernel of $\tilde{L}_{\Psi}$ in $\IH$ has the following poly-homogeneous expansion near the knot:

$$\psi \sim \sum\limits_{\lambda} \sum\limits_{p = 0}^{n_{\lambda}} \psi_{\lambda, p}\rho^{\lambda}(\log \rho)^p,$$
where the summation is taken over all $\lambda$ such that $$ \lambda = n + 2j\gamma| n \in \IZ, ~ j = 0, \pm 1, ~ \text{and}~ \lambda > -\dfrac{1}{2}.$$ 
Each $n_{\lambda}$ is a non-negative integer. 
All $\psi_{\lambda, p}$ are smooth in $\theta$ and $s$. But we'll give a different proof in this paper which requires zero background knowledge in micro-local analysis. \\

Note that if $\psi_1 \sim \rho^{\lambda_1}$, $\psi_2 \sim \rho^{\lambda_2}$, then the quadratic operator $Q$ sends $(\psi_1, \psi_2)$ to $Q(\psi_1,  \psi_2) \sim \rho^{\lambda_1 + \lambda_2}$. If both $\lambda_1, \lambda_2 > \dfrac{1}{2}$, then $\lambda_1 + \lambda_2 > -1$ and intuitively $\tilde{L}_{\Psi}^{-1} (Q(\psi_1, \psi_2)) \sim \rho^{\lambda_1 + \lambda_2 + 1}$ with $\lambda_1 + \lambda_2 + 1 > 0$. The author suspects that conjecture \ref{the conjecture} may be proved using a similar argument with Mazzeo and Witten's argument for the extended Bogomolny equations in \cite{Mazzeo2013TheCondition}.

\subsection{The norms $\IH_{\delta}^k$ and $\IL_{\delta}$}

This subsection defines some norms that can be viewed as generilizations of $\IH$ and $\IL$ locally on $N_{2\epsilon} \backslash K$. However, the readers should note that these norms are only defined on functions/configurations that are supported near the knot. Here are the definitions:

\begin{definition}
Suppose $\psi$ is a smooth function or a smooth configuration on $\IR^3 \backslash K$ and is supported in $N_{2\epsilon}$. Then we can define the following norms:
\begin{itemize}
    \item The $\IL_{\delta}$ norm of $\psi$ is the square root of
    $$\int_0^{l}\int_0^{2\pi}\int_{0}^{2\epsilon} \rho^{-1 - 2\delta} |\psi|^2 d\rho d\theta ds. $$
    \item The $\IH_{\delta}^{k}$ norm of $\psi$, where $k \geq 1$, is the square root of
    $$\int_0^{l}\int_0^{2\pi}\int_{0}^{2\epsilon} \rho^{ - 1 - 2\delta} (\sum\limits_{j = 0}^k \rho^{2j}|\nabla^j \psi|^2) d\rho d\theta ds.$$
    Note that we have assumed that $\psi$ is supported in $N_{2\epsilon}$. So $\IH^k_{\delta}$ is indeed a norm (not just a semi-norm).
    \item More generally the $\IL^{p}_{\delta}$ norm is the $p$th root of
     $$\int_0^{l}\int_0^{2\pi}\int_{0}^{2\epsilon} \rho^{-1 - p\delta} |\psi|^p d\rho d\theta ds. $$
\end{itemize}
\end{definition}

We also use $\IH_{\delta}$ to represent $\IH^1_{\delta}$. Recall that $|\nabla \psi|^2$ is equivalent with $$|\partial_{\rho}\psi|^2 + \dfrac{1}{\rho^2}|\partial_{\theta}\psi|^2 + |\partial_s \psi|^2,$$ which we will frequently use as a substitution of $|\nabla \psi|^2$. These norms are equivalent with the weighted Sobolev norms defined in Mazzeo's paper \cite{Mazzeo1991EllipticI}, although in our situation they are only defined on functions/configurations that are supported near the knot.\\

We have the following basic properties of these norms: (Recall that $\chi$ is the cut-off function which is $1$ on $N_{\epsilon}$ and is $0$ on $\IR^3 \backslash N_{2\epsilon}$.)\\

 The subscript $\delta$ has the following property which can be verified directly:
\begin{proposition}
    $\chi \rho^{\mu}$ is in $\IL_{\delta}$ (that is to say, has finite $\IL_{\delta}$ norm) if and only if $\mu > \delta$.  
\end{proposition}

The following proposition studies the limit behavior of $\psi$ as $\rho \rightarrow 0$ in these norms.

\begin{proposition}\label{proposition of limit behavior of H and L} We assume that $\psi$ is a configuration supported in $N_{2\epsilon} \backslash K$.
\begin{itemize}
    \item Suppose $\delta < 0$. Suppose $\rho \partial_{\rho} \psi$ has finite $\IL_{\delta}$ norm. (Note that this is a weaker condition than $\psi$ has finite $\IH_{\delta}$ norm.) Then $$\lim\limits_{\rho \rightarrow 0} \rho^{-\delta}|\psi| = 0.$$ Moreover, $\psi$ has finite $\IL_{\delta}$ norm and 
    $$||\rho \partial_{\rho}(\psi)||_{\IL_{\delta}} \geq |\delta| \cdot ||\psi||_{\IL_{\delta}}.$$
    \item Suppose $\psi$ has finite $\IL_{\delta}$ norm. Then
    $$\liminf\limits_{\rho \rightarrow 0} \rho^{-\delta} |\psi| = 0. $$
    Moreover, if we also have $\rho \partial_{\rho}(\psi)$ has finite $\IL_{\delta}$ norm. (We don't need to assume $\delta < 0$ in this case.) Then we actually have
    $$\lim\limits_{\rho \rightarrow 0} \rho^{-\delta}|\psi| = 0.$$
    \item Suppose $\rho \partial_{\rho} \psi$ has finite $\IL_{\delta}$ norm  and suppose $\delta > 0$. Then there is a configuration/function  $\psi_0 \in L^2(d\theta ds)$ which is independent with $\rho$  such that
    $$\lim\limits_{\rho \rightarrow 0} \int_0^l \int_0^{2\pi} |\psi - \psi_0|^2 d\theta ds = 0.$$Moreover, we have
    $$\int_0^l \int_0^{2\pi} |\psi - \psi_0 |^2 d\theta ds = O(\rho^{\delta}). $$
    Note that $\psi_0$ may not be smooth. Also note that when $\psi \in \IL_{\delta}$, $\psi_0 = 0$.
\end{itemize}

\end{proposition}

\begin{proof}
To start, we prove the first bullet.\\

Suppose $$f_{\delta}(\rho) = \int_0^l \int_0^{2\pi} \rho^{-2\delta} |\psi|^2 d\theta ds.$$ Then for any positive constant $C$ and any two values $\rho_1 < \rho_2 < 2\epsilon$, we have
$$f_{\delta}(\rho_1) - f_{\delta}(\rho_2) = -\int_{\rho_1}^{\rho_2} \int_0^l \int_0^{2\pi} \partial_{\rho}(\rho^{-2\delta} |\psi|^2) d\theta ds d\rho $$
$$ = \int_{\rho_1}^{\rho_2} \int_0^l \int_0^{2\pi} (-2\rho^{-2\delta} <\partial_{\rho}\psi, \psi> + 2\delta \rho^{-2\delta - 1} |\psi|^2) d\theta ds d\rho $$
$$\leq \int_{\rho_1}^{\rho_2} \int_0^l \int_0^{2\pi} (C\rho^{1-2\delta}|\partial_{\rho}\psi|^2 + \dfrac{1}{C}\rho^{-2\delta - 1}|\psi|^2 + 2\delta \rho^{-2\delta - 1} |\psi|^2) d\theta ds d\rho.$$
If $\delta < 0$, we may choose $\delta_1$ such that $\delta < \delta_1 < 0$ and $C = \dfrac{1}{|2\delta_1|}$. Then
$$f_{\delta_1}(\rho_1) - f_{\delta_1}(\rho_2) \leq \int_{\rho_1}^{\rho_2} \int_0^l \int_0^{2\pi} C\rho^{1-2\delta_1}|\partial_{\rho}\psi|^2  d\theta ds d\rho \leq (\rho_1)^{2\delta - 2\delta_1} \int_{\rho_1}^{\rho_2} \int_0^l \int_0^{2\pi} C\rho^{1-2\delta}|\partial_{\rho}\psi|^2  d\theta ds d\rho.$$

So
$$f_{\delta}(\rho_1) - (\dfrac{\rho_2}{\rho_1})^{2\delta - 2\delta_1} f_{\delta}(\rho_2) \leq C \int_{\rho_1}^{\rho_2} \int_0^l \int_0^{2\pi} C\rho^{1-2\delta}|\partial_{\rho}\psi|^2  d\theta ds d\rho. $$

Note that $(\dfrac{\rho_2}{\rho_1})^{2\delta - 2\delta_1} < 1$ and
$$\lim\limits_{\rho_1, \rho_2 \rightarrow 0} \int_{\rho_1}^{\rho_2} \int_0^l \int_0^{2\pi} C\rho^{1-2\delta}|\partial_{\rho}\psi|^2  d\theta ds d\rho = 0.$$

We may fix $\dfrac{\rho_2}{\rho_1}$ and let $\rho_1 \rightarrow 0$. Then 
$\limsup\limits_{\rho \rightarrow 0} f_{\delta}(\rho) \geq 0$ exists and satisfies
$$\limsup\limits_{\rho \rightarrow 0} f_{\delta}(\rho) \leq (\dfrac{\rho_2}{\rho_1})^{2\delta - 2\delta_1} \limsup\limits_{\rho \rightarrow 0} f_{\delta}(\rho). $$

This implies that $$\lim\limits_{\rho \rightarrow 0} f_{\delta}(\rho) = 0. $$

Moreover, 
$$0 = \int_{0}^{2\epsilon} \int_0^l \int_0^{2\pi} \partial_{\rho}(\rho^{-2\delta} |\psi|^2) d\theta ds d\rho  = \int_{2\epsilon}^{0} \int_0^l \int_0^{2\pi} (2\rho^{-2\delta} <\partial_{\rho}\psi, \psi> - 2\delta \rho^{-2\delta - 1} |\psi|^2) d\theta ds d\rho $$
$$\geq \int_{\epsilon}^{0} \int_0^l \int_0^{2\pi} (-\dfrac{1}{|\delta|}\rho^{1-2\delta}|\partial_{\rho}\psi|^2 - |\delta|\rho^{-2\delta - 1}|\psi|^2 - 2\delta \rho^{-2\delta - 1} |\psi|^2) d\theta ds d\rho.$$

Because $-2\delta > 0$, $-|\delta| - 2\delta = |\delta| > 0$. So 
$$\int_{\rho_1}^{\rho_2} \int_0^l \int_0^{2\pi} \rho^{-2\delta - 1} |\psi|^2 d\theta ds d\rho \leq \dfrac{1}{|\delta|^2} \int_0^{2\epsilon} \int_0^l \int_0^{2\pi} \rho^{1 - 2\delta} |\partial_{\rho}\psi|^2 d\theta ds d\rho < +\infty, $$
which means, in fact, $\psi \in \IL_{\delta}$ and
$$||\psi||_{\IL_{\delta}} \leq \dfrac{1}{|\delta|} ||\rho (\partial_{\rho}\psi)||_{\IL_{\delta}}. $$

On the other hand, if $\psi \in \IL_{\delta}$, then it is obvious that 
$$\liminf\limits_{\rho \rightarrow 0} \rho^{-\delta} |\psi| = 0.$$
If moreover $\psi \in \IH_{\delta}$, then we may choose $\delta_1 > \delta$ and

$$f_{\delta_1}(\rho_1) - f_{\delta_1}(\rho_2) \leq C \int_{\rho_1}^{\rho_2} \int_0^l \int_0^{2\pi} (\rho^{1-2\delta_1}|\partial_{\rho}\psi|^2  +  \rho^{-2\delta_1 - 1} |\psi|^2) d\theta ds d\rho.$$

Still we have

$$f_{\delta}(\rho_1) - (\dfrac{\rho_2}{\rho_1})^{2\delta - 2\delta_1} f_{\delta}(\rho_2) \leq  C \int_{\rho_1}^{\rho_2} \int_0^l \int_0^{2\pi} (\rho^{1-2\delta}|\partial_{\rho}\psi|^2  +  \rho^{-2\delta - 1} |\psi|^2) d\theta ds d\rho. $$

And from the fact that
$$\lim\limits_{\rho_1, \rho_2 \rightarrow 0} C\int_{\rho_1}^{\rho_2} \int_0^l \int_0^{2\pi} (\rho^{1-2\delta}|\partial_{\rho}\psi|^2 + \rho^{-2\delta - 1}|\psi|^2)  d\theta ds d\rho = 0,$$
we get 
$$\lim\limits_{\rho \rightarrow 0} f_{\delta}(\rho) = 0.$$

Now we prove the third bullet. Suppose $0 < \rho_1 < \rho_2$ are small but $\dfrac{\rho_2}{\rho_1}$ is bounded from above. Suppose $C$ is a large constant which may depend on this bound (and may change from line to line as usual). Then
$$|\psi(\rho_2) - \psi(\rho_1)|^2 \leq  (\int_{\rho_1}^{\rho_2} |\partial_{\rho} \psi| d\rho)^2 \leq \rho_2(\int_{\rho_1}^{\rho_2} |\partial_{\rho} \psi|^2 d\rho) \leq C \int_{\rho_1}^{\rho_2} \rho |\partial_{\rho} \psi|^2 d\rho. $$
So
$$\int_0^l\int_0^{2\pi} \rho_2^{-2\delta}|\psi(\rho_1) - \psi(\rho_2)|^2 d\theta ds \leq C\int_0^l\int_0^{2\pi}\int_{\rho_1}^{\rho_2} \rho^{1-2\delta} |\partial_{\rho} \psi|^2 d\rho d\theta ds \leq C. $$

If $\delta > 0$,  we actually have
$$\int_0^l\int_0^{2\pi} |\psi(\rho_1) - \psi(\rho_2)|^2 d\theta ds \leq C\rho_2^{2\delta }.$$
So
$$(\int_0^l\int_0^{2\pi} |\psi(\rho_1) - \psi(\rho_2)|^2 d\theta ds)^{\frac{1}{2}} \leq C\rho_2^{\delta}.$$
So for any geometric sequence $\{\rho_n\}$ that converges to $0$ (or more generally, any sequence that  converges to $0$ and such that $\dfrac{\rho_{n+1}}{\rho_{n}}$ is bounded above), $\psi(\rho_n)$ is a Cauchy sequence in $L^2(d\theta ds)$. So there exists an element $\psi_0$ as a function/configuration in $\theta$ and $s$ (that is, independent with $\rho$) and in $L^2(d\theta ds)$  such that
$$\lim\limits_{\rho \rightarrow 0} \int_0^l \int_0^{2\pi} |\psi - \psi_0|^2 d\theta ds = 0.$$

In fact, we have
$$(\int_0^l\int_0^{2\pi} |\psi(\rho_1) - \psi(\rho_2)|^2 d\theta ds)^{\frac{1}{2}} \leq C\rho_2^{\delta}.$$
So for any geometric sequence $\{\rho_n\}$ that converges to $0$, we have
$$(\int_0^l \int_0^{2\pi} |\psi(\rho_1) - \psi(\rho_n)|^2 d\theta ds)^{\frac{1}{2}} \leq C \sum\limits_{i = 1}^n \rho_{i}^{\delta} \leq C \rho_1^{\delta}. $$
Here the last $C$ above depends on the common ratio of the sequence, but doesn't depend on $n$. We may fix the common ratio to be $\dfrac{1}{2}$. Then this implies that
$$(\int_0^l \int_0^{2\pi} |\psi(\rho) - \psi_0|^2 d\theta ds)^{\frac{1}{2}} \leq  C \rho^{\delta}. $$

Note that $\psi_0$ is not in $\IL_{\delta}$ for any $\delta \geq 0$. So if $\psi \in \IL_{\delta}$, then $\psi_0 = 0$.

\end{proof}

\begin{proposition}\label{proposition of generilized weighted Sobolev embedding} Suppose $\mu$ is a constant such that $0 \leq \mu \leq 1$. There exists a constant $C$ which depends on $\delta$ and $N_{2\epsilon}$, such that,
 the $\IL_{\delta - \frac{1}{3}(1 - \mu)}^6$ norm of $\psi$ is bounded above by $C (||\psi||_{\IH_{\delta}} + ||\partial_s \psi||_{\IL_{\delta + \mu - 1}})$. More precisely,
$$ (\int_0^l \int_0^{2\pi} \int_0^{2\epsilon} \rho^{1 - 2\mu} |\rho^{-\delta} \psi|^6 d\rho d\theta ds)^{\frac{1}{6}} \leq C(||\psi||_{\IH_{\delta}} + ||\partial_s \psi||_{\IL_{\delta + \mu - 1}}). $$
An immediate corollary is, the $\IL_{\delta- \frac{1}{4}(1 - \mu)}^4$ norm  of $\psi$ is bounded above by $C (||\psi||_{\IH_{\delta }} + ||\partial_s \psi||_{\IL_{\delta + \mu - 1}})$. More precisely,
$$ (\int_0^l \int_0^{2\pi} \int_0^{2\epsilon} \rho^{- 4\delta  - \mu} | \psi|^4 d\rho d\theta ds)^{\frac{1}{4}} \leq C(||\psi||_{\IH_{\delta}} + ||\partial_s \psi||_{\IL_{\delta + \mu - 1}}). $$
\end{proposition}
Note that when $\mu = 0$, then $||\partial_s \psi||_{\IL_{\delta - 1}}$ can be bounded above by $C ||\psi||_{\IH_{\delta}}$. So lemma \ref{weight version of Sobolev embedding lemma} can be viewed as a special case of this proposition such that $\mu = 0, \delta = -\dfrac{1}{2}$.
\begin{proof} 

Suppose $A(\theta, s)$ is a function that doesn't depend on $\rho$, $B(s, \rho)$ is a function that doesn't depend on $\theta$, $C(\theta, \rho)$ is a function that doesn't depend on $s$. Then
\begin{equation*}
    \begin{array}{ll}
         ~ & ~~~~   \displaystyle \int_0^l \int_0^{2\pi} \int_0^{2\epsilon}A(\theta, s)B(s, \rho) C(\theta, \rho) d\rho d\theta ds\\
         & \leq \displaystyle \int_0^l \int_0^{2\pi} A(\theta, s)(\int_0^{2\epsilon}B(s, \rho)^2d\rho)^{\frac{1}{2}} (\int_0^{2\epsilon} C(\theta, \rho)^2  d\rho)^{\frac{1}{2}} d\theta ds\\
         & \leq \displaystyle \int_0^{2\pi} (\int_0^l  A(\theta, s)^2 ds)^{\frac{1}{2}}(\int_0^l\int_0^{2\epsilon}B(s, \rho)^2d\rho ds)^{\frac{1}{2}} (\int_0^{2\epsilon} C(\theta, \rho)^2  d\rho)^{\frac{1}{2}} d\theta \\
         & \leq \displaystyle  (\int_0^l \int_0^{2\pi} A(\theta, s)^2 d\theta ds)^{\frac{1}{2}}(\int_0^l\int_0^{2\epsilon}B(s, \rho)^2d\rho ds)^{\frac{1}{2}} (\int_0^{2\pi} \int_0^{2\epsilon} C(\theta, \rho)^2  d\rho d\theta)^{\frac{1}{2}}
    \end{array}
\end{equation*}
We assume
\begin{equation*}
\begin{array}{ll}
     \rho^{-4\delta- \mu} |\psi|^4 &  \leq  \displaystyle  \int_0^{2\pi}\rho^{-4\delta - \mu }\partial_{\theta}(|\psi|^4) d\theta  + \dfrac{1}{2\pi}\int_0^{2\pi} \rho^{-4\delta - \mu} |\psi|^4 d\theta\\
          & \displaystyle \leq C(\int_0^{2\pi} (\rho^{-1 - 2\delta}|\partial_{\theta} \psi|^2 + \rho^{-1 - 2\delta}|\psi|^2) d{\theta})^{\frac{1}{2}} ( \int_0^{2\pi} \rho^{1 - 2\mu - 6\delta} |\psi|^6 d{\theta})^{\frac{1}{2}} = B(s, \rho)^2.
\end{array}
\end{equation*}
Similarly,
\begin{equation*}
           \rho^{1 - 4\delta - 2\mu}|\psi|^4  \leq  C \displaystyle (\int_0^{l} (\rho^{1 - 2\mu - 2\delta}|\partial_{s} \psi|^2 + \rho^{-1 - 2\delta} |\psi|^2) ds)^{\frac{1}{2}} ( \int_0^{l} \rho^{1 - 6\delta - 2\mu} |\psi|^6  ds)^{\frac{1}{2}} = C(\theta, \rho)^2.
\end{equation*}
In the above inequality, we used the fact that $\rho \leq 2\epsilon$ is bounded.
\begin{equation*}
    \begin{array}{ll}
          \rho^{1 - 4\delta - \mu} |\psi|^4  & \leq \displaystyle \int_0^{2\epsilon} (|\partial_{\rho} (\rho^{1 - 4\delta - 2\mu}|\psi|^4) | + \dfrac{1}{2\epsilon}\rho^{1 - 4\delta - 2\mu}|\psi|^4) d{\rho}\\
          & \leq \displaystyle C (\int_0^{2\epsilon} \rho^{1 - 2\delta}|\partial_{\rho}\psi|^2 + \rho^{-1 - 2\delta}|\psi|^2d\rho )^{\frac{1}{2}}(\int_0^{2\epsilon} \rho^{1 - 6\delta - 2\mu}|\psi|^6 d\rho)^{\frac{1}{2}} = A(\theta, s)^2
    \end{array}
\end{equation*}
So
\begin{equation*}
    \begin{array}{ll}
         & \displaystyle ~~~~ \int_0^l\int_0^{2\pi} \int_0^{2\epsilon} \rho^{1 - 6\delta - 2\mu}|\psi|^6 d\rho d\theta ds   \leq \displaystyle \int_0^l \int_0^{2\pi} \int_0^{2\epsilon}A(\theta, s)B(s, \rho) C(\theta, \rho) d\rho d\theta ds\\
         & \leq \displaystyle  (\int_0^l \int_0^{2\pi} A(\theta, s)^2 d\theta ds)^{\frac{1}{2}}(\int_0^l\int_0^{2\epsilon}B(s, \rho)^2d\rho ds)^{\frac{1}{2}} (\int_0^{2\pi} \int_0^{2\epsilon} C(\theta, \rho)^2  d\rho d\theta)^{\frac{1}{2}}\\    \end{array}
\end{equation*}
         $$\leq C \displaystyle (\int_0^l\int_0^{2\pi}\int_0^{2\epsilon} \rho^{1 - 2\delta}(|\nabla \psi|^2 + \rho^{-2\mu}|\partial_s \psi|^2 + \dfrac{1}{\rho^2}|\psi|^2) d\rho d\theta ds)^{\frac{3}{4}}(\int_0^l\int_0^{2\pi}\int_0^{2\epsilon} \rho^{1 - 2\mu - 6\delta}|\psi|^6 d\rho d\theta ds)^{\frac{3}{4}}.$$

This implies 
$$     (\int_0^l\int_0^{2\pi} \int_0^{2\epsilon} \rho^{1 - 2\mu - 6\delta}|\psi|^6 d\rho d\theta ds)^{\frac{1}{6}}   \leq C  (\int_0^l\int_0^{2\pi}\int_0^{2\epsilon} \rho^{1 - 2\delta}(|\nabla \psi|^2 + \rho^{-2\mu}|\partial_s \psi|^2 + \dfrac{1}{\rho^2}|\psi|^2) d\rho d\theta ds)^{\frac{1}{2}}.$$
\end{proof}

Here is an immediate corollary

\begin{corollary}
Suppose $Q(\psi_1, \psi_2)$ is a bi-linear map with smooth and bounded coefficients. Then $Q$ maps $\IH_{\delta_1} \times \IH_{\delta_2}$ into $\IL_{\delta_1 + \delta_2 - \frac{1}{2}}$.\\

Moreover, suppose $\mu_1 $ and $\mu_2$ are constants in $[0, 1]$. Then
$$||Q(\psi_1, \psi_2)||_{\IL_{\delta_1 + \delta_2 - \frac{1}{2}(2 - \mu_1 - \mu_2) }} \leq C (||\psi_1||_{\IH_{\delta_1}} + ||\partial_s \psi_1||_{\IL_{\delta_1 - 1 + \eta_1}}) \cdot (||\psi_2||_{\IH_{\delta_2}} + ||\partial_s \psi_2||_{\IL_{\delta_2 - 1 + \eta_2}}). $$
\end{corollary}

\subsection{Properties of the normal operator $\tilde{N}$}

Recall that each Fourier component of $\psi$ can be written as 
$$\psi_m = \begin{pmatrix}
 \alpha \\ \beta
\end{pmatrix} = \begin{pmatrix}
 u(\rho) e^{-i(m + 1)\theta - iks} \\  v(\rho) e^{-im \theta - iks}
\end{pmatrix},$$
where $m$ is an integer and $k$ is a real number. The normal operator $\tilde{N}$ is written as

$$\tilde{N}\begin{pmatrix}
\alpha \\ \beta
\end{pmatrix}= \begin{pmatrix}
 - (\partial_s \alpha)  + (\partial_{\rho} - \dfrac{i\partial_{\theta} }{\rho} - \dfrac{(2\gamma j - 1)}{\rho}) \beta  \\
      (\partial_{\rho} + \dfrac{i\partial_{\theta}}{\rho} + \dfrac{( 2\gamma j + 1)}{\rho})\alpha + (\partial_s \beta) 
\end{pmatrix}.$$
Or more precisely, acting on $\begin{pmatrix}u \\ v\end{pmatrix}$ it is convenient to write $\tilde{N}$ as the following (if there is no ambiguity)

$$\tilde{N} \begin{pmatrix} u \\ v \end{pmatrix} = \begin{pmatrix} 
(\partial_{\rho} + \dfrac{(1 - \lambda)}{\rho}) v + ik u\\
(\partial_{\rho} + \dfrac{\lambda }{\rho})u - ik v \\ \end{pmatrix},$$
where $\lambda =  m + 2\gamma j + 1$.\\

If either $\lambda + 1 + \delta$ or $2 - \lambda + \delta$ equals $0$, $\delta + 1$ is called the ``indicial weight". This term seems to be originated from  the theory of regular singular ODE. It is also used in micro-local analysis for edge type of operators. See for example \cite{Mazzeo1991EllipticI}.

\begin{proposition}\label{proposition of generilized local Fredholmness}
Suppose $\delta \geq \delta' \geq -\dfrac{3}{2}$. Suppose there are no indicial weights in the closed interval $[\delta' + 1, \delta + 1]$. (That is to say, for any $\eta \in [\delta' + 1, \delta + 1]$, both $\lambda + 1 + \eta$ and $2 - \lambda + \eta$ take values in a discrete closed subset of $\IR$ which doesn't include $0$.)  Suppose $\psi$ has finite $\IL_{\delta' + 1}$ norm and suppose $\tilde{N}(\psi)$ has finite $\IL_{\delta}$ norm. Then
$$C ||\tilde{N}(\psi)||_{\IL_{\delta}} \geq  ||\psi||_{\IH_{\delta + 1}}. $$
\end{proposition}

Note: this proposition can be viewed as a generalized version of proposition \ref{proposition of local Fredholmness}. In fact, proposition \ref{proposition of local Fredholmness} corresponds to the case that $-\dfrac{1}{2}$ is not an indicial weight (since $\gamma \neq \dfrac{1}{4}$) and $\delta = \delta' = -\dfrac{3}{2}$. (Strictly speaking, we used $\tilde{N}$ instead of $N$ here. But they are essentially equivalent.)

\begin{proof}
The proof has three steps.
\paragraph{Step 1:} We prove a special case when $\delta' = \delta$.\\

We only need to prove the inequality on each Fourier component. So we may assume $\psi$ has only one Fourier component.\\

Note that the fact that $\psi$ has finite $\IL_{\delta}$ norm implies that $$\liminf\limits_{\rho \rightarrow 0}\rho^{-2 - 2\delta}(|u|^2 + |v|^2) = 0.$$

This limit behavior guarantees that in the following arguments, the boundary terms all vanish in integration by parts. Using integration by parts, we have

$$ \int_0^{\epsilon}\rho^{-1 - 2\delta}|\tilde{N}(\psi)|^2 d\rho   = \int_0^{\epsilon} \rho^{-1 - 2\delta}(|\partial_{\rho} u + \dfrac{\lambda}{\rho} u + ik v|^2 + |\partial_{\rho} v + \dfrac{(1 - \lambda)}{\rho} v - ik u|^2 )d\rho $$
$$ = \int_0^{\epsilon} \rho^{-1 - 2\delta} (|\partial_{\rho}u|^2 + |\partial_{\rho}v|^2 + k^2(|u|^2 + |v|^2) + \dfrac{\lambda^2}{\rho^2}|u|^2 + \dfrac{(1 - \lambda)^2 }{\rho^2}|v|^2)     d\rho $$
$$ +  \int_0^{\epsilon} (2 + 2\delta)\rho^{-3 - 2\delta} (\lambda |u|^2 + (1 - \lambda)|v|^2) - (4 + 4 \delta)\rho^{-2 - 2\delta}<v, ik u> )     d\rho. $$

Since  both $\lambda + 1 + \delta$ and $2 - \lambda + \delta$ take values in a discrete closed subset of $\IR$ which doesn't include $0$. We only need to prove that there is a large enough constant $C$ uniformly in $m$ and $k$ such that all of the three inequalities are satisfied:

$$ C \int_0^{\epsilon}\rho^{-1 - 2\delta}|\tilde{N}(\psi)|^2 d\rho   \geq \int_0^{\epsilon} \rho^{-3 - 2\delta} ((\lambda + 1 + \delta)^2 |u|^2 + (2 - \lambda + \delta)^2 |v|^2)     d\rho ;$$
$$ C \int_0^{\epsilon}\rho^{-1 - 2\delta}|\tilde{N}(\psi)|^2 d\rho   \geq k^2 \int_0^{\epsilon} \rho^{-1 - 2\delta} ( | u|^2 + |v|^2)     d\rho ;$$
$$ C\int_0^{\epsilon}\rho^{-1 - 2\delta}|\tilde{N}(\psi)|^2 d\rho  \geq  \int_0^{\epsilon} \rho^{-1 - 2\delta} ( |\partial_{\rho} u|^2 + |\partial_{\rho} v|^2)     d\rho.$$

For each real number $t$ (to be determined),

$$0 \leq \int_0^{2\epsilon} \rho^{-1 - 2\delta}(|\partial_{\rho} u - \dfrac{(1 + \delta)}{\rho} u + itk v|^2 + |\partial_{\rho} v - \dfrac{(1 + \delta)}{\rho} v - itk u|^2)d\rho = $$
$$ \int_0^{2\epsilon} \rho^{-1 - 2\delta} (|\partial_{\rho} u|^2 + |\partial_{\rho}v|^2  +  t^2 k^2( |u|^2 + |v|^2))  + 2 t\rho^{-2 - 2\delta}  <v, ik u> d\rho $$ $$ -(1 + \delta)^2 \int_0^{2\epsilon}\rho^{-3 - 2\delta} (|u|^2  + |v|^2) d\rho.$$

So for each real number $t$,

$$ \int_0^{\epsilon}\rho^{-1 - 2\delta}|\tilde{N}(\psi)|^2 d\rho  $$
$$ \geq \int_0^{\epsilon} \rho^{-1 - 2\delta} ((1 - t^2)k^2( | u|^2 + | v|^2) + \dfrac{(\lambda + 1 + \delta)^2}{\rho^2}|u|^2 + \dfrac{(2 - \lambda  + \delta)^2 }{\rho^2}|v|^2)     d\rho $$
$$ -  \int_0^{\epsilon}\rho^{-1 - 2\delta}(4 + 4 \delta + 2t)\rho^{-2 - 2\delta}<v, ik u>      d\rho. $$

\paragraph{Case 1:} Suppose $-\dfrac{3}{2} \leq \delta \leq -\dfrac{1}{2}$. Then we may take $t = -2(1 + \delta)$. There is a constant $C$ which doesn't depend on $m$ or $k$ such that

$$ C \int_0^{2\epsilon}\rho^{-1 - 2\delta}|\tilde{N}(\psi)|^2 d\rho  \geq  \int_0^{2\epsilon} \rho^{-3 - 2\delta} ((\lambda + 1 + \delta)^2 |u|^2 + (2 - \lambda + \delta)^2 |v|^2)     d\rho.$$

Then we take $t = 0$ and get that for a possibly larger $C$,

$$ C \int_0^{2\epsilon}\rho^{-1 - 2\delta}|\tilde{N}(\psi)|^2 d\rho  \geq  k^2 \int_0^{2\epsilon} \rho^{-1 - 2\delta} ( | u|^2 + | v|^2)     d\rho.$$

Finally all these inequalities imply that 
$$ C \int_0^{2\epsilon}\rho^{-1 - 2\delta}|\tilde{N}(\psi)|^2 d\rho \geq \int_0^{2\epsilon} \rho^{-1 - 2\delta} ( |\partial_{\rho} u|^2 + |\partial_{\rho} v|^2)     d\rho.$$

\paragraph{Case 2:} Suppose $\delta > - \dfrac{1}{2}$. Then we take $t = 0$ and get

$$  \int_0^{2\epsilon}\rho^{-1 - 2\delta}|\tilde{N}(\psi)|^2 d\rho  $$
$$ \geq  \int_0^{2\epsilon} \rho^{-1 - 2\delta}( k^2( |u|^2 + | v|^2) + \dfrac{(\lambda + 1 + \delta)^2}{\rho^2}|u|^2 + \dfrac{(2 - \lambda  + \delta)^2 }{\rho^2}|v|^2)     d\rho$$
$$ - \int_0^{2\epsilon}\rho^{-1 - 2\delta}(4 + 4 \delta)\rho^{-2 - 2\delta}<v, ik u>     d\rho . $$

For any positive real number $c$, we have
$$|\int_0^{2\epsilon}\rho^{-2 - 2\delta}<v, ik u>  d\rho| \leq \int_0^{2\epsilon}\rho^{-1 - 2\delta}(c|v|^2 + \dfrac{k^2}{4c \rho^2} |u |^2) d\rho.$$

Note that when $\delta > -\dfrac{1}{2}$, we always have

$$2|\lambda + 1 + \delta| + 2|2 - \lambda + \delta| \geq 2|2 \delta + 3| > |4 + 4\delta|.  $$

So there is a constant $C$ such that
$$ C \int_0^{2\epsilon}\rho^{-1 - 2\delta}|\tilde{N}(\psi)|^2 d\rho $$
$$ \geq \ \int_0^{2\epsilon} \rho^{-1 - 2\delta} (k^2( | u|^2 + | v|^2) + \dfrac{(\lambda + 1 + \delta)^2}{\rho^2}|u|^2 + \dfrac{(2 - \lambda  + \delta)^2 }{\rho^2}|v|^2)     d\rho .$$

Finally this implies that for a possibly different $C$,
$$ C \int_0^{2\epsilon}\rho^{-1 - 2\delta}|\tilde{N}(\psi)|^2 d\rho \geq \int_0^{2\epsilon} \rho^{-1 - 2\delta} ( |\partial_{\rho} u|^2 + |\partial_{\rho} v|^2)     d\rho.$$

\paragraph{Step 2:} This step proves a slightly more general case than the previous step. We assume $\delta - \delta'$ is small enough, whose precise meaning is: For each valid $\delta$, there is a $\delta' < \delta$ depending on $\delta$ such that $\delta - \delta'$ or any smaller number is small enough. On the other hand, for each valid $\delta'$, there is a $\delta > \delta'$ which depends on $\delta'$ such that $\delta - \delta'$ or any smaller number is small enough.\\

From the previous step, we know $||\psi||_{\IH_{\delta'}}$ is finite. In fact, we can modify the proof of the previous step in the new situation. The new difficulty is that we do not have the following limit behavior anymore
$$\liminf\limits_{\rho \rightarrow 0} \rho^{-2 - 2\delta} (|u|^2 + |v|^2) = 0.$$
So the boundary terms from $\rho \rightarrow 0$ from integration by parts cannot be ignored. Instead, our new assumption is: For slightly smaller $ \delta' < \delta$,
$$ \int_0^{2\epsilon} \rho^{-1 - 2(\delta' + 1)}((k^2\rho^2 + m^2 + 1)(|u|^2 + |v|^2) + \rho^2(|\partial_{\rho}u|^2 + |\partial_{\rho}v|^2)) d\rho < +\infty. $$
The above inequality is uniform for different Fourier components (that is, for different $k$, $m$). So let $$F(\rho) = \int_{\rho}^{2\epsilon}  \rho^{-1 - 2(\delta + 1)} ((k^2\rho^2 + m^2 + 1)(|u|^2 + |v|^2) + \rho^2(|\partial_{\rho}u|^2 + |\partial_{\rho}v|^2) d\rho.$$ We have 
$$ \lim\limits_{\rho \rightarrow 0} \rho^{\delta - \delta'} F(\rho) = 0.$$

The above limit is uniform for different Fourier components. \\

In the arguments in the previous step, if we don't ignore the boundary terms in the integration by parts and if we take the integration by parts from $\rho$ to $2\epsilon$. Then the boundary term is bounded above by $C \rho |F'(\rho)|$, where the constant $C$ doesn't depend on $m, k$. So taking the boundary term into account, the inequality should be
$$F(\rho) \leq  C \rho |F'(\rho)| + C||\tilde{N}(\psi)||^2_{\IL_{\delta}}. $$
Note that $F(\rho)$ is a decreasing function in $\rho$. So in fact $|F'(\rho)| = - F'(\rho)$. And 
$$(\rho^{\frac{1}{C}}F(\rho))' \leq C^2 \rho^{\frac{1}{C} - 1} ||\tilde{N}(\psi)||_{\IL_{\delta}}^2.  $$
Since $\delta - \delta'$ is small enough, we may assume $\delta - \delta' < \dfrac{1}{C}$. Then
$$\lim\limits_{\rho \rightarrow 0} \rho^{\frac{1}{C}}F(\rho) = 0. $$

So
$$\rho^{\frac{1}{C}}F(\rho) \leq C^2\int_0^{\rho} \rho^{\frac{1}{C} - 1}||\tilde{N}(\psi)||^2_{\IL_{\delta}} = C^3 \rho^{\frac{1}{C}} ||\tilde{N}(\psi)||^2_{\IL_{\delta}}. $$

Thus 
$$F(\rho) \leq  C^3 ||\tilde{N}(\psi)||^2_{\IL_{\delta}}. $$
Let $\rho \rightarrow 0$. Note that $\lim\limits_{\rho \rightarrow 0} F(\rho)$
is equivalent with the $\IH_{\delta + 1}$ norm of $\psi$. We get (for a possibly larger $C$)
$$||\psi||_{\IH_{\delta + 1}} \leq C ||\tilde{N}(\psi)||_{\IL_{\delta}}.$$

\paragraph{Step 3:} Finally we prove the proposition without any additional assumptions. The previous step implies that all the weights $\eta \in [\delta', \delta]$ such that
$$C ||\tilde{N}(\psi)||_{\IL_{\eta}} \geq  ||\psi||_{\IH_{\eta + 1}} $$
is an relatively open and closed subset of $[\delta', \delta]$ which is non-empty. Thus this subset is $[\delta', \delta]$ itself. In particular, we may choose $\eta = \delta$ and get
$$C ||\tilde{N}(\psi)||_{\IL_{\delta}} \geq  ||\psi||_{\IH_{\delta + 1}}. $$

\end{proof}

The reader should be aware of the subtle point in the above proposition: We have assumed that $\psi$ has finite  $\IL_{\delta' + 1}$ norm in the first place. Here is an counterexample when $\psi$ is not assumed to be in $\IL_{\delta' + 1}$: We may find a solution to $\tilde{N}(\psi) = 0$ near the knot using Bessel's function which has an infinite  $\IL_{\delta + 1}$ norm. And use a cut-off function to make it supported in $N_{2\epsilon}$. Then clearly we cannot use the $\IL_{\delta}$ norm of $\tilde{N}(\psi)$ (which is finite) to bound its $\IH_{\delta + 1}$ norm (which is infinite).\\

The following proposition studies the case where either  $\lambda + 1 + \delta$ or $2 - \lambda + \delta$ equals $0$. (In another word to say, $\delta + 1$ is an indicial weight.) Moreover, $m$ is determined from $\delta$. So we may assume that there is only one Fourier component with respect to $\theta$. Or more precisely, $\psi$ can be written as $\psi_0 e^{-im\theta}$ for some $\psi_0$, $m$, where $\psi_0$ is independent with $\theta$. And here it is:

\begin{proposition}\label{proposition of generilized local Fredholmness indicial version}
Suppose either $\lambda + 1 + \delta = 0$ or $2 - \lambda + \delta = 0$. Suppose $\delta > -\dfrac{3}{2}$. We may assume that $\psi$ can be written as $\psi_0 e^{-im\theta}$, where $m$ is determined from $\delta$. Suppose $\psi$ is supported in $N_{\epsilon}$.  Suppose $\psi$ has finite $\IH_{\delta' + 1}$ norm, where $\delta'$ is slightly less than $\delta + \mu$. Suppose $0 \leq \mu < 1$. Then there is a constant $C$ such that
$$C ||\tilde{N}(\psi)||_{\IL_{\delta + \mu}} \geq  ||\partial_{\rho}(\rho^{-(\delta + 1)} \psi)||_{\IL_{\mu}} + ||\partial_s\psi||_{\IL_{\mu + \delta}}. $$
\end{proposition}

\begin{proof}
The idea is similar with the proof of proposition \ref{proposition of generilized local Fredholmness}. Without loss of generality, we may assume $\lambda + 1 + \delta = 0$. First of all, if we can ignore the boundary terms and let $\tilde{u} = \rho^{-(\delta + 1)}u, ~~ \tilde{v} = \rho^{-(\delta + 1)} v$, then from integration by parts, we have

$$ \int_0^{\epsilon}\rho^{-1 - 2\delta- 2\mu}|\tilde{N}(\psi)|^2 d\rho   = \int_0^{\epsilon} \rho^{-1 - 2\delta - 2\mu}(|\partial_{\rho} u - \dfrac{(\delta + 1)}{\rho} u + ik v|^2 + |\partial_{\rho} v + \dfrac{(2 + \delta)}{\rho} v - ik u|^2 )d\rho $$
$$ = \int_0^{\epsilon} \rho (|\partial_{\rho}\tilde{u} + ik\tilde{v}|^2 + |\partial_{\rho}\tilde{v} + \dfrac{(3 + 2\delta)}{\rho}\tilde{v} - ik\tilde{u}|^2)     d\rho $$
$$ =  \int_0^{\epsilon} \rho^{1 - 2\mu}(|\partial_{\rho}\tilde{u}|^2 + |\partial_{\rho}\tilde{v}|^2 + k^2(|\tilde{u}|^2 + |\tilde{v}|^2) +  \dfrac{2(2 + 2\delta  + 2\mu)}{\rho}<\tilde{u}, ik\tilde{v}> + \dfrac{(3 + 2\delta + 2\mu)(3 + 2\delta)}{\rho^2}|\tilde{v}|^2)    d\rho. $$

For each $-1 \leq t \leq 1$. We have
$$0 \leq \int_0^{\epsilon}\rho^{1 - 2\mu} (|\partial_{\rho}\tilde{u} + itk\tilde{v}|^2 + |\partial_{\rho}\tilde{v} - itk \tilde{u}|^2) d\rho$$ $$=  \int_0^{\epsilon} \rho^{1 - 2\mu}(|\partial_{\rho}\tilde{u}|^2 + |\partial_{\rho}\tilde{v}|^2 + t^2k^2(|\tilde{u}|^2 + |\tilde{v}|^2) + \dfrac{(4\mu - 2) t}{\rho}<\tilde{u}, itk\tilde{v}>)   d\rho.$$

Using the above inequality, we get, for any $|t| \leq 1$,

$$\int_0^{\epsilon}\rho^{-1 - 2\delta- 2\mu}|\tilde{N}(\psi)|^2 d\rho  \geq \int_0^{\epsilon} \rho^{1 - 2\mu}((1 - t^2)k^2(|\tilde{u}|^2 + |\tilde{v}|^2) + \dfrac{|(4\mu - 2)t|}{\rho}|<\tilde{u}, ik\tilde{v}>|) d\rho$$ $$+ \int_0^{\epsilon}\rho^{1 - 2\mu}( \dfrac{2(2 + 2\delta  + 2\mu)}{\rho}<\tilde{u}, ik\tilde{v}> + \dfrac{(3 + 2\delta + 2\mu)(3 + 2\delta)}{\rho^2}|\tilde{v}|^2)    d\rho$$
$$\geq (2\sqrt{(1 - t^2)(3 + 2\delta + 2\mu)(3 + 2\delta)} + |(4\mu - 2)t| - 2(2 + 2\delta + 2\mu) ) \int_0^{\epsilon} \rho^{- 2\mu}k|\tilde{u}||\tilde{v}| d\rho. $$

We may choose some $|t| \leq 1$ such that
$$(2\sqrt{(1 - t^2)(3 + 2\delta + 2\mu)(3 + 2\delta)} + |(4\mu - 2)t| = \sqrt{4(3 + 2\delta + 2\mu)(3 + 2\delta) + (4\mu - 2)^2}. $$

When $\delta > -\frac{3}{2}$, $0 \leq \mu < 1$, it is not hard to prove
$$ \sqrt{4(3 + 2\delta + 2\mu)(3 + 2\delta) + (4\mu - 2)^2} > 2|2 + 2\delta + 2\mu|. $$

This implies that there is a constant $C$ which doesn't depend on $k$ such that

$$\int_0^{\epsilon}\rho^{-1 - 2\delta- 2\mu}|\tilde{N}(\psi)|^2 d\rho \geq \dfrac{1}{C} \int_0^{\epsilon} \rho^{ - 1 - 2\mu} |\tilde{v}|^2 d\rho, $$
which further implies for a possibly larger $C$,
$$ C\int_0^{\epsilon}\rho^{-1 - 2\delta- 2\mu}|\tilde{N}(\psi)|^2 d\rho \geq  \int_0^{\epsilon} \rho^{ 1 - 2\mu} (|\partial_{\rho}\tilde{u}|^2 + |\partial_{\rho}\tilde{v}|^2 + k^2(|\tilde{u}|^2 + |\tilde{v}|^2) + \dfrac{1}{\rho^2}|\tilde{v}|^2) d\rho.  $$

So we have
$$C ||\tilde{N}(\psi)||_{\IL_{\delta + \mu}} \geq ||\rho\partial_{\rho}(\rho^{-(\delta + 1)} \psi)||_{\IL_{\mu}} + ||\rho\partial_s(\rho^{-(\delta + 1)}\psi)||_{\IL_{\mu}} = ||\rho \partial_{\rho}(\rho^{-(\delta + 1)} \psi)||_{\IL_{\mu}} + ||\partial_s\psi||_{\IL_{\mu + \delta}}. $$

To deal with the boundary term, we can use the same method as in the proof of proposition \ref{proposition of generilized local Fredholmness} by setting 
$$F(\rho) =  \int_{\rho}^{\epsilon} \rho^{ 1 - 2\mu} (|\partial_{\rho}\tilde{u}|^2 + |\partial_{\rho}\tilde{v}|^2 + k^2(|\tilde{u}|^2 + |\tilde{v}|^2) + \dfrac{1}{\rho^2}|\tilde{v}|^2) d\rho. $$
And the boundary term doesn't affect the result.

\end{proof}

\subsection{The proof of proposition \ref{proposition of polyhomogeneous expansion}}
Suppose $\Psi$ is a solution to the Bogomolny equations on $\IR^3 \backslash K$ such that $$\psi = \Psi - \Psi_{\gamma, M, k} \in \IH.$$

We use $\tilde{L}$ to represent $\tilde{L}_{\gamma, M, k}$. Then in a neighbourhood of the knot except the knot itself (which may be assumed to be $N_{2\epsilon} \backslash K$), $\Psi_{\gamma, M, k}$ is an actual solution to the Bogomolny equations. So the Bogomolny equations for $\Psi$ can be written as

$$\tilde{L}(\psi) + Q(\psi, \psi) = 0.$$

We may further write $\tilde{L}$ as $$\tilde{L} = \tilde{N} + \tilde{R}, $$
where $\tilde{N}$ is the normal operator that we have studied in the previous subsection. And $\tilde{R}$ is an bounded operator from $\IH_{\delta}$ to $\IL_{\delta}$ for any $\delta$. \\

Start from $\psi \in \IH$. We know $ \chi \psi \in \IH_{-\frac{1}{2}}$. Since we only care about the behavior of $\psi$ near the knot. We use the notation $\psi \in \IH_{-\frac{1}{2}}$ to represent that $\chi \psi \in \IH_{-\frac{1}{2}}$ for convenience even if $\psi$ is not supported in $N_{2\epsilon} \backslash K$.\\

We always assume $\gamma \neq \dfrac{1}{2}$.

\begin{proposition}
For small enough $\delta > 0$, we have $\psi \in \IH_{-\frac{1}{2} + \delta}$. That is to say, $\chi \psi$ has finite $\IH_{-\frac{1}{2} + \delta}$ norm.
\end{proposition}

\begin{proof}

Let $$F(\rho) = \int_0^l \int_0^{2\pi} \int_{\rho}^{2\epsilon} \rho^{2 - 2\delta}|\nabla \psi|^2  + \rho^{- 2\delta} |\psi|^2 d\rho d\theta ds.$$
Then we just need to show that $\lim\limits_{\rho \rightarrow 0} F(\rho) < +\infty$.\\

We have $\lim\limits_{\rho \rightarrow 0} \rho^{2\delta}F(\rho) = 0$. We assume $C$ is a constant that doesn't depend on $\rho$ and $\psi$. Using the same integration by part as in the proof of proposition \ref{proposition of generilized local Fredholmness}, we know
$$F(\rho) \leq C \int_0^l \int_0^{2\pi}\int_{\rho}^{2\epsilon}  \rho^{2- 2\delta} |\tilde{N}(\psi)|^2 d\rho d\theta ds + C \rho F'(\rho) + C $$
$$\leq C \int_0^l \int_0^{2\pi} \int_{\rho}^{2\epsilon}  \rho^{2- 2\delta} (|\tilde{R}(\psi)|^2 + |Q(\psi, \psi)|^2) d\rho d\theta ds + C \rho F'(\rho) + C.$$

The last constant $C$ above comes from the fact that $\psi$ may not be supported in $N_{2\epsilon} \backslash K$, so there is another bounded boundary term on the end $\rho = 2\epsilon$.\\

Since $R(\psi)$ is in $\IL_{-\frac{1}{2}}$,
$$ \int_0^l \int_0^{2\pi} \int_{\rho}^{2\epsilon}  \rho^{2- 2\delta} |\tilde{R}(\psi)|^2  d\rho d\theta ds  < C. $$ By proposition \ref{proposition of generilized weighted Sobolev embedding},
$$\int_0^l \int_0^{2\pi} \int_{\rho}^{2\epsilon}  \rho^{2- 2\delta} |Q(\psi, \psi)|^2 d\rho d\theta ds \leq CF(\rho) \cdot \int_0^l \int_0^{2\pi} \int_{\rho}^{2\epsilon} (\rho^{2 }|\nabla \psi|^2  + |\psi|^2) d\rho d\theta ds.$$

Since $\psi \in \IH_{-\frac{1}{2}}$, we may choose $\epsilon$ small enough such that

$$C\cdot \int_0^l \int_0^{2\pi} \int_{0}^{2\epsilon} (\rho^{2 }|\nabla \psi|^2  + |\psi|^2) d\rho d\theta ds < \dfrac{1}{2}. $$

Thus what we get is
$$F(\rho) \leq C\rho F'(\rho) + C. $$

Since $\delta$ is small, using the same argument as in the proof of proposition \ref{proposition of generilized local Fredholmness}, we get
$$\lim\limits_{\rho \rightarrow 0} F(\rho) < +\infty.$$
\end{proof}

Note that the above proposition doesn't give an a prior bound on $||\psi||_{-\frac{1}{2} + \delta}$, although proving it to be finite.\\

Now we can start our boot-strapping and get
\begin{proposition}
Suppose $\alpha = \min\{2\gamma, 1 - 2\gamma \}$. Then for any $\delta' > 0$, we have $\psi \in \IH_{-\alpha - \delta'}$. Moreover, $\rho \partial_{\rho} (\psi) \in \IL_{-\alpha + 1 - \delta'}$. 
\end{proposition}

\begin{proof}
From the previous proposition, we can start from $\psi \in \IH_{-\frac{1}{2} + \delta}$ for some small positive $\delta$. Clearly $\delta < 1$. So $-\dfrac{1}{2}+ \delta > -\dfrac{3}{2} + 2\delta$. We get

$$\tilde{N}( \psi) = -(\tilde{R}(\psi) + Q(\psi, \psi)) \in \IL_{ -\frac{3}{2} + 2\delta}.$$
If $-\dfrac{3}{2} + 2\delta < -\alpha$, then by proposition \ref{proposition of generilized local Fredholmness}, we have
$$\psi \in \IH_{-\frac{1}{2} + 2\delta}. $$
And we can replace $\delta$ by $2\delta$ and continue the same step. After finitely many iterations, we get  $-\dfrac{3}{2} + 2\delta \geq - \alpha$. We may use a smaller $\delta$ such that $-\dfrac{3}{2} + 2\delta$ is slightly smaller than $-\alpha$ and do the iteration one more time. Then still by proposition \ref{proposition of generilized local Fredholmness}, we have for any $\delta' > 0$, $\psi \in \IH_{-\alpha - \delta'}$. Moreover, $\rho \partial_{\rho} (\psi) \in \IL_{-\alpha + 1 - \delta'}$.
\end{proof}
Immediately from proposition \ref{proposition of limit behavior of H and L}, we have
\begin{corollary}
There is a function $\psi_{-\alpha}$ in $L^2(d\theta ds)$ which doesn't depend on $\rho$ such that for any $\delta > 0$,
$$\psi - \psi_{-\alpha}\rho^{-\alpha} \in \IL_{-\alpha + 1 - \delta}. $$
\end{corollary}

In fact, we can say something more about $\psi_{-\alpha}$. Recall that in the proof of proposition \ref{proposition of generilized local Fredholmness indicial version}, $\lambda + 1 + \alpha = 0$ and $m$ is fixed. So in fact, $\psi_{-\alpha}$ can be written as $\psi_{-\alpha} = a(s) e^{-im\theta}$ for some integer $m$ on each sub-bundle $l_{I, j}$ of $E$. Different sub-bundle may give different $m$ and $a(s)$. But in any case, $a(s)$ is an $L^2(ds)$ function of $s$. In particular, $\psi_{-\alpha}$ is smooth in the $\theta$ direction. (All its $\theta$ direvatives are in $L^2$.)\\

Note that this finishes the proof of the second bullet of proposition \ref{proposition of polyhomogeneous expansion}. Now we prove the first bullet.\\

\begin{proposition}\label{proposition of smooth of a(s)}
Suppose $0 \leq \alpha < \dfrac{1}{6}$. Then the function $\psi_{-\alpha}$ is also smooth in $s$ direction.
\end{proposition}
\begin{proof}
Start from $\psi \in \IH_{-\alpha - \delta}$ for any $\delta > 0$. We know that
$$\tilde{N}(\psi) = - \tilde{R}(\psi) - Q(\psi, \psi) \in \IL_{-2\alpha - 2\delta - \frac{1}{2}}.$$
So $\tilde{N}(\psi) \in \IL_{-2\alpha - \frac{1}{2} - \delta}$ for any $\delta > 0$. Using proposition \ref{proposition of generilized local Fredholmness indicial version} and the fact that $\alpha + \dfrac{1}{2} < 1$, we know that for small enough $\mu > 0$, $\partial_s \psi \in \IL_{-\alpha - 1 + \mu}$.\\

The next step is an iteration. For any $0 \leq \mu \leq 1.$ Suppose we have
$$\partial_s \psi \in \IL_{-\alpha - 1 + \mu}.$$
Then proposition \ref{proposition of generilized weighted Sobolev embedding} implies that
$$\tilde{N}(\psi) = - \tilde{R}(\psi) - Q(\psi, \psi) \in \IL_{-2\alpha - \delta - \frac{1}{2}(1 - \mu)} $$
for any $\delta > 0$, which further implies that
$$\partial_s \psi \in \IL_{-2\alpha - \delta -\frac{1}{2}(1 - \mu)}.$$
Note that $-2\alpha - \delta - \dfrac{1}{2}(1 - \mu) = -\alpha - 1 + (\dfrac{1}{2} - \alpha) + \dfrac{1}{2}\mu - \delta$. And we may use $\dfrac{1}{2} - \alpha + \dfrac{1}{2}\mu - \delta$ as a new definition of $\mu$ and do the same step. After finitely many iterations, we know that for any $\delta > 0$, $$\partial_s \psi \in \IL_{-3 \alpha - \delta}.$$
Since we have assumed $\alpha < \dfrac{1}{6}$. So $\partial_s \psi \in \IL_{\frac{1}{2}}$.\\

Let $\xi = \partial_s \psi$. Note that $\tilde{N}$ commutes with $\partial_s$. So applying $ \partial_s$ to the identity
$$\tilde{N}(\psi) + \tilde{R}(\psi) + Q(\psi, \psi) = 0,$$
we get
$$\tilde{N}(\xi)  + \tilde{R}(\xi) + [\partial_s, \tilde{R}](\psi) + 2Q(\xi, \psi) = 0. $$

Using the same arguments as in the proof of proposition \ref{proposition of generilized local Fredholmness indicial version}, we know that there is a configuration $\xi_{-\alpha}$ which doesn't depend on $\rho$, smooth in $\theta$ and $L^2$ in $s$ such that
$$\xi - \xi_{-\alpha}\rho^{-\alpha} \in \IL_{-\alpha + \delta} $$
for some small enough $\delta > 0$.\\

Taking the integration in $s$ direction implies that $\xi_{-\alpha} = \partial_s \psi_{-\alpha}$. So in fact $\psi_{-\alpha}$ is in $L^{1,2}$ in the $s$ direction.\\

We can continue the same arguments and get that $\psi_{-\alpha}$ is in $L^{m, 2}$ for any positive integer $m$. Thus $\psi_{-\alpha}$ is fully smooth in $s$ direction.
\end{proof}

\bibliography{references.bib}
\bibliographystyle{ieeetr}

\end{document}